%% file: main.tex
\pdfoutput=1
\documentclass[11pt]{article}
\usepackage{graphicx} %
\usepackage{tikz-cd}
\usetikzlibrary{arrows.meta,positioning,calc,decorations.pathreplacing,fit}
\usepackage{subcaption}

\usepackage[margin=1in]{geometry}
\usepackage{amsmath}
\usepackage{amsfonts}
\usepackage{amssymb}
\usepackage{amsthm}

\usepackage{xcolor}
\usepackage{xurl} %
\providecommand{\href}[2]{#2} %

\usepackage[sort]{cite}

\title{Averting `Gravity': Launch Control and the Limits of Collision Avoidance in Low Earth Orbit Operations}
\author{Aditya S. Gopalan \and Harsha Honnappa}
\date{}

\newcommand{\kfrags}{\mathrm{kfrags}}
\newcommand{\Qsat}{Q_\mathrm{sat}}

\newcommand{\crit}{\mathrm{crit}}

\input{math}

\begin{document}

\maketitle

\begin{abstract}
As large constellations fill low Earth orbit, the risk of a runaway collision cascade has become a central concern of space operations.
Debris population models treat collision avoidance as a fixed input, although the capacity to screen potential collisions and plan maneuvers degrades as debris accumulates.
Using a stylized multi-layer queueing network model, we show that collision avoidance capabilities that saturate as the in-custody fragment count increases, can make the orbital shell bistable, with a low-debris and a high-debris equilibrium, and hysteretic. This has two major operational consequences.
First, leaving the high-debris equilibrium requires a large cut in launch rates. Second, in the low-debris equilibrium, fragment counts respond with a delay as the launch rate rises past a tipping point.~
Removing bistability is desirable, and we show that fragment count-indexed admission control can achieve this.~
The least restrictive such policy lengthens the lag, whereas an exponential tilt of this policy makes the fragment count track its equilibrium with a bounded lag.
Furthermore, admitting more launches can reduce the rate of completed missions, and under a unimodality condition a threshold rule on the debris fragment count maximizes this rate at equilibrium.

\end{abstract}

\input{intro}

\section{Model and Fluid Limit}\label{sec:model}
We consider a layered queueing network, where the first layer is a processor sharing queue with deadlines and the second layer is a queueing network which captures the evolution of intact satellite, derelict, and debris fragment counts in Low Earth Orbit.
Intact satellites can perform collision avoidance maneuvers and the processor sharing queue captures conjunction avoidance.

\input{avoidance_layer}

\subsection{Orbital Shell Queueing Network}
We consider a three station queueing network model.
Station~1 models the orbital shell as an Erlang-B $M/M/n/n$ loss system with $n$ orbital slots, Poisson launch arrivals at rate $\lambda$, and exponential mission durations with mean duration $1/\mu_I$. Launches finding all $n$ slots occupied are blocked: the quantity $n$ is the \emph{capacity} of the orbital shell.
Stations~2 and~3 model the derelict and debris fragment populations, respectively, and are infinite-server queues with exponential service at rates $\mu_D$ and $\mu_F$, respectively.
For stations~2 and~3, there is no blocking and each derelict or fragment immediately finds a free server.
Thus, the number of intact satellites is constrained, but there is no limit to the number of derelicts or fragments that can be in a particular shell.
Customers in service at stations~1 and~2 may depart sooner than their exponential service distribution: this is a state-dependent dynamics and represents catastrophic collisions before planned mission lifetime.

\paragraph{Routing Among the Queues.}
Routing among stations is state-dependent, so the process is Markov on the state space $\{(Q_I, Q_D, Q_F) \in \bZ_+^3\}$, which represents the counts of intacts, derelicts, and fragments, respectively.

Before we exactly specify the routing, we describe the interpretation of routing events.
Satellites can end their lifetime in three ways:
\begin{itemize}
    \item Catastrophic collision with a fragment: produces a new batch of fragments.
    \item Post-mission disposal (PMD) success: departure from the system.
    \item PMD failure: produces a derelict.
\end{itemize}

Derelicts can end their lifetime in two ways:
\begin{itemize}
    \item Catastrophic collision with a fragment: produces a new batch of fragments.
    \item De-orbit due to drag: departure from the system.
\end{itemize}

Fragments de-orbit only due to drag.

We now specify the routing explicitly.
PMD success occurs with probability $P$, independently for each satellite that does not undergo catastrophic collision.
Denote by $\kfrags$ the mean number of fragments produced by a catastrophic collision of a satellite or a derelict.
We assume that each catastrophic collision produces a new batch of fragments which is distributed as $\Poi(\kfrags)$, independently of everything else.
Intact satellites are able to perform collision avoidance, as discussed in the previous section.
Notice that the cumulative instantaneous rate at which potential conjunctions arrive to \emph{each} intact is $cQ_F(t)$.

For a random variable $X$, we write $Y = X + \Poi(k)$ if $Y$ is distributed as the sum of $X$ and a $\Poi(k)$ random variable which is independent of everything else.

Satellites are routed as follows:
\begin{itemize}
    \item At rate $\alpha(Q_F(t))cQ_F(t)$, a catastrophic collision occurs and $Q_F(t) \sim Q_F(t-) + \Poi(\kfrags)$.
    \item At rate $(1-P)\mu_I$, the satellite becomes a derelict and $Q_D(t) = Q_D(t-) + 1$.
    \item Otherwise, the satellite departs the system at rate $P\mu_I$.
\end{itemize}

Derelicts are routed as follows:
\begin{itemize}
    \item At rate $cQ_F(t)$, a catastrophic collision occurs and $Q_F(t) \sim Q_F(t-) + \Poi(\kfrags)$.
    \item Otherwise, the derelict departs the system at rate $\mu_D$.
\end{itemize}

All customers depart the system after completing service in the third station at rate $\mu_F$.

\begin{figure}
    \centering
    \includegraphics[width=0.75\linewidth]{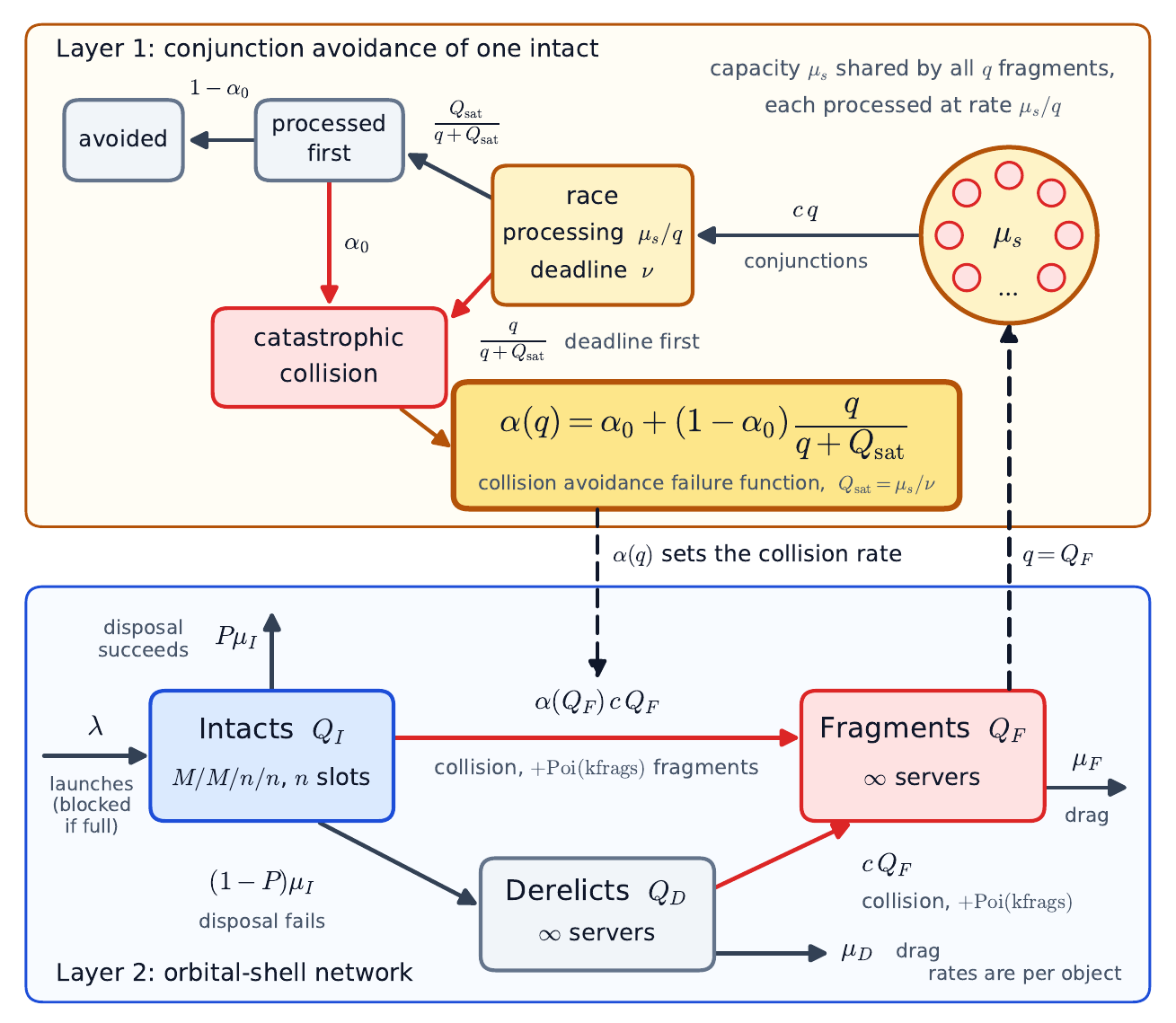}
    \caption{Pictorial depiction of the layered queueing model.}
    \label{fig:layered-model}
\end{figure}

The full model is depicted visually in Figure~\ref{fig:layered-model}.

We refer the reader to Mandelbaum, Massey, and Reiman~\cite{mandelbaum1998b} for a careful treatment of Markov queueing networks with state-dependent routing.

\paragraph{Coupling Between Queueing Network Layers. }
Let $(Q_I, Q_D, Q_F) \in \bZ_+^3$ denote the intact, derelict, and fragment counts in the second layer.
Let $Q_T$ denote the number of fragments in custody by a typical intact in the first layer, that is,
$Q_T := \frac{1}{Q_I}\sum_{i=1}^{Q_I}C_i,$
where $C_i$ is the number of fragments in custody for the $i$-th intact.
Our model assumes that $C_i = Q_F$ for all $i$.
It is immediate that
$Q_T(t) = \frac{1}{Q_I(t)}\sum_{i=1}^{Q_I(t)}C_i(t) = Q_F(t)$ whenever $Q_I(t) > 0$.
Thus, the joint process for both layers $(Q_T, Q_I, Q_D, Q_F)$ is exactly determined by the state $(Q_I, Q_D, Q_F)$ of the second layer.

\subsection{Fluid Limit}
Exact analysis of the Markov chain is intractable for constellation-scale $n$.
Instead, we approximate the large-$n$ model with its limit as $n\to\infty$.
We consider a mean-field \emph{fluid limit} which captures the aggregate behavior of the system and takes the form of ordinary differential equations.
To do this, we scale the constant $\Qsat$ to depend on the constellation size $n$ as follows, so that the relative ratios are invariant to scaling by $n$:
\[\Qsat^n = n\Qsat = \frac{n\mu_s}{\nu},\]
with
\[\alpha^n(q) := \alpha_0 + (1-\alpha_0)\frac{q}{q + \Qsat^n}.\]
In addition, assume that in the $n$-th system, launches arrive at rate $n\lambda$, and the conjunctions occur at rate $c/n$, so that the rate $\alpha^n(Q_F)\tfrac{c}{n}Q_F$ at which an intact suffers a catastrophic collision depends on the fragment count only through the density $Q_F/n$.
All other parameters are held fixed, and $\mathcal{S} := [0,1]\times[0,\infty)^2$ denotes the state space of the fluid-scaled process.

The limiting equations follow from a strong-approximation coupling of the Markov chain to a Brownian system,
stated as Theorem~\ref{thm:strong-approx} and proved in Section~\ref{sec:limiting}.

\begin{theorem}
    \label{thm:fluid-limit}
    Suppose $\lambda \leq \mu_I$.
    Let $Q: = (Q_I, Q_D, Q_F)$, and let the superscript $n$ denote the number of servers for intacts.
    Suppose that $\frac{1}{n}Q^n(0) \to q(0) \in \mathcal{S}$.
    Then there exists a probability space upon which, almost surely, the following convergence holds, uniformly on compacts:
    \[\frac{1}{n}Q^n(t) \to q(t),\]
    where $q(t) := (q_I(t), q_D(t), q_F(t))$ is given by:
    \begin{align}
        \frac{dq_I}{dt} &:= \lambda - \mu_I q_I - \alpha(q_F)cq_Iq_F
        \label{eq:FL-1}
    \\
        \frac{dq_D}{dt} &:= (1-P)\mu_I q_I - \mu_D q_D - cq_Dq_F
        \label{eq:FL-2}
    \\
        \frac{dq_F}{dt} &:= \kfrags\ \alpha(q_F) c q_F q_I + \kfrags c q_F q_D - \mu_F q_F.
        \label{eq:FL-3}
        \end{align}
\end{theorem}

The hypothesis $\lambda \leq \mu_I$ is a sufficient condition for this form of the fluid limit but is not necessary.
In this regime, the system is below its critical threshold and the $M/M/n/n$ loss dynamics do not play a role.
We use the assumption $\lambda \leq \mu_I$ to demonstrate that even the subcritical system possesses interesting behaviors with challenging policy implications.

The fixed point of the fluid limits are solutions to the following equations:
\begin{align}
    \lambda - \mu_I q^*_I - \alpha(q^*_F)cq^*_Iq^*_F &= 0,
    \label{eq:FL-FP-1}
    \\
    (1-P)\mu_I q^*_I - \mu_Dq^*_D - cq^*_Dq^*_F & = 0,
    \label{eq:FL-FP-2}
    \\
    \kfrags\ \alpha(q^*_F)cq^*_Iq^*_F + \kfrags\ cq^*_Dq^*_F - \mu_F q^*_F &= 0.
    \label{eq:FL-FP-3}
\end{align}

Notice that there may be multiple fixed points for a given set of parameters $\lambda, c, \mu_I, \mu_D, \mu_F, \kfrags, \Qsat$.
The following is the \emph{trivial fixed point}:
\begin{align*}
    q^*_I = \frac{\lambda}{\mu_I},
    \\
    q^*_D = \frac{\lambda}{\mu_D}(1-P),
    \\
    q^*_F = 0.
\end{align*}
We will refer to any fixed point in which $q^*_F > 0$ as a \emph{non-trivial fixed point}.

\subsection{Remarks on Non-Saturating Collision Avoidance}
For this section only, let $\alpha(q) := \alpha_0$, which is the predominant assumption in the literature.
We treat the equilibrium counts of intacts, derelicts, and fragments as functions of the launch rate $\lambda$.
Notice that if $\alpha(q) = \alpha_0$, then the parameter $\Qsat$ is irrelevant.
In this setting, there is a transcritical bifurcation between the trivial branch and a branch with a positive fragment count which occurs at the supremum value of $\lambda$ for which the trivial fixed point is the only fixed point.

It is easy to see by inspection that in a non-trivial fixed point, we have
\begin{align*}
    q^*_I &= \frac{\lambda}{\mu_I + \alpha_0 c q^*_F},
    \\
    q^*_D &= \frac{(1-P)\lambda \mu_I}{(\mu_I + \alpha_0 c q^*_F)(\mu_D + cq^*_F)},
    \\
    0 &= \kfrags\ \alpha_0 c q^*_I + \kfrags\ c q^*_D - \mu_F.
\end{align*}

The equilibrium fragment count $q^*_F$ solves
\[\mu_F(\mu_I + \alpha_0 c q^*_F)(\mu_D + c q^*_F) = \kfrags\ c\lambda [\alpha_0(\mu_D + c q^*_F) + (1-P)\mu_I].\]

To understand the bifurcation structure of the system with non-saturating collision avoidance, it suffices to understand the bifurcation structure of $q^*_F$ as a function of $\lambda$.
\begin{proposition}
    The system with $\alpha(q) = \alpha_0$ undergoes a transcritical bifurcation at
    \[\lambda_\crit = \frac{\mu_I\mu_D\mu_F}{\kfrags\ c (\alpha_0 \mu_D + (1-P)\mu_I)}\]
    if $\lambda_\crit < \mu_I$.
    If $\lambda_\crit < \mu_I$, then for $\mu_I \geq \lambda > \lambda_\crit$ there is a unique non-trivial stable branch on which $q^*_F$ and $q^*_I$ are increasing in $\lambda$, and $q^*_D$ is decreasing in $\lambda$.
\end{proposition}
\begin{proof}
    If $q^*_F \geq 0$, notice that $\lambda$ is an increasing function of $q^*_F$, which can be seen by direct evaluation.
    Thus, for any $\lambda > \lambda_\crit$, there is at most one positive root for $q^*_F$, and the $q^*_F$ is increasing in $\lambda$.
    If $\lambda < \lambda_\crit$, then there are no positive solutions for $q^*_F$.
    Thus, for $\lambda > \lambda_\crit$, there is a positive branch for $q^*_F$.
    The trivial branch is stable if and only if $\lambda < \lambda_\crit$, and the positive branch is stable: both facts are special cases of the stability analysis in Section~\ref{sec:proof-bifurcation} (the trivial fixed point is treated in the Proof of Theorem~\ref{thm:bifurcation-structure}, and for the positive branch apply Lemma~\ref{lem:linearization} with $\alpha' \equiv 0$, so that $e = w = 0$ and $H' < 0$).

    For the monotonicity of $q^*_I$ and $q^*_D$, the last of the fixed-point equations above, together with $q^*_D = (1-P)\mu_I q^*_I/(\mu_D + cq^*_F)$, gives
    \[
        q^*_I = \frac{\mu_F}{\kfrags\ c\left(\alpha_0 + \frac{(1-P)\mu_I}{\mu_D + cq^*_F}\right)},
        \qquad
        q^*_D = \frac{(1-P)\mu_I\,\mu_F}{\kfrags\ c\big(\alpha_0(\mu_D + cq^*_F) + (1-P)\mu_I\big)}.
    \]
    Since $P < 1$, the first is strictly increasing and the second strictly decreasing in $q^*_F$, and $q^*_F$ is strictly increasing in $\lambda$.
    In particular, along this branch more launches raise the equilibrium intact count, even though they also raise the equilibrium debris.
\end{proof}

The key operational insights come from the fact that $q^*_F$ is a strictly increasing function of $\lambda$, whenever $\lambda > \lambda_\crit$ (recall that $q^*_F = 0$ when $\lambda \leq \lambda_\crit$).
Indeed, if the equilibrium fragment count slightly exceeds a given threshold, the system can be restored to an equilibrium with a smaller fragment count simply by slightly lowering the launch rate.
Note, however, that the equilibrium intact count is increasing in $\lambda$ along this branch, so lowering the launch rate to reduce the debris also lowers the equilibrium intact count.

The linearized system only sees the dynamics under a sufficiently small launch rate and cannot capture the full system dynamics.
\begin{proposition}
    For $\alpha(q) = \alpha_0$, the linearized dynamics have a stable fixed point iff $\lambda < \lambda_\crit$.
\end{proposition}
\begin{proof}
    Linearizing around the trivial equilibrium point gives the following dynamics.
    Denote by $\Delta q_I := q_I - q^*_I, \Delta q_D := q_D - q^*_D, \Delta q_F := q_F - q^*_F$.
    The linearized fluid limit at the trivial fixed point is given by:
    \begin{align*}
        \frac{d \Delta q_I}{dt} & := -\mu_I \Delta q_I - \frac{\alpha_0 c \lambda}{\mu_I} \Delta q_F,
        \\
        \frac{d \Delta q_D}{dt} & := \mu_I(1 - P)\Delta q_I - \mu_D \Delta q_D - \frac{c(1-P)\lambda}{\mu_D}\Delta q_F,
        \\
        \frac{d \Delta q_F}{dt} & := \left(\kfrags\ c\,\lambda\left(\frac{\alpha_0}{\mu_I} + \frac{1-P}{\mu_D}\right) - \mu_F\right)\Delta q_F = \mu_F\left(\frac{\lambda}{\lambda_\crit} - 1\right)\Delta q_F.
    \end{align*}
    The terms in $\Delta q_I$ and $\Delta q_D$ are absent from the last equation because the collision rates are proportional to $q_F$, which vanishes at the trivial fixed point.
    The Jacobian is therefore block upper triangular, with eigenvalues $s_I = -\mu_I$, $s_D = -\mu_D$, and $s_F = \mu_F(\lambda/\lambda_\crit - 1)$, which are all negative if and only if $\lambda < \lambda_\crit$. This establishes the result.
\end{proof}

\begin{proposition}
    For $\alpha(q) = \alpha_0$, the linearized system exhibits a critical slowing down as $\lambda \uparrow \lambda_\crit$.
\end{proposition}
\begin{proof}
    It suffices to show that the restoring rate approaches zero at a linear rate in $\lambda_\crit - \lambda$.
    By the proof of the previous proposition, the Jacobian at the trivial fixed point is block upper triangular, with eigenvalues
    \[
        s_I = -\mu_I, \qquad s_D = -\mu_D, \qquad s_F = \mu_F\left(\frac{\lambda}{\lambda_\crit} - 1\right) = -\mu_F\left(1 - \frac{\lambda}{\lambda_\crit}\right).
    \]
    For $\lambda$ close to $\lambda_\crit$, $s_F$ is the eigenvalue closest to zero (the dominant eigenvalue), so the restoring rate is
    \[
        |s_F| = \frac{\mu_F}{\lambda_\crit}\,(\lambda_\crit - \lambda),
    \]
    which is linear in $\lambda_\crit - \lambda$.
    At $\lambda = \lambda_\crit$ the dominant eigenvalue vanishes, while $s_I$ and $s_D$ remain negative.
    Equivalently, $\det(\cL) = s_Is_Ds_F = -\mu_I\mu_D\mu_F(1 - \lambda/\lambda_\crit)$.
\end{proof}

It follows that when $\lambda < \lambda_\crit$, the restoring rate is positive and that any fluctuation in the number of fragments will vanish at an exponential rate.
However, as $\lambda \uparrow \lambda_\crit$, the linearized system takes an increasing time to equilibrate from small perturbations, before ultimately becoming unstable.

\section{Bistability, Hysteresis, and Slow Passage}\label{sec:bistability}

In this section, we analyze the fixed points of the fluid limit equations~\eqref{eq:FL-1},~\eqref{eq:FL-2},~\eqref{eq:FL-3} as functions of the launch rate $\lambda$.
We identify the regimes in terms of $\Qsat$, assuming the remaining parameters are fixed, with $\alpha_0, P \in (0,1)$.
When $\Qsat$ is large, there is a transcritical bifurcation and the system behavior is as above with the non-saturating collision avoidance.
When $\Qsat$ is small, there are saddle-node bifurcations, with bistability and hysteresis.
There are two qualitatively different regimes which produce bistability and hysteresis.
In the first, there is a supercritical transcritical bifurcation, and the jump between the bistable branches occurs when the lower branch is positive.
In the second, the jump to the upper branch occurs from the trivial branch itself.
In addition, when fragments decay quickly relative to intacts and derelicts, a fixed point which the fold structure alone would predict to be stable can lose its stability through a Hopf bifurcation, and the system then settles on a periodic orbit of recurring debris outbreaks instead of a fixed point (Section~\ref{sec:hopf}).
Finally, when the launch rate increases slowly through the upper fold, the jump from the low-debris to the high-debris branch is delayed, and occurs only after the launch rate has overshot the fold by an amount of order $\varepsilon^{2/3}$, where $\varepsilon$ is the rate of increase (Section~\ref{sec:slow-passage-result}).

Solving~\eqref{eq:FL-FP-1} and~\eqref{eq:FL-FP-2} for $q^*_I$ and $q^*_D$ in terms of $q^*_F$ gives
\begin{equation}\label{eq:fp-intact-derelict}
    q^*_I = \frac{\lambda}{\mu_I + \alpha(q^*_F) c q^*_F},
    \qquad
    q^*_D = \frac{(1-P)\mu_I}{\mu_D + c q^*_F}\, q^*_I,
\end{equation}
and substituting into~\eqref{eq:FL-FP-3} shows that the non-trivial fixed points are exactly the solutions $q^*_F > 0$ of
\[\kfrags\ \lambda c H(q^*_F) = \mu_F,\]
where
\begin{equation}\label{eq:H}
    H(q) := \frac{\alpha(q)}{\mu_I + \alpha(q) c q} + \frac{(1-P)\mu_I}{(\mu_I + \alpha(q) c q)(\mu_D + cq)}.
\end{equation}
The quantity $\kfrags\ \lambda c H(q_F)$ is the rate at which a single fragment generates new fragments when the intact and derelict populations are at their equilibrium values given $q_F$: the first term of $H$ is the contribution of collisions with intacts, and the second is the contribution of collisions with derelicts.
Equivalently, the non-trivial fixed points are the graph of
\begin{equation}\label{eq:cM}
    \lambda = \cM(q_F) := \frac{\mu_F}{c\,\kfrags\,H(q_F)}, \qquad q_F > 0,
\end{equation}
with $\cM(0) = \lambda_\crit$.
A \emph{fold} is a positive root of $H'$ at which $H'$, and hence $\cM'$, changes sign, \textit{i.e.}, a local extremum of $\cM$; every fold is a critical point of $\cM$, and every simple positive root of $H'$ is a fold.

\paragraph{Standing assumption on the launch rate.}
Station~1 is a loss system, so the fluid-scaled intact count is confined to $q_I \leq 1$, and~\eqref{eq:FL-1} describes the intacts only as long as no launches are blocked.
Throughout this section we assume $\lambda \leq \mu_I,$
which is the hypothesis of Theorem~\ref{thm:fluid-limit}.
Under this assumption, there is no blocking at the fluid scale, by Proposition~\ref{prop:reflected-fluid}(iii).
If instead $\lambda > \mu_I$, then the fluid limit is the reflected system~\eqref{eq:RFL} of Proposition~\ref{prop:reflected-fluid} and not~\eqref{eq:FL-1}--\eqref{eq:FL-3}, and the results of this section do not apply; see Remark~\ref{rem:capacity}.

Accordingly, all launch rates in this section are in $(0,\mu_I]$, every statement of the form ``for every $\lambda$ in a set $\Lambda$'' is to be read as ``for every $\lambda \in \Lambda\cap(0,\mu_I]$'', and stability is relative to $\mathcal{S}$.
The trivial fixed point has $q^*_I = \lambda/\mu_I \leq 1$, and by~\eqref{eq:fp-intact-derelict} a non-trivial fixed point has $q^*_I < \lambda/\mu_I \leq 1$, so it lies strictly below the wall.
In terms of the fragment count, there is a non-trivial fixed point with $q^*_F = z$ at an admissible launch rate if and only if $\cM(z) \leq \mu_I$, that is, if and only if
\begin{equation}\label{eq:admissible}
    \kfrags\, c\,\mu_I H(z) \geq \mu_F.
\end{equation}
In particular, a bifurcation of the fixed points is seen as the launch rate is increased only if it occurs below $\mu_I$: the transcritical bifurcation at $\lambda_\crit = \cM(0)$ requires $\lambda_\crit < \mu_I$, that is,
\[
    \kfrags\, c\,\big(\alpha_0\mu_D + (1-P)\mu_I\big) > \mu_D\mu_F,
\]
and a fold $q_\dagger$ requires $\cM(q_\dagger) < \mu_I$.
If $\lambda_\crit > \mu_I$, then the trivial fixed point is stable for every admissible launch rate.

Define the thresholds
\[
    \overline{Q} := \frac{(1-\alpha_0)\mu_I}{c\,\alpha_0^2},
    \qquad
    \underline{Q} := \frac{(1-\alpha_0)\mu_I\mu_D^2}{c\left(\alpha_0^2\mu_D^2 + (1-P)\mu_I^2 + (1-P)\alpha_0\mu_I\mu_D\right)},
\]
and the constant
\[
    \varepsilon_0 := \sup_{q > 0}\frac{q\,\alpha'(q)}{\alpha(q)} = \frac{1-\sqrt{\alpha_0}}{1+\sqrt{\alpha_0}} \in (0,1).
\]

\begin{theorem}
    \label{thm:bifurcation-structure}
    We have $0 < \underline{Q} < \overline{Q}$.
    Assume that $\lambda \leq \mu_I$, and read every range of launch rates below as intersected with $(0,\mu_I]$.
    The trivial fixed point is locally asymptotically stable if $\lambda < \lambda_\crit$ and unstable if $\lambda > \lambda_\crit$.
    The non-trivial fixed points of the fluid limit~\eqref{eq:FL-1},~\eqref{eq:FL-2},~\eqref{eq:FL-3} exhibit the following:
    \begin{enumerate}
        \item If $\Qsat > \overline{Q}$, then $\cM$ is a strictly increasing bijection from $(0,\infty)$ onto $(\lambda_\crit,\infty)$: there is a unique non-trivial fixed point for every $\lambda > \lambda_\crit$ and none for $\lambda \leq \lambda_\crit$, and $q^*_F$ is strictly increasing in $\lambda$. There are no folds.
        \item If $\Qsat > \underline{Q}$, there is a supercritical transcritical bifurcation at $\lambda_\crit$: the non-trivial branch emerges from the trivial one into $\lambda > \lambda_\crit$.
        There are either no folds or exactly two folds.
        \item If \[2c\,\underline{Q}\,(\beta_0 - \gamma_0 c\alpha_0\underline{Q}) > 0, \qquad \gamma_0 := \alpha_0 \mu_D + (1-P)\mu_I, \qquad \beta_0 := (\mu_I+\mu_D)(P\mu_I-\alpha_0 \mu_D) - \mu_I^2, \] then there exists $\delta > 0$ such that for $\Qsat \in (\underline{Q}, \underline{Q} + \delta)$ there are exactly two folds $q_{\min} < q_{\max}$. With $\lambda_+ := \cM(q_{\min})$ and $\lambda_- := \cM(q_{\max}) < \lambda_+$, there are three non-trivial fixed points for every $\lambda \in (\max\{\lambda_-,\lambda_\crit\}, \lambda_+)$: a low-debris one with $q^*_F < q_{\min}$, a middle one with $q_{\min} < q^*_F < q_{\max}$, and a high-debris one with $q^*_F > q_{\max}$. The middle one is unstable, and under condition~\eqref{eq:rate-condition} below the other two are stable; the trivial fixed point is unstable for these launch rates.
        \item If $0 < \Qsat < \underline{Q}$, the transcritical bifurcation is subcritical: the non-trivial branch emerges into $\lambda < \lambda_\crit$.
        There is exactly one fold $q_{\max}$, and with $\lambda_- := \cM(q_{\max}) < \lambda_\crit$ there are two non-trivial fixed points for every $\lambda \in (\lambda_-, \lambda_\crit)$: a middle one with $q^*_F < q_{\max}$, and a high-debris one with $q^*_F > q_{\max}$. The middle one is unstable, and under condition~\eqref{eq:rate-condition} below the high-debris one is stable; the trivial fixed point is stable for these launch rates. For $\lambda \geq \lambda_\crit$ there is one non-trivial fixed point, the high-debris one, and for $\lambda < \lambda_-$ there is none.
    \end{enumerate}
    The stability of a non-trivial fixed point $q^*$ is as follows:
    \begin{enumerate}
        \item[(a)] If $H'(q^*_F) > 0$, then $q^*$ is unstable. In particular, the middle branch between two folds, and the branch which emerges from a subcritical transcritical bifurcation, are unstable.
        \item[(b)] If $H'(q^*_F) < 0$ and
        \begin{equation}\label{eq:rate-condition}
            \varepsilon_0\,\mu_F \leq \min\{\mu_I, \mu_D\},
        \end{equation}
        then $q^*$ is locally asymptotically stable.
        \item[(c)] Without condition~\eqref{eq:rate-condition}, a fixed point with $H'(q^*_F) < 0$ is locally asymptotically stable whenever $q^*_F$ is sufficiently small or sufficiently large, but it can be unstable for intermediate values of $q^*_F$.
    \end{enumerate}
    Consequently, under~\eqref{eq:rate-condition}: in Statement~1 the non-trivial branch is stable; in Statement~3 the system is bistable for $\lambda \in (\max\{\lambda_-,\lambda_\crit\}, \lambda_+)$, with saddle-node bifurcations at $\lambda_\pm$ and hysteresis; and in Statement~4 the system is bistable for $\lambda \in (\lambda_-, \lambda_\crit)$, where the trivial fixed point coexists with a stable high-debris fixed point.
\end{theorem}
The proof is in Section~\ref{sec:proof-bifurcation}.

\begin{remark}[Launch rates above $\mu_I$]\label{rem:capacity}
    The assumption that $\lambda \leq \mu_I$ is used in an essential way.
    If $\lambda > \mu_I$, the drift of $q_I$ at the wall $\{q_I = 1\}$ is $\lambda - \mu_I - \alpha(q_F)cq_F$, which is positive when there are few fragments: launches are blocked, and by Proposition~\ref{prop:reflected-fluid} the fluid limit is the system~\eqref{eq:RFL}, in which the intact equation~\eqref{eq:FL-1} is reflected at $q_I = 1$.
    Theorem~\ref{thm:bifurcation-structure} does not cover this regime, and its conclusions change there.
    The trivial fixed point of~\eqref{eq:FL-1}--\eqref{eq:FL-3} has $q^*_I = \lambda/\mu_I > 1$ and is not in the state space.
    The clean state is instead the saturated state $(1,\ (1-P)\mu_I/\mu_D,\ 0)$, at which a small fragment population grows at rate
    \[
        \kfrags\, c\left(\alpha_0 + \frac{(1-P)\mu_I}{\mu_D}\right) - \mu_F = \mu_F\left(\frac{\mu_I}{\lambda_\crit} - 1\right),
    \]
    which does not depend on $\lambda$.
    Hence, if $\lambda_\crit > \mu_I$, no increase of the launch rate destabilizes the clean state, and $\lambda_\crit$ is a bifurcation value only if $\lambda_\crit < \mu_I$.
    Similarly, a non-trivial fixed point of~\eqref{eq:FL-1}--\eqref{eq:FL-3} is relevant only if $q^*_I = \lambda/(\mu_I + \alpha(q^*_F)cq^*_F) \leq 1$; this can fail for a low-debris fixed point and hold for a high-debris one at the same launch rate, since the intact count is smaller on the high-debris branch.
    We do not analyze the regime $\lambda > \mu_I$ here.
\end{remark}

The implications of Theorem~\ref{thm:bifurcation-structure} are significant.
If $\Qsat > \overline{Q}$, the fixed-point structure is the same as with non-saturating collision avoidance.
However, $\Qsat$ is an unmeasurable aggregate quantity, and the launch rate to LEO is increasing significantly.
Thus, it is crucial to understand the dynamics when $\Qsat < \overline{Q}$.
In this regime, the system can exhibit bistability and hysteresis.
If the debris count exceeds a threshold, a small launch-rate backoff will not revert the system to the branch with lesser debris: the launch rate needs to decrease below the hysteresis window.
This is a difficult challenge, and we address it in Section~\ref{sec:admission-control}.

\begin{remark}[The role of the rate condition]\label{rem:hopf}
    Condition~\eqref{eq:rate-condition} says that fragments do not decay much faster than intacts are retired and derelicts decay.
    It cannot simply be dropped.
    At a non-trivial fixed point, the fragment equation~\eqref{eq:FL-3} has a \emph{positive} self-feedback $\partial \dot q_F/\partial q_F = \kfrags\, c\, q^*_I q^*_F \alpha'(q^*_F) > 0$, because additional fragments degrade collision avoidance; the fixed point is stabilized only through the depletion of intacts and derelicts, which acts at rates $\mu_I$ and $\mu_D$.
    When fragments are fast relative to these rates, the stabilization arrives too late, and a fixed point with $H'(q^*_F) < 0$ loses stability through a pair of complex eigenvalues crossing the imaginary axis.
    Section~\ref{sec:hopf} makes this precise: it identifies the sharp threshold on $\mu_F$, shows that the loss of stability is a Hopf bifurcation, and gives an example.
    This phenomenon is absent when $\alpha(\cdot) = \alpha_0$, since then the self-feedback vanishes.
\end{remark}

\subsection{Oscillatory Instability: a Hopf Bifurcation}\label{sec:hopf}

Statement~(c) of Theorem~\ref{thm:bifurcation-structure} leaves open what happens at a fixed point with $H'(q^*_F) < 0$ when the rate condition~\eqref{eq:rate-condition} fails.
In this section we characterize exactly when such a fixed point is unstable, and we show that it loses its stability through a Hopf bifurcation: the low-debris equilibrium is replaced by a periodic orbit, along which debris outbreaks recur.

For $q > 0$, define
\begin{gather*}
    \sigma_I(q) := \mu_I + \alpha(q) c q, \qquad \sigma_D(q) := \mu_D + cq, \qquad e(q) := \frac{q\alpha'(q)}{\alpha(q)}, \qquad \kappa(q) := 1 - \alpha(q) - q\alpha'(q),
    \\
    \theta(q) := \frac{\alpha(q)\sigma_D(q)}{\alpha(q)\sigma_D(q) + (1-P)\mu_I}.
\end{gather*}
By~\eqref{eq:fp-intact-derelict}, $\theta(q^*_F) = \alpha(q^*_F)q^*_I/(\alpha(q^*_F)q^*_I + q^*_D)$ is the fraction of the fragment production which is due to intacts at a non-trivial fixed point.
These functions, and $H$, do not depend on $\mu_F$, $\kfrags$, or $\lambda$: changing $\mu_F$ or $\kfrags$ only rescales $\cM$.
In particular, for every $z > 0$ and every $\mu_F > 0$ there is exactly one launch rate, $\lambda = \cM(z)$, for which there is a non-trivial fixed point with $q^*_F = z$, and whether $H'(z) < 0$ does not depend on $\mu_F$.
This launch rate satisfies the standing assumption $\lambda \leq \mu_I$ if and only if~\eqref{eq:admissible} holds.
We can therefore ask how large $\mu_F$ can be before this fixed point becomes unstable.
Suppressing the argument $z$, let
\begin{align*}
    p_0(z) &:= (\sigma_I + \sigma_D)\,\sigma_I\sigma_D,
    \\
    p_1(z) &:= -e\theta(\sigma_I + \sigma_D)^2 + cz\big(\theta(1-\kappa)\sigma_I + (1-\theta)\kappa\sigma_D\big),
    \\
    p_2(z) &:= e\theta\big(e\theta(\sigma_I + \sigma_D) - cz(1-\theta\kappa)\big).
\end{align*}

\begin{proposition}[Hopf threshold]\label{prop:hopf}
    Fix $\mu_I, \mu_D, c, \alpha_0, \Qsat, P$, and let $z > 0$ with $H'(z) < 0$.
    The polynomial $m \mapsto p_0(z) + p_1(z)m + p_2(z)m^2$ has exactly one root $\mu_F^H(z)$ in the interval $\big(0, \frac{\sigma_I(z)+\sigma_D(z)}{e(z)\theta(z)}\big)$, and
    \begin{equation}\label{eq:hopf-bounds}
        \frac{\min\{\mu_I,\mu_D\}}{e(z)} < \mu_F^H(z) < \frac{\sigma_I(z)+\sigma_D(z)}{e(z)\theta(z)}.
    \end{equation}
    Let $\mu_F$ and $\kfrags$ satisfy~\eqref{eq:admissible}, so that $\lambda = \cM(z) \leq \mu_I$, let $q^*$ be the non-trivial fixed point with $q^*_F = z$, and let $c_1, c_2$ be as in Lemma~\ref{lem:linearization}.
    \begin{enumerate}
        \item If $\mu_F < \mu_F^H(z)$, then $q^*$ is locally asymptotically stable.
        \item If $\mu_F = \mu_F^H(z)$, then $c_1, c_2 > 0$ and the eigenvalues of the Jacobian at $q^*$ are $-c_2$ and $\pm i\omega$, with $\omega = \sqrt{c_1}$.
        \item If $\mu_F > \mu_F^H(z)$, then $q^*$ is unstable: the Jacobian at $q^*$ has one negative eigenvalue and two eigenvalues with positive real part.
    \end{enumerate}
    Moreover, fix also $\mu_F$ and $\kfrags$, and let $z_H > 0$ satisfy $H'(z_H) < 0$, $\mu_F^H(z_H) = \mu_F$, $(\mu_F^H)'(z_H) \neq 0$, and $\cM(z_H) < \mu_I$.
    Then on the branch of non-trivial fixed points through $q^*_F = z_H$, parametrized by $\lambda$ near $\lambda_H := \cM(z_H)$, the Jacobian has a pair of eigenvalues $x(\lambda) \pm i\omega(\lambda)$ with
    \[
        x(\lambda_H) = 0, \qquad \omega(\lambda_H) = \sqrt{c_1} > 0, \qquad \operatorname{sign}x'(\lambda_H) = -\operatorname{sign}\,(\mu_F^H)'(z_H) \neq 0,
    \]
    that is, the pair crosses the imaginary axis transversally as $\lambda$ passes through $\lambda_H$.
\end{proposition}
The proof is in Section~\ref{sec:proof-hopf}.

Proposition~\ref{prop:hopf} has three consequences.
First, it identifies the sharp version of the rate condition.
Let
\[
    \mu_F^\dagger := \inf\big\{\mu_F^H(z) : z > 0,\ H'(z) < 0\big\} \geq \frac{\min\{\mu_I,\mu_D\}}{\varepsilon_0},
\]
where the inequality is~\eqref{eq:hopf-bounds} and~\eqref{eq:e-bound}.
If $\mu_F < \mu_F^\dagger$, then every non-trivial fixed point with $H'(q^*_F) < 0$ is locally asymptotically stable, and the conclusions which Theorem~\ref{thm:bifurcation-structure} draws under~\eqref{eq:rate-condition} hold.
The threshold $\mu_F^\dagger$ does not depend on $\kfrags$ or on $\lambda$.

Second, the instability cannot be excluded by the choice of $\Qsat$, but it must occur below capacity to be relevant.
Since $\mu_F^H(z) < \infty$, we have $\mu_F^\dagger < \infty$, and if $\mu_F > \mu_F^\dagger$ then there is a non-empty open set of values of $z$ at which $H'(z) < 0$ and $\mu_F > \mu_F^H(z)$; by Statement~(c) of Theorem~\ref{thm:bifurcation-structure}, this set is bounded away from $0$ and from $\infty$.
The fixed point with $q^*_F = z$ is present at a launch rate which satisfies $\lambda \leq \mu_I$ exactly when~\eqref{eq:admissible} holds.
Hence there is an unstable fixed point with $H'(q^*_F) < 0$ at some launch rate $\lambda \leq \mu_I$ if and only if there is $z > 0$ with $H'(z) < 0$ and
\begin{equation}\label{eq:hopf-window}
    \mu_F^H(z) < \mu_F \leq \kfrags\, c\,\mu_I H(z).
\end{equation}
Since $\mu_F^H$ and $H$ do not depend on $\kfrags$, for every $\mu_F > \mu_F^\dagger$ condition~\eqref{eq:hopf-window} holds for some $z$ once $\kfrags$ is large enough; conversely, if $\kfrags\, c\,\mu_I H(z) \leq \mu_F^H(z)$ for every $z$ with $H'(z) < 0$, then the instability does not occur for any $\lambda \leq \mu_I$.
When $\Qsat > \overline{Q}$ and~\eqref{eq:hopf-window} holds at $z$, the launch rate $\cM(z)$ exceeds $\lambda_\crit$, so the trivial fixed point is unstable as well: at this launch rate the fluid limit has \emph{no} stable fixed point.

\begin{figure}[h!]
    \centering
    \includegraphics[width=\linewidth]{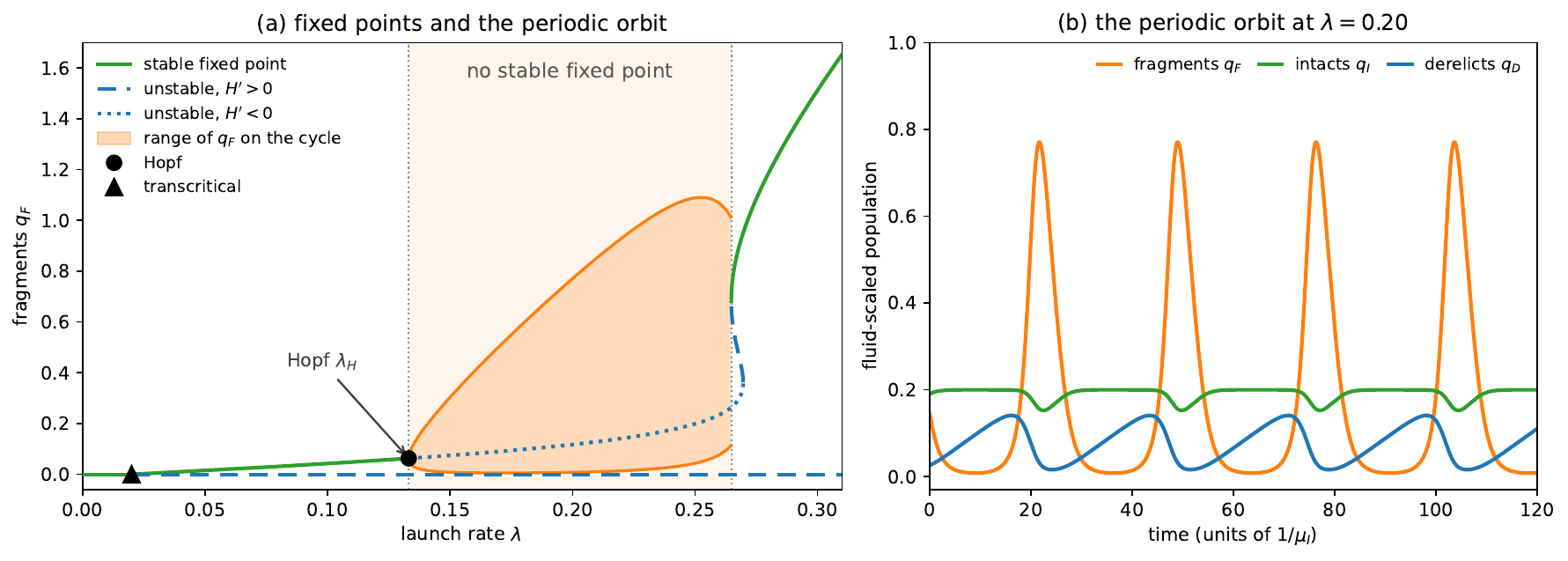}
    \caption{The Hopf bifurcation on the low-debris branch, with $\mu_I = \mu_F = c = 1$, $\mu_D = 0.01$, $\alpha_0 = 0.01$, $\Qsat = 1$, $P = 0.95$, $\kfrags = 10$. (a) Fixed points as functions of the launch rate, and the range of the fragment count on the stable periodic orbit which exists for $\lambda_H < \lambda < \lambda_-$. (b) The periodic orbit at $\lambda = 0.20$.}
    \label{fig:hopf}
\end{figure}

Third, the loss of stability is a Hopf bifurcation (see Figure~\ref{fig:hopf}).
At a boundary point $z_H$ of the unstable set the hypotheses of the second part of Proposition~\ref{prop:hopf} hold generically, and then the Hopf bifurcation theorem~\cite[Chapters~3 and~5]{kuznetsov2023elements} applies at $\lambda_H$ whenever the first Lyapunov coefficient $\ell_1$ at $(q^*, \lambda_H)$ is non-zero: a family of periodic orbits with period close to $2\pi/\omega$ bifurcates from the branch of fixed points.
If $\ell_1 < 0$ the bifurcation is supercritical, and the periodic orbits are stable and exist on the side of $\lambda_H$ on which the fixed point is unstable.

\paragraph{An example.}
Let $\mu_I = \mu_F = c = 1$, $\mu_D = 0.01$, $\alpha_0 = 0.01$, $\Qsat = 1$, $P = 0.95$, and $\kfrags = 10$.
Then $\underline{Q} \approx 0.002 < \Qsat < \overline{Q} = 9900$, and $H'$ has two positive roots, $q_{\min} \approx 0.363$ and $q_{\max} \approx 0.682$, so that
\[
    \lambda_\crit \approx 0.0200, \qquad \lambda_- \approx 0.2648, \qquad \lambda_+ \approx 0.2696.
\]
All the launch rates which we consider in this example are at most $0.31 < \mu_I$, so the standing assumption $\lambda \leq \mu_I$ holds, and $q_I(t) \leq \max\{q_I(0), \lambda/\mu_I\}$ stays well below the wall $q_I = 1$.
Evaluating $\mu_F^H$ numerically gives $\mu_F^\dagger \approx 0.650$, attained at $z \approx 0.175$ on the low-debris branch, so $\mu_F = 1$ exceeds the threshold.
(The sufficient condition~\eqref{eq:rate-condition} is far from sharp here, since $\min\{\mu_I,\mu_D\}/\varepsilon_0 \approx 0.012$: the bound is driven by $\mu_D$, whereas the instability sets in only when $\mu_F$ is comparable to $\mu_I$.)
The equation $\mu_F^H(z) = 1$ has a single solution on the low-debris branch, $z_H \approx 0.0629$, with $(\mu_F^H)'(z_H) \approx -11.5 < 0$.
The corresponding launch rate is $\lambda_H \approx 0.1331$, the critical eigenvalues are $\pm 0.229\,i$, and the third eigenvalue is $-1.004$.
Evaluating the formula of~\cite[Section~5.4]{kuznetsov2023elements} numerically gives $\ell_1 \approx -4.5 < 0$, so the Hopf bifurcation is supercritical.

The consequences are shown in Figure~\ref{fig:hopf}, where the periodic orbit is computed by numerical integration of~\eqref{eq:FL-1}--\eqref{eq:FL-3}.
The low-debris branch is stable only for $\lambda < \lambda_H$, and it is unstable on the whole interval $(\lambda_H, \lambda_+)$: the fold at $\lambda_+$, at which the low-debris branch ends, plays no role in the dynamics.
For $\lambda \in (\lambda_H, \lambda_-)$ there is no stable fixed point, and trajectories converge to a stable periodic orbit with period approximately $28/\mu_I$.
At $\lambda = 0.20$, the fragment count on this orbit ranges over $[0.008, 0.77]$, whereas the unstable fixed point has $q^*_F \approx 0.117$.
As $\lambda \uparrow \lambda_-$ the period of the orbit grows without bound (it is approximately $87$ at $\lambda = 0.2640$ and $535$ at $\lambda = 0.26475$), and the orbit disappears at $\lambda_-$, where the saddle and the high-debris fixed point are born; this is consistent with a saddle-node bifurcation on an invariant circle.
For $\lambda > \lambda_-$, trajectories converge to the high-debris fixed point.
In particular, in this example there is no bistability and no hysteresis: in the window $(\lambda_-, \lambda_+)$, in which $H$ alone would predict two stable fixed points, only the high-debris fixed point is stable.

The mechanism is that of a predator--prey cycle, in which the fragments prey on the derelicts and the intacts.
While there are few fragments, derelicts accumulate, slowly when $\mu_D$ is small, until the fragment generation rate exceeds $\mu_F$.
The outbreak which follows is amplified by the degradation of collision avoidance, since $\alpha(q_F)$ increases with $q_F$, and it consumes the derelicts and a part of the intact population (Figure~\ref{fig:hopf}(b)).
Deprived of targets, the fragments then decay at rate $\mu_F$, the derelict population is replenished by failed post-mission disposals, and the cycle repeats.
When $\mu_F$ is small, the fragments outlast the depletion of their targets, and the system settles at a fixed point instead.

The operational implications differ from those of bistability.
In the oscillatory regime, a low fragment count is a phase of the cycle and not evidence of a safe operating point, and the fragment count at the peak of the cycle is several times its value at the fixed point.
On the other hand, the onset is not abrupt: since the Hopf bifurcation is supercritical, the amplitude of the oscillation grows continuously from zero as $\lambda$ increases past $\lambda_H$, and lowering the launch rate below $\lambda_H$ restores the low-debris fixed point without hysteresis.
In the example, the regime requires $\mu_F > \mu_F^\dagger \approx 0.65\,\mu_I$, that is, fragments which decay on the time scale of a mission; this suggests that it concerns shells in which atmospheric drag removes fragments within roughly one mission lifetime.
When $\mu_F < \mu_F^\dagger$, as for fragments which are long-lived relative to missions in the example, the picture of Theorem~\ref{thm:bifurcation-structure} under~\eqref{eq:rate-condition} applies.

\subsection{Slow Passage Through the Saddle Node}\label{sec:slow-passage-result}
The analysis so far concerns a constant launch rate.
A constellation does not change its launch rate discontinuously; it grows, so the natural forcing is a slow ramp
\begin{equation}\label{eq:ramp}
    \lambda(t) = \lambda_0 + \varepsilon t, \qquad 0 < \varepsilon \ll 1.
\end{equation}
As long as $\lambda(t) \leq \mu_I$, the fluid limit under a time-varying launch rate is~\eqref{eq:FL-1}--\eqref{eq:FL-3} with $\lambda$ replaced by $\lambda(t)$ (Remark~\ref{rem:time-varying}).
This subsection describes what the fluid limit does when the ramp carries the launch rate across the upper fold $\lambda_+$ in the two-fold regime.
The answer is a \emph{bifurcation delay}: the trajectory does not leave the low-debris branch when $\lambda(t)$ crosses $\lambda_+$, but lingers near the vanished fixed point, and jumps to the high-debris branch only after the launch rate has overshot $\lambda_+$ by an amount of order $\varepsilon^{2/3}$.

We work in the two-fold regime, where $\Qsat > \underline{Q}$, and $H'$ has two distinct positive roots $q_{\min} < q_{\max}$.
By the Proof of Statement~2 of Theorem~\ref{thm:bifurcation-structure} these roots are simple, so $H''(q_{\min}) > 0 > H''(q_{\max})$, and $\cM$ increases on $(0, q_{\min})$ from $\lambda_\crit$ to $\lambda_+ = \cM(q_{\min})$, decreases on $(q_{\min}, q_{\max})$ to $\lambda_- = \cM(q_{\max})$, and increases on $(q_{\max}, \infty)$.
For $\lambda \in (\lambda_\crit, \lambda_+)$, let $q^{\mathrm{lo}}(\lambda)$ be the non-trivial fixed point with $q^*_F < q_{\min}$ (the low debris branch); for $\lambda > \lambda_-$, let $q^{\mathrm{hi}}(\lambda)$ be the one with $q^*_F > q_{\max}$ (the high debris branch); and let $q^\dagger$ be the fixed point at $\lambda = \lambda_+$ with $q^*_F = q_{\min}$ (the fold point).
With $\sigma_I, \sigma_D$ and $c_1$ as in Lemma~\ref{lem:linearization}, evaluated at the fold point, define
\begin{equation}\label{eq:sp-coefficients}
    \beta_\lambda := \frac{\mu_F\, q_{\min}\,\sigma_I\sigma_D}{\lambda_+\, c_1},
    \qquad
    \beta_2 := \beta_\lambda\,\frac{\lambda_+ H''(q_{\min})}{2H(q_{\min})}.
\end{equation}
These are the coefficients of the saddle-node normal form $\dot u = \beta_\lambda(\lambda - \lambda_+) + \beta_2 u^2$ of the dynamics on a center manifold at the fold point, where $u$ is the displacement of $q_F$ from $q_{\min}$ (Section~\ref{sec:slow-passage}); under Assumption~(A1) below, $c_1 > 0$, so $\beta_\lambda, \beta_2 > 0$.

We make two assumptions.
\begin{itemize}
    \item[(A1)] $\lambda_+ < \mu_I$, and $c_2 > 0$ and $c_1c_2 - c_0 > 0$ at every non-trivial fixed point with $q^*_F \in (0, q_{\min}] \cup (q_{\max}, \infty)$ and launch rate at most $\mu_I$.
    \item[(A2)] At $\lambda = \lambda_+$, the orbit which leaves the fold point $q^\dagger$ in the direction of increasing $q_F$ converges to $q^{\mathrm{hi}}(\lambda_+)$.
\end{itemize}
Rate condition~\eqref{eq:rate-condition} implies the second part of~(A1) (see Lemma~\ref{lem:linearization}(iii)), and it also ensures that every fixed point on the low- and high-debris branches is locally asymptotically stable (see Lemma~\ref{lem:linearization}(ii)).
In~(A2), the orbit is well defined: at $\lambda = \lambda_+$ the fold point has a one-dimensional center manifold, on which the flow is $\dot u = \beta_2 u^2 + O(u^3)$, and~(A2) concerns its half $\{u > 0\}$.
Assumption~(A2) is a condition on a single orbit.
It does not follow from the fixed-point structure alone, because the fluid limit is not a monotone system and can have periodic orbits (Section~\ref{sec:hopf}); it is checked by integrating~\eqref{eq:FL-1}--\eqref{eq:FL-3} once.

\begin{theorem}\label{thm:slow-passage}
    Assume the two-fold regime,~(A1), and~(A2).
    Let $\lambda_0 \in (\max\{\lambda_-, \lambda_\crit\}, \lambda_+)$ and $\Lambda \in (\lambda_+, \mu_I]$, and let $q$ be the solution of~\eqref{eq:FL-1}--\eqref{eq:FL-3} with $\lambda$ replaced by the ramp~\eqref{eq:ramp} and $q(0) = q^{\mathrm{lo}}(\lambda_0)$, for $t \in [0, (\Lambda - \lambda_0)/\varepsilon]$.
    Let $\Omega_0$ be the smallest positive zero of $s \mapsto \mathrm{Ai}(-s)$, where $\mathrm{Ai}$ is the Airy function.
    Then, as $\varepsilon \to 0$:
    \begin{enumerate}
        \item[\textup{(i)}] for every $\delta > 0$, $\sup\big\{|q(t) - q^{\mathrm{lo}}(\lambda(t))| : \lambda(t) \leq \lambda_+ - \delta\big\} = O(\varepsilon)$;
        \item[\textup{(ii)}] for every $K > 0$ and $\delta \in (0, \Omega_0)$,
        \[
            \sup\Big\{|q(t) - q^\dagger| : -K\varepsilon^{2/3} \leq \lambda(t) - \lambda_+ \leq (\Omega_0 - \delta)(\beta_\lambda\beta_2)^{-1/3}\varepsilon^{2/3}\Big\} = O(\varepsilon^{1/3});
        \]
        in particular, at the time $t_+$ at which $\lambda(t_+) = \lambda_+$, $|q(t_+) - q^\dagger| = O(\varepsilon^{1/3})$.
        In time, since $\lambda(t) - \lambda_+ = \varepsilon(t - t_+)$, the window above is $-K\varepsilon^{-1/3} \leq t - t_+ \leq (\Omega_0 - \delta)(\beta_\lambda\beta_2)^{-1/3}\varepsilon^{-1/3}$;
        \item[\textup{(iii)}] for every sufficiently small $\Delta > 0$, the first time $t_{\Delta}$ at which $q_F(t_{\Delta}) = q_{\min} + \Delta$ satisfies
        \begin{equation}\label{eq:sp-law}
            \lambda(t_{\Delta}) - \lambda_+ = \Omega_0\,(\beta_\lambda\beta_2)^{-1/3}\,\varepsilon^{2/3} + O\big(\varepsilon\log(1/\varepsilon)\big);
        \end{equation}
        equivalently, the jump comes a time
        \begin{equation}\label{eq:sp-law-time}
            t_{\Delta} - t_+ = \Omega_0\,(\beta_\lambda\beta_2)^{-1/3}\,\varepsilon^{-1/3} + O\big(\log(1/\varepsilon)\big)
        \end{equation}
        after the launch rate crosses $\lambda_+$;
        \item[\textup{(iv)}] for every $\delta > 0$, $\sup\big\{|q(t) - q^{\mathrm{hi}}(\lambda(t))| : \lambda_+ + \delta \leq \lambda(t) \leq \Lambda\big\} = O(\varepsilon)$.
    \end{enumerate}
    The constants depend on $\delta$, $\Delta$, $K$, $\lambda_0$, $\Lambda$ and the parameters of the model, but not on $\varepsilon$.
\end{theorem}
The proof is in Section~\ref{sec:slow-passage}.

Theorem~\ref{thm:slow-passage} says that the tipping launch rate observed under a slow ramp exceeds the static threshold $\lambda_+$ by $\Omega_0(\beta_\lambda\beta_2)^{-1/3}\varepsilon^{2/3}$.
In time, the jump occurs $\Omega_0(\beta_\lambda\beta_2)^{-1/3}\varepsilon^{-1/3}$ after the launch rate crosses $\lambda_+$ (equation~\eqref{eq:sp-law-time}).
The overshoot in the launch rate shrinks as the ramp slows, but the waiting time grows.
Until $\lambda(t) - \lambda_+$ reaches $(\Omega_0 - \delta)(\beta_\lambda\beta_2)^{-1/3}\varepsilon^{2/3}$, for any fixed $\delta > 0$, the state stays within $O(\varepsilon^{1/3})$ of the fold point (Statement~(ii)), so the system has passed its tipping point while its observables remain close to their values at the threshold.
Once $q_F$ has left the fold region, the transition to the high-debris branch takes a time of order one, during which the launch rate changes by $O(\varepsilon)$.
These statements concern the fluid limit; at finite $n$, fluctuations of order $n^{-1/2}$ can trigger the escape earlier, and the deterministic delay is the relevant description when $\varepsilon^{2/3}$ is large compared with $n^{-2/3}$, that is, when $n\varepsilon$ is large~\cite{berglund2002hysteresis}.

\section{Admission Control}\label{sec:admission-control}

Theorem~\ref{thm:bifurcation-structure} shows that this constant-input system can be bistable.
The natural mitigation is \emph{admission control}: a launch-licensing authority that blocks launches as a function of the observed state.
This raises a sharp design question.
The most operationally natural signal is the shell's own occupancy---the number $Q_I$ of occupied intact slots---since that is what a launch operator already meters through blocking.
The alternative is to sense the debris, or the conjunction load, directly.
We first show that blocking launches as the shell fills is insufficient to avoid bistability, and then that debris-indexed blocking can remove it.
Throughout this section, $\lambda \leq \mu_I$, so that by Proposition~\ref{prop:reflected-fluid} the fluid limit is not reflected.

Recall from Section~\ref{sec:hopf} that
\[
    \sigma_I(q) := \mu_I + \alpha(q)\,c\,q, \qquad \sigma_D(q) := \mu_D + c\,q,
\]
with these definitions now extended to cover $q \geq 0$.
These are the total rates at which a single intact and a single derelict leave the shell when the fragment count is $q$.
We have $$H = \alpha/\sigma_I + (1-P)\mu_I/(\sigma_I\sigma_D).$$
Recall from Section~\ref{sec:bistability} that at a non-trivial fixed point, the intact and derelict counts are determined by the fragment count, $q^*_I = J/\sigma_I(q^*_F)$ and $q^*_D = (1-P)\mu_I q^*_I/\sigma_D(q^*_F)$, where $J$ is the \emph{admitted} launch rate at the fixed point, and that the fragment balance~\eqref{eq:FL-FP-3} then reads
\[
    c\,\kfrags\,J\,H(q^*_F) = \mu_F, \qquad\text{that is,}\qquad J = \cM(q^*_F) = \frac{\mu_F}{c\,\kfrags\,H(q^*_F)}.
\]
For constant launches $J = \lambda$, and the fixed points are the intersections of the horizontal line at height $\lambda$ with the graph of $\cM$.
Bistability is thus the non-monotonicity of $\cM$.
Admission control replaces the constant $J$ by a state-dependent schedule.
We address the design question of whether such a schedule can remove the bistability.

\begin{definition}\label{def:admissible}
    An \emph{admission policy} blocks launches as a function of the observed state.
    We distinguish two classes of admission policies: \emph{occupancy feedback}, in which launches are offered at rate $n\bar\lambda(Q_I/n)$ for a Lipschitz function $\bar\lambda : [0,1] \to (0,\mu_I]$, reacting only to the number of occupied intact slots; and \emph{debris feedback}, in which launches are offered at rate $n\lambda\,\rho(Q_F/n)$, a throttle $\rho$ of the uncontrolled launch rate $\lambda$ indexed by the fragment count.
    A debris-feedback policy is \emph{admissible} if
    \begin{enumerate}
        \item[\textup{(i)}] it only blocks launches: $0 < \rho(q_F) \le 1$ for all $q_F \ge 0$; and
        \item[\textup{(ii)}] $\rho$ is Lipschitz on $[0,\infty)$.
    \end{enumerate}
\end{definition}
Condition~\textup{(i)} is the operational constraint that a launch-licensing authority may block launches but not create them; condition~\textup{(ii)} is the regularity under which the throttled process has a fluid limit (Lemma~\ref{lem:rho-lip}, Theorem~\ref{thm:clip-fluid}).

Under debris feedback the fluid limit is~\eqref{eq:FL-1}--\eqref{eq:FL-3} with the launch term $\lambda$ of~\eqref{eq:FL-1} replaced by the admitted rate $\lambda\rho(q_F)$ as in~\eqref{eq:cFL-1}--\eqref{eq:cFL-3} below, which by Theorem~\ref{thm:clip-fluid} is the fluid limit of the throttled process when $\lambda \le \mu_I$:
\begin{align}
    \dot q_I &= \lambda\,\rho(q_F) - \mu_I q_I - \alpha(q_F)cq_Iq_F, \label{eq:cFL-1}\\
    \dot q_D &= (1-P)\mu_I q_I - \mu_D q_D - cq_Dq_F, \label{eq:cFL-2}\\
    \dot q_F &= \kfrags\,\alpha(q_F) c q_F q_I + \kfrags\, c q_F q_D - \mu_F q_F. \label{eq:cFL-3}
\end{align}

In state $(Q_I, Q_D, Q_F)$ the launch stream of the $n$-th system is \emph{thinned} to rate $n\lambda\,\rho(Q_F/n)$---equivalently, each offered launch is blocked with probability $1 - \rho(Q_F/n)$.
All other transitions (missions, collisions, disposals, decay) are unchanged.
Denote the resulting process $Q^{n,\rho}$.

\begin{theorem}\label{thm:clip-fluid}
    Let $\rho$ be admissible, $\lambda > 0$, and suppose $Q^{n,\rho}(0)/n \to q(0) \in \mathcal S$.
    Then almost surely, uniformly on compacts, $\tfrac1n Q^{n,\rho}(t) \to q^\rho(t)$, the unique solution of~\eqref{eq:cFL-1}--\eqref{eq:cFL-3} with~\eqref{eq:cFL-1} reflected at $q_I = 1$; if $\lambda \leq \mu_I$ the reflection is inactive and $q^\rho$ solves~\eqref{eq:cFL-1}--\eqref{eq:cFL-3}.
    For $\lambda \leq \mu_I$, the non-trivial fixed points satisfy $q_I^* = \lambda\rho(q_F^*)/\sigma_I(q^*_F)$ and
    \[
        c\,\kfrags\,\lambda\,G(q_F^*) = \mu_F, \qquad G = \rho H.
    \]

\end{theorem}

For an admissible policy we write
    \[
        G := \rho H, \qquad \cM_\rho(q_F) := \frac{\mu_F}{c\,\kfrags\,G(q_F)},
    \]
so that the non-trivial fixed points under the launch rate $\lambda$ are the solutions of $\lambda = \cM_\rho(q^*_F)$ (see Lemma~\ref{lem:G-bifurcation} below).
An admissible policy is \emph{single-branch} if in addition $G$ is non-increasing, so that $\cM_\rho$ is non-decreasing and its graph has no fold.
We now show that $G$ is the bifurcation function of the controlled system in the same sense in which $H$ is that of the uncontrolled one: its value at the origin locates the transcritical bifurcation, its level sets are the non-trivial fixed points, and the sign of $G'$ is the sign of the determinant of the Jacobian, so that the folds are critical points of $G$.

\begin{lemma}\label{lem:G-bifurcation}
Let $\rho$ be admissible and $\lambda > 0$, and write $f = (f_I, f_D, f_F)$ for the right-hand side of~\eqref{eq:cFL-1}--\eqref{eq:cFL-3}.
\begin{enumerate}
    \item[\textup{(i)}] The point $q^0 := \big(\lambda\rho(0)/\mu_I,\ (1-P)\lambda\rho(0)/\mu_D,\ 0\big)$ is a fixed point. If $\rho$ is differentiable at $0$, the eigenvalues of the Jacobian $Df(q^0)$ are $-\mu_I$, $-\mu_D$ and $c\,\kfrags\,\lambda\,G(0) - \mu_F$, so $q^0$ is locally asymptotically stable for $\lambda < \cM_\rho(0)$ and unstable for $\lambda > \cM_\rho(0)$.
    \item[\textup{(ii)}] The fixed points with $q_F > 0$ are exactly the points
    \[
        q^*(z) := \Big(\frac{\lambda\rho(z)}{\sigma_I(z)},\ \frac{(1-P)\mu_I\,\lambda\rho(z)}{\sigma_I(z)\,\sigma_D(z)},\ z\Big),
        \qquad z > 0, \qquad c\,\kfrags\,\lambda\,G(z) = \mu_F,
    \]
    one for each solution $z$ of the last equation; equivalently, $\lambda = \cM_\rho(z)$.
    \item[\textup{(iii)}] If $\rho$ is differentiable at a solution $z$ of~\textup{(ii)}, then
    \[
        \det Df\big(q^*(z)\big) = \sigma_I(z)\,\sigma_D(z)\,\mu_F\,z\,\frac{G'(z)}{G(z)}.
    \]
    Hence $Df(q^*(z))$ has a zero eigenvalue if and only if $G'(z) = 0$, that is, if and only if $\cM_\rho'(z) = 0$, so every fold of the branch $\lambda = \cM_\rho(z)$ is a critical point of $G$; and if $G'(z) > 0$, then $Df(q^*(z))$ has a positive real eigenvalue and $q^*(z)$ is unstable.
\end{enumerate}
For $\rho \equiv 1$ we have $G = H$, and \textup{(i)--(iii)} are the corresponding facts for the uncontrolled fluid limit: $\lambda_\crit = \cM(0)$, the fixed-point relation~\eqref{eq:FL-FP-3}, and the constant term $c_0 = -\det Df$ of Lemma~\ref{lem:linearization}.
\end{lemma}

\begin{proof}
    Write $f_F = q_F\gamma_F$ with $\gamma_F(q) := c\,\kfrags\,\big(\alpha(q_F)\,q_I + q_D\big) - \mu_F$.
    For fixed $q_F = z$, the equations $f_I = 0$ and $f_D = 0$ are linear in $(q_I, q_D)$:
    \[
        A(z)\begin{pmatrix} q_I \\ q_D \end{pmatrix} = -\begin{pmatrix} \lambda\rho(z) \\ 0 \end{pmatrix},
        \qquad
        A(z) := \begin{pmatrix} -\sigma_I(z) & 0 \\ (1-P)\mu_I & -\sigma_D(z) \end{pmatrix}.
    \]
    Since $\det A(z) = \sigma_I(z)\sigma_D(z) > 0$, they have the unique solution $X(z) := \lambda\rho(z)/\sigma_I(z)$, $Y(z) := (1-P)\mu_I X(z)/\sigma_D(z)$.
    The throttle enters only through the admitted rate $\lambda\rho(z)$:
    \begin{equation}\label{eq:G-reduction}
        \gamma_F\big(X(z), Y(z), z\big)
        = c\,\kfrags\,\lambda\rho(z)\Big(\frac{\alpha(z)}{\sigma_I(z)} + \frac{(1-P)\mu_I}{\sigma_I(z)\,\sigma_D(z)}\Big) - \mu_F
        = c\,\kfrags\,\lambda\,G(z) - \mu_F,
    \end{equation}
    by the definition of $H$.

    \emph{(i).} At $z = 0$, $(X(0), Y(0)) = (\lambda\rho(0)/\mu_I, (1-P)\lambda\rho(0)/\mu_D)$, and $f_F = 0$ because $q_F = 0$; so $q^0$ is a fixed point.
    At $q^0$ we have $\partial f_F/\partial q_I = c\,\kfrags\,\alpha_0\,q_F = 0$ and $\partial f_F/\partial q_D = c\,\kfrags\,q_F = 0$, while $\partial f_F/\partial q_F = \gamma_F(q^0) =c\,\kfrags\,\lambda\,G(0) - \mu_F$ by~\eqref{eq:G-reduction}.
    The Jacobian is therefore block upper triangular, with diagonal blocks $A(0)$, whose eigenvalues are $-\sigma_I(0) = -\mu_I$ and $-\sigma_D(0) = -\mu_D$, and the scalar $c\,\kfrags\,\lambda\,G(0) - \mu_F$, which is negative exactly when $\lambda < \cM_\rho(0)$.

    \emph{(ii).} A fixed point with $q_F = z > 0$ has $(q_I, q_D) = (X(z), Y(z))$, and then $f_F = 0$ if and only if $\gamma_F(X(z), Y(z), z) = 0$, which by~\eqref{eq:G-reduction} is $c\,\kfrags\,\lambda\,G(z) = \mu_F$.

    \emph{(iii).} The upper-left $2 \times 2$ block of $Df$ is $A(z)$, which does not involve the throttle, so by the Schur-complement formula
    \[
        \det Df = \det A(z)\, s(z), \qquad
        s(z) := \partial_{q_F} f_F - \nabla_{(q_I, q_D)} f_F\; A(z)^{-1}\, \partial_{q_F}(f_I, f_D)^{\top},
    \]
    with all partial derivatives evaluated at $q^*(z)$.
    The identities $f_I(X(\zeta), Y(\zeta), \zeta) = 0 = f_D(X(\zeta), Y(\zeta), \zeta)$ hold for every $\zeta \ge 0$, and $X$ and $Y$ are differentiable at $z$ because $\rho$ is; differentiating the identities at $\zeta = z$ gives $A(z)\,(X'(z), Y'(z))^{\top} = -\partial_{q_F}(f_I, f_D)^{\top}$.
    Thus $A(z)^{-1}\,\partial_{q_F}(f_I, f_D)^{\top} = -(X'(z), Y'(z))^{\top}$. %
    Substituting this into $s(z)$ and using the chain rule,
    \[
        s(z) = \partial_{q_F} f_F + \nabla_{(q_I, q_D)} f_F \cdot \big(X'(z), Y'(z)\big)^{\top}
        = \frac{d}{d\zeta}\, f_F\big(X(\zeta), Y(\zeta), \zeta\big)\Big|_{\zeta = z},
    \]
    and by~\eqref{eq:G-reduction} the function being differentiated is $\zeta\,\big(c\,\kfrags\,\lambda\,G(\zeta) - \mu_F\big)$.
    At $z$, where $c\,\kfrags\,\lambda\,G(z) = \mu_F$, its derivative is $c\,\kfrags\,\lambda\,z\,G'(z) = \mu_F\,z\,G'(z)/G(z)$, which proves the formula.
    As $\sigma_I, \sigma_D, \mu_F, z$ and $G(z)$ are positive, $\det Df = 0$ exactly when $G'(z) = 0$, and $\cM_\rho' = -\mu_F G'/(c\,\kfrags\,G^2)$ vanishes at the same points.
    Finally, $\det Df$ is the product of the three eigenvalues, and a pair of complex-conjugate eigenvalues contributes a positive factor; so if $\det Df > 0$, at least one eigenvalue is real and positive.
\end{proof}

The controlled bifurcation diagram is therefore the graph of $\cM_\rho = \mu_F/(c\,\kfrags\,G)$, exactly as the uncontrolled one is the graph of $\cM = \mu_F/(c\,\kfrags\,H)$.
Every statement of Section~\ref{sec:bistability} that is phrased through $H$ carries over to debris feedback with $H$ replaced by $G$, and designing a debris-feedback policy amounts to choosing $G$ subject to $0 < G \le H$ (admissibility~(i)).
The explicit product structure $G = \rho H$ is special to debris feedback: the throttle depends only on $q_F$, which is held fixed in the reduction, so it factors out of the fragment balance.
Under occupancy feedback the admitted rate depends on $q_I$, which is itself determined by $q_F$ through $\bar\lambda(q_I) = \sigma_I(q_F)\,q_I$.
One instead needs to consider $J(z)H(z)$, with the implicitly defined admitted rate $J$ of Proposition~\ref{prop:occupancy}.

\begin{proposition}
    \label{prop:occupancy}
    Let $\bar\lambda : [0,1] \to (0,\mu_I]$ be Lipschitz and non-increasing, so that launches are blocked as the shell fills.
    \begin{enumerate}
        \item[\textup{(i)}] For every $z \geq 0$ the equation $\bar\lambda(q) = \sigma_I(z)\,q$ has a unique root $q_I^\circ(z) \in (0,1]$, and the non-trivial fixed points of the fluid limit under $\bar\lambda$ are the points $(q_I^\circ(z),\ (1-P)\mu_I q_I^\circ(z)/\sigma_D(z),\ z)$ with $z > 0$ and
        \[
            c\,\kfrags\,J(z)\,H(z) = \mu_F, \qquad J(z) := \bar\lambda\big(q_I^\circ(z)\big) = \sigma_I(z)\,q_I^\circ(z),
        \]
        the admitted launch rate. The function $J$ is non-decreasing.
        \item[\textup{(ii)}] On every interval on which $H$ is increasing, so is $JH$, and every fixed point with $q^*_F$ in such an interval, near whose intact count $\bar\lambda$ is continuously differentiable, is unstable.
        In particular, in the two-fold regime, whenever $c\,\kfrags\,J(q_{\min})H(q_{\min}) < \mu_F < c\,\kfrags\,J(q_{\max})H(q_{\max})$, there are at least two non-trivial fixed points: an unstable one with $q^*_F \in (q_{\min}, q_{\max})$ and one with $q^*_F > q_{\max}$.
        \item[\textup{(iii)}] Let $\bar\lambda = s\,g$ for a fixed non-increasing profile $g$ and a level $s \in (0, \bar s]$, $\bar s := \mu_I/g(0)$, and write $J_s$ for the corresponding admitted rate. If $c\,\kfrags\,J_{\bar s}(q_{\max})H(q_{\max}) > \mu_F$, then the condition of~(ii) holds for every level $s$ in a non-empty interval.
    \end{enumerate}
\end{proposition}

The content of Proposition~\ref{prop:occupancy} is that an occupancy cap acts on the wrong side of the debris feedback.
In this model the occupancy does carry information about the debris: at a fixed point $q^*_I = J/\sigma_I(q^*_F)$, and collisions deplete the shell as fragments accumulate.
But a cap which blocks launches as the shell fills responds to that depletion by \emph{admitting more} launches, so the admitted rate $J$ rises with the debris level and the fold of $H$ survives, and can only be shifted along the launch axis.
An occupancy rule could remove the fold only by \emph{reducing} launches as the shell empties, that is, by using $Q_I$ as a proxy for the debris; such a rule is a form of debris feedback, which we now study directly.

\begin{proposition}
    \label{prop:debris}
    Consider debris feedback with a strictly decreasing $C^1$ throttle $\rho$.
    Then the non-trivial fixed points of the fluid limit solve $c\,\kfrags\,\lambda\,G(q_F) = \mu_F$ with $G = \rho H$.
    If $G$ is strictly decreasing on $[0,\infty)$, the fluid limit has at most one non-trivial fixed point at every launch rate, so the graph of $\cM_\rho$ has no fold.
    Since $H'(q_F) > 0$ only on a bounded set $\cI$ (elsewhere $H$ is non-increasing), it suffices that
    \[
    -\frac{\rho'(q_F)}{\rho(q_F)} > \frac{H'(q_F)}{H(q_F)}
    \qquad\text{for all } q_F \in \cI,
    \]
    equivalently $\dfrac{d}{dq_F}\log \rho(q_F) < -\dfrac{d}{dq_F}\log H(q_F)$ on $\cI$.
\end{proposition}

\begin{corollary}
\label{cor:exp}
Because $\cI$ is bounded and $H\in C^1$ with $H>0$, the quantity
$\eta^\ast := \sup_{q_F \in \cI}\frac{d}{dq_F}\log H(q_F)$ is finite.
The exponential throttle $\rho(q_F) = e^{-\eta q_F}$ satisfies $-\frac{d}{dq_F}\log \rho = \eta$, so for any $\eta > \eta^\ast$ the hypothesis of Proposition~\ref{prop:debris} holds and the fluid limit has at most one non-trivial fixed point at every launch rate.
\end{corollary}

Removing the fold settles the number of fixed points, but not their stability.
The throttle closes a second feedback loop, from the fragments to the launches, and a steep throttle acting with the lag of the intact and derelict populations can itself destabilize a fixed point.
We address this next.
Write
\[
    \psi_\rho(z) := -\frac{\rho'(z)}{\rho(z)}
\]
for the logarithmic slope of the throttle.

At a non-trivial fixed point $q^*$ of the controlled fluid limit with $q^*_F = z > 0$, we use the notation of Section~\ref{sec:proof-bifurcation}:
\begin{gather*}
    b := c\,z, \qquad \sigma_I := \sigma_I(z) = \mu_I + \alpha(z)\,c\,z, \qquad \sigma_D := \sigma_D(z) = \mu_D + b,
    \\
    e := \frac{z\,\alpha'(z)}{\alpha(z)}, \qquad \kappa := 1 - \alpha(z) - z\,\alpha'(z), \qquad \theta := \frac{\alpha(z)\,q^*_I}{\alpha(z)\,q^*_I + q^*_D}, \qquad w := e\,\theta\,\mu_F,
    \\
    c_2 := \sigma_I + \sigma_D - w.
\end{gather*}
Here $\theta$ is the fraction of the fragment production which is due to intacts.
The throttle does not enter the derelict equation, so $q^*_D/q^*_I = (1-P)\mu_I/\sigma_D(z)$ as without control, and $\theta = \theta(z)$ is the function of Section~\ref{sec:hopf}; moreover $\theta \in (0,1)$, and $e, \kappa \in (0,1)$ by~\eqref{eq:e-kappa}.
With $f = (f_I, f_D, f_F)$ the right-hand side of~\eqref{eq:cFL-1}--\eqref{eq:cFL-3}, the quantity $w$ is the entry $\partial f_F/\partial q_F$ of the Jacobian at $q^*$: with $q_I$ and $q_D$ held fixed, extra fragments raise the avoidance-failure probability $\alpha$, and with it their own production.
The throttle enters the Jacobian only through the entry $\partial f_I/\partial q_F$. %

\begin{lemma}\label{lem:clip-linearization}
    Let $\rho$ be admissible and differentiable at $z > 0$, and let $q^*$ be a non-trivial fixed point of the controlled fluid limit with $q^*_F = z$.
    With $\sigma_I, \sigma_D, w, \theta, \kappa$ and $c_2$ as in Lemma~\ref{lem:linearization}, the characteristic polynomial of the Jacobian at $q^*$ is $s^3 + c_2s^2 + c_1^\rho s + c_0^\rho$, with
    \[
        c_0^\rho = -\sigma_I\sigma_D\,\mu_F\,z\,\frac{G'(z)}{G(z)},
        \qquad
        c_1^\rho c_2 - c_0^\rho = (c_1c_2 - c_0) + \psi_\rho(z)\,\mu_F\,z\,\sigma_I\big(\theta(\sigma_I - w) - (1-\theta)\sigma_D\big),
    \]
    where $c_1c_2 - c_0$ is the uncontrolled expression of Lemma~\ref{lem:linearization}, given by~\eqref{eq:RH-gap}.
    Consequently:
    \begin{enumerate}
        \item[\textup{(i)}] if $G'(z) > 0$, then $q^*$ is unstable;
        \item[\textup{(ii)}] if $G'(z) < 0$, then $q^*$ is locally asymptotically stable if and only if $c_2 > 0$ and $c_1^\rho c_2 - c_0^\rho > 0$;
        \item[\textup{(iii)}] if the rate condition~\eqref{eq:rate-condition} holds and $|\psi_\rho(z)|\,\mu_F \leq c\,\alpha(z)$, then $c_2 > 0$ and $c_1^\rho c_2 - c_0^\rho > 0$; in particular, $q^*$ is locally asymptotically stable if, in addition, $G'(z) < 0$.
    \end{enumerate}
\end{lemma}

The sign of $\psi_\rho\big(\theta(\sigma_I - w) - (1-\theta)\sigma_D\big)$ has an important interpretation: cutting launches stabilizes when the intacts dominate fragment production ($\theta$ near $1$), and can destabilize when the derelicts dominate, because a derelict population responds to a cut in launches only after the intacts have aged into it.
Condition~(iii) bounds the steepness of the throttle.
It is sufficient, not necessary. %

Corollary~\ref{cor:exp} exhibits \emph{a} fold-free throttle, but a whole family of throttles removes the fold, and they differ in how many launches they block.
We now identify the \emph{extremal clip}, which is the pointwise-largest throttle which removes the fold.
Recall that $\rho \le 1$.
By Theorem~\ref{thm:clip-fluid}, for every launch rate $\lambda \leq \mu_I$ the non-trivial fixed points of the controlled fluid limit are the solutions $q^*_F > 0$ of
\begin{equation}\label{eq:controlled-fp}
    c\,\kfrags\,\lambda\,G(q^*_F) = \mu_F, \qquad G = \rho H.
\end{equation}
We call an admissible throttle \emph{fold-free} if, for every $\lambda > 0$, the set of solutions of~\eqref{eq:controlled-fp} is connected (a single point or, degenerately, an interval of fixed points) and \emph{monostable} if it has at most one point for every $\lambda > 0$.
Under the rate condition~\eqref{eq:rate-condition}, the fixed points of a monostable throttle are locally asymptotically stable whenever the steepness condition of Lemma~\ref{lem:clip-linearization}(iii) holds; for the extremal clip this is part of Theorem~\ref{prop:heaviest}.

\begin{theorem}
\label{prop:heaviest}
Let $\bar H(q) := \inf_{s \le q} H(s)$ and $\rho^\star(q) := \bar H(q)/H(q)$.
\begin{enumerate}
    \item[\textup{(i)}] An admissible throttle is fold-free if and only if $G = \rho H$ is non-increasing, and monostable if and only if $G$ is strictly decreasing.
    In particular, a throttle is fold-free if and only if it is single-branch.
    \item[\textup{(ii)}] $\rho^\star$ is admissible and fold-free, and every fold-free throttle satisfies $\rho \leq \rho^\star$ pointwise: $\rho^\star$ is the largest fold-free throttle, and at every fragment level it blocks the fewest launches.
    \item[\textup{(iii)}] In the two-fold regime, for the extremal clip $\rho = \rho^\star$, $G = \bar H$: there is no non-trivial fixed point for $\lambda \leq \lambda_\crit$, exactly one for $\lambda \in (\lambda_\crit, \mu_I]$ with $\lambda \neq \lambda_\flat$, where $\lambda_\flat := \mu_F/(c\,\kfrags\,H(q_{\min})) = \lambda_+$, and at $\lambda = \lambda_\flat$ the non-trivial fixed points are those with $q^*_F$ in the window $[q_{\min}, q^\ddagger]$, where $q^\ddagger$ is the unique point $> q_{\max}$ with $H(q^\ddagger) = H(q_{\min})$.
    If the rate condition~\eqref{eq:rate-condition} holds, every non-trivial fixed point with $\lambda \neq \lambda_\flat$ is locally asymptotically stable: the controlled fluid limit has a \emph{single stable branch}.
    \item[\textup{(iv)}] If $H$ is non-increasing, then $\rho^\star \equiv 1$.
    Otherwise (in particular in the one- and two-fold regimes), $\rho^\star$ is fold-free but not monostable, there is no largest monostable throttle, and $\rho^\star$ is the pointwise supremum of the monostable throttles: every monostable throttle satisfies $\rho \leq \rho^\star$, and the tilted clips $\rho^\star(q)e^{-\eta q}$, $\eta > 0$, are monostable and increase to $\rho^\star$ as $\eta \downarrow 0$.
\end{enumerate}
In the two-fold regime, let $q_{\min} < q_{\max}$ be the interior minimizer and maximizer of $H$, and $q^\ddagger$ the unique point $>q_{\max}$ with $H(q^\ddagger) = H(q_{\min})$.
Then $\rho^\star$ is the \emph{plateau clip}
\[
\rho^\star(q_F) =
\begin{cases}
1, & q_F \notin [q_{\min}, q^\ddagger] \quad(\text{no throttle}),\\[3pt]
\dfrac{H(q_{\min})}{H(q_F)}, & q_F \in [q_{\min}, q^\ddagger] \quad(\text{clip}),
\end{cases}
\]
which throttles only inside the fold window and by at most the factor $H(q_{\min})/H(q_{\max})$, attained at the peak $q_{\max}$.

In the one-fold regime ($0 < \Qsat < \underline{Q}$), $H'(0) > 0$, so $H$ rises from the origin to one peak $q_{\max}$ before decaying and has no interior minimizer.
The minimizer then degenerates to $q_{\min} = 0$, the clip window becomes $[0, q^\ddagger]$ with $H(q^\ddagger) = H(0)$, and $\rho^\star = H(0)/H(q_F)$.
Thus the running-infimum characterization $G^\star = \bar H$ holds in either regime.
\end{theorem}

\begin{remark}
\label{rem:frontier}
Under the condition~\eqref{eq:rate-condition}, for every launch rate $\lambda \in (\lambda_\crit, \mu_I]$ other than the critical rate $\lambda_\flat$, the controlled fluid limit under $\rho^\star$ has exactly one non-trivial fixed point, and it is locally asymptotically stable.
The throttle $\rho^\star$ is the boundary of the fold-free region: $G^\star = \bar H$ is only \emph{weakly} monotone.
It takes the constant value $H(q_{\min})$ across the fold window, meeting the rising $H$ at $q_{\min}$ and $q^\ddagger$.
Thus, the non-trivial fixed point is unique at every launch rate except the single critical rate
\[
    \lambda_\flat = \frac{\mu_F}{c\,\kfrags\,H(q_{\min})} = \lambda_+,
\]
at which the fixed points fill the whole window, a degenerate cusp.
Operationally, $\rho^\star$ removes the high-debris branch while retaining the critical slowing.
\end{remark}

\begin{remark}
\label{rem:high-branch}
Monotonicity of $G$ in either direction would make every level set $\{G = L\}$ connected, so one may ask why Theorem~\ref{prop:heaviest}(i) admits only non-increasing $G$, and whether a throttle could remove the fold by keeping the high-debris branch instead of the low-debris one.
Under admissibility neither is possible.
\begin{enumerate}
    \item[\textup{(a)}] \emph{$G$ cannot rise for good.}
    Since $\rho \le 1$, we have $G = \rho H \le H$, and $H(q) \to 0$ as $q \to \infty$; so every admissible $G$ is positive and tends to $0$, and no admissible $G$ is non-decreasing.
    If $G$ rises on some interval, its eventual decay produces a level which $G$ crosses on both sides of the rise, and the throttle is not fold-free (proof of Theorem~\ref{prop:heaviest}(i)).
    Moreover, a fixed point at which $\rho$ is differentiable and $G' > 0$ is unstable by Lemma~\ref{lem:G-bifurcation}(iii), so a rising stretch of $G$ never carries a stable branch.
    \item[\textup{(b)}] \emph{Keeping the high-debris branch requires boosting launches.}
    The mirror image of $\rho^\star$ is built from the smallest non-increasing majorant of $H$,
    \[
        \hat H(q) := \sup_{s \ge q} H(s), \qquad \hat\rho := \hat H / H .
    \]
    Since $\hat H$ is non-increasing, $\hat\rho$ has no fold, and it keeps the high-debris branch.
    Indeed $\hat H \ge H(q_{\max})$ on $[0, q_{\max}]$ and $\hat H = H$ on $[q_{\max}, \infty)$, so for every $\lambda > \lambda_- = \cM(q_{\max})$ the only non-trivial fixed point is the high-debris one, with $q^*_F > q_{\max}$, and the clean state is unstable, because $\cM_{\hat\rho}(0) = \mu_F/(c\,\kfrags\,\hat H(0)) \le \lambda_-$ (Lemma~\ref{lem:G-bifurcation}(i), whose proof does not use $\rho \le 1$).
    But $\hat\rho \ge 1$, with strict inequality wherever $H < \hat H$---across the fold window, and down to the origin if $H(0) < H(q_{\max})$---so $\hat\rho$ \emph{boosts} launches instead of blocking them, and it violates admissibility~(i).
    \item[\textup{(c)}] \emph{Every admissible fold-free throttle keeps the low-debris branch.}
    If $\rho$ is fold-free, then $G \le \bar H$ by Theorem~\ref{prop:heaviest}(ii).
    Fix $\lambda$ and let $L := \mu_F/(c\,\kfrags\,\lambda)$.
    If $L > H(0) = \bar H(0)$, then $G \le \bar H < L$, and neither $\rho$ nor $\rho^\star$ has a non-trivial fixed point.
    Otherwise, if $G(z) = L$, then $\bar H(z) \ge L$, so $z \le \hat z_L := \max\{q \ge 0 : \bar H(q) \ge L\}$, and $\hat z_L$, when positive, is the largest non-trivial fixed point under $\rho^\star$.
    Thus, at every launch rate, every non-trivial fixed point of every admissible fold-free throttle lies at or below the largest one under $\rho^\star$: among the fold-free throttles, $\rho^\star$ holds the most debris in equilibrium, as befits the largest.
    In particular, for every $\lambda < \lambda_+$ in the two-fold regime, and every $\lambda < \lambda_\crit$ in the one-fold regime---ranges which contain the bistable windows of Theorem~\ref{thm:bifurcation-structure}---the non-trivial fixed points of $\rho^\star$ have $q^*_F < q_{\min}$ or do not exist, and so the same holds for every admissible fold-free throttle.
\end{enumerate}
Admission control can therefore only remove the high-debris fixed point, leaving the system on the low-debris branch or at a trivial fixed point.
\end{remark}

\begin{figure}[htbp]
\centering
\includegraphics[width=\textwidth]{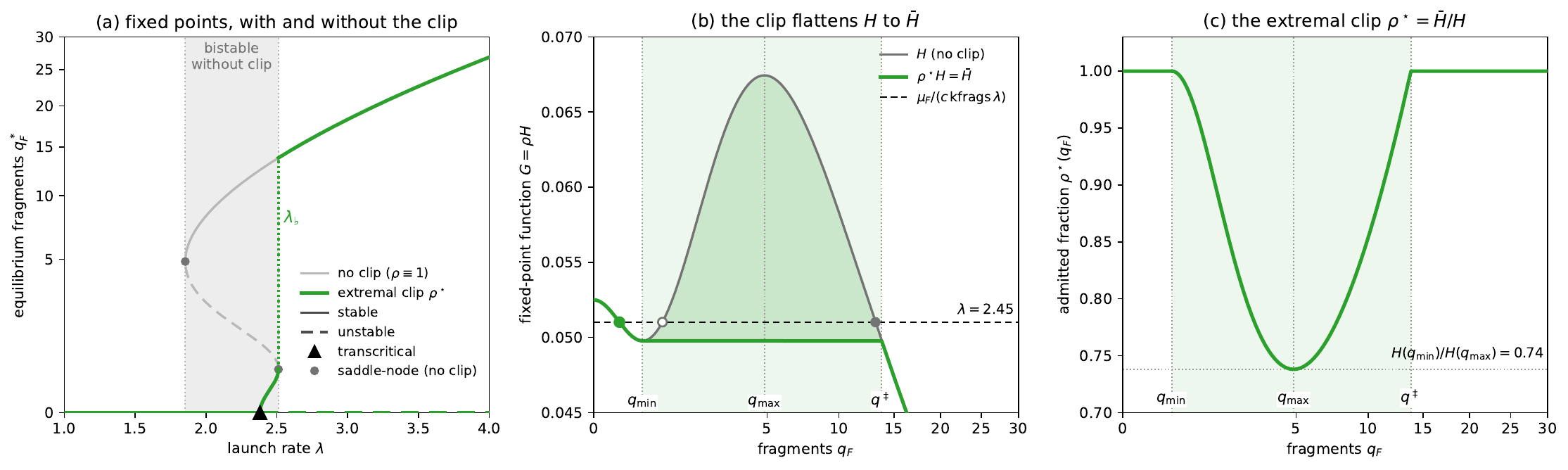}
\caption{\textbf{(a)}~Without the clip (gray) the fixed points form the S-curve of Figure~\ref{fig:fragment-bistability} panel (b). Under $\rho^\star$ (green) the bistability is removed. \textbf{(b)}~A comparison of $\bar{H}$ against $H$. \textbf{(c)}~$\rho^\star$ as a function of the fragment count $q_F$.}
\label{fig:clip-compare}
\end{figure}

Figure~\ref{fig:clip-compare} contrasts the uncontrolled system with that of the extremal clip.

\begin{theorem}\label{thm:tilted-clip}
Fix $\eta > 0$ and take the tilted clip $\rho_\eta(q_F) = \frac{\bar H(q_F)}{H(q_F)}\,e^{-\eta q_F}$, so that the controlled fixed-point function is
\[
    G_\eta(q_F) := \rho_\eta(q_F)\,H(q_F) = \bar H(q_F)\,e^{-\eta q_F}.
\]
Then $\rho_\eta$ is admissible, and, for launch rates $\lambda \leq \mu_I$:
\begin{enumerate}
\item[\textup{(i)}] for $\lambda \le \lambda_\crit$ there is no non-trivial fixed point, and the clean state $q_F = 0$ is locally asymptotically stable for $\lambda < \lambda_\crit$; for $\lambda > \lambda_\crit$ the clean state is unstable and there is exactly one non-trivial fixed point $q_F^*(\lambda) > 0$;
\item[\textup{(ii)}] on $\lambda > \lambda_\crit$ the branch $\lambda \mapsto q_F^*(\lambda)$ is continuous and \emph{strictly} increasing, with no saddle-node; it emerges continuously from $0$ at $\lambda_\crit$. In particular the plateau and jump of the extremal clip at $\lambda_\flat$ (Remark~\ref{rem:frontier}) are removed;
\item[\textup{(iii)}] at every non-trivial fixed point at which $\rho_\eta$ is differentiable,
\[
    c_0^{\rho_\eta} \;\ge\; \eta\,\sigma_I\sigma_D\,\mu_F\,q_F^* \;>\; 0,
\]
so the branch is uniformly bounded away from a saddle-node; and if the rate condition~\eqref{eq:rate-condition} holds and $|\psi_{\rho_\eta}(q^*_F)|\,\mu_F \leq c\,\alpha(q^*_F)$, the fixed point is locally asymptotically stable.
\end{enumerate}
The fold and plateau statements use no assumption on the shape of $H$; in particular they hold in the bistable regime.
\end{theorem}

\begin{figure}[htbp]
\centering
\includegraphics[width=\textwidth]{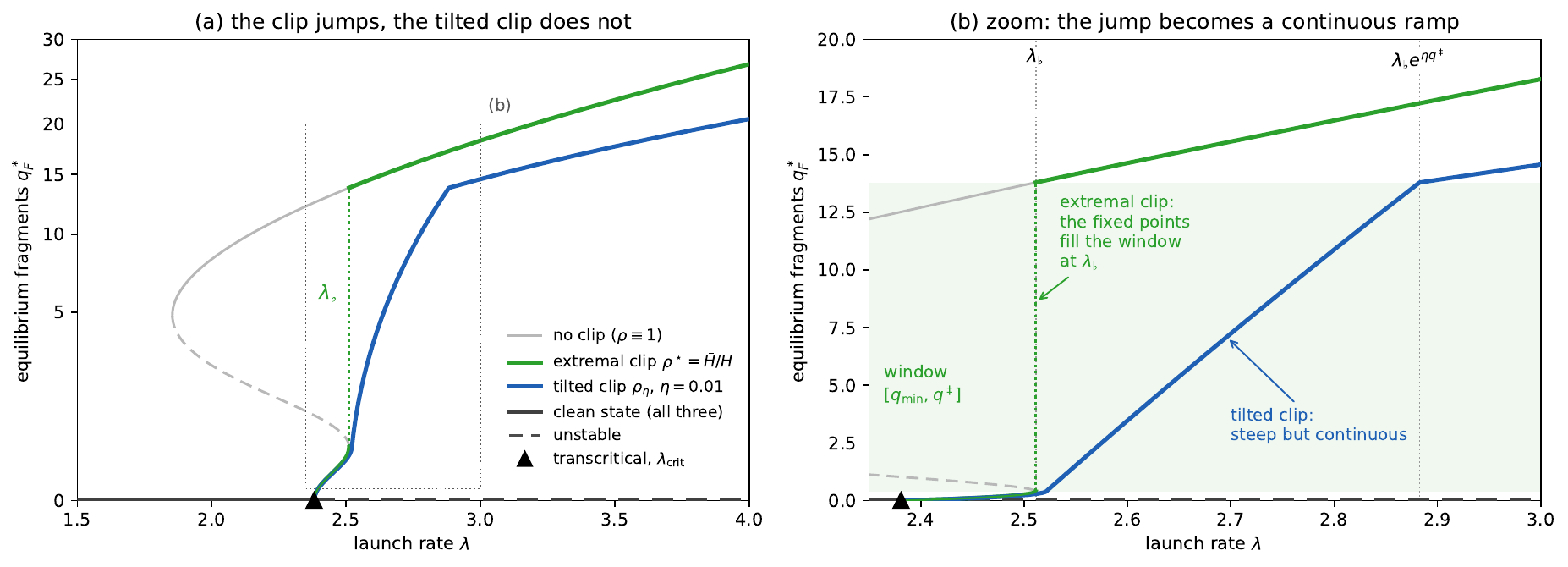}
\caption{The extremal clip's branch jumps; the tilted clip's is continuous. Two-fold regime, with the parameters of Figure~\ref{fig:clip-compare}. \textbf{(a)}~Under the extremal clip $\rho^\star$ (green) the branch is non-decreasing, but $q_F^*$ jumps at $\lambda_\flat$. The tilted clip $\rho_\eta = (\bar H/H)\,e^{-\eta q_F}$ of Theorem~\ref{thm:tilted-clip}, with $\eta = 0.01$ (blue), has a strictly increasing, continuous branch and holds less debris. The uncontrolled S-curve is shown in gray. \textbf{(b)}~Zooming in on the transition, under the tilt the jump becomes a continuous ramp. Apart from at $\lambda_\flat$ for extremal clip, every equilibrium point is stable.}
\label{fig:clip-continuous}
\end{figure}

Figure~\ref{fig:clip-continuous} shows this residual jump and a continuous alternative.

\subsection{Slow Passage}
The clip also reshapes the \emph{dynamic} transition.
Without launch control, under the slow ramp~\eqref{eq:ramp}, a growing constellation crosses the critical rate $\lambda_\flat = \lambda_+$ through the delayed capture fold of Theorem~\ref{thm:slow-passage}, overshooting by $O(\varepsilon^{2/3})$ as in~\eqref{eq:sp-law}, that is, jumping a time of order $\varepsilon^{-1/3}$ after the crossing, as in~\eqref{eq:sp-law-time}.
Replacing $H$ by the clipped fixed-point function changes the \emph{type} of the crossing, and with it the delay law.

\begin{theorem}\label{prop:clip-slow-passage}
Assume the two-fold regime, $\lambda_\flat < \mu_I$, and the rate condition~\eqref{eq:rate-condition}.
Let $\lambda_0 \in (\max\{\lambda_-, \lambda_\crit\}, \lambda_\flat)$, and let $t_\flat := (\lambda_\flat - \lambda_0)/\varepsilon$ be the time at which the ramp~\eqref{eq:ramp} crosses $\lambda_\flat$, so that $\lambda(t) - \lambda_\flat = \varepsilon(t - t_\flat)$.
\begin{enumerate}
\item[\textup{(i)}] \emph{Extremal clip.} Assume in addition that $|H'(z)/H(z)|\,\mu_F \leq c\,\alpha(z)$ for every $z \in [q_{\min}, q^\ddagger]$, and let $q$ be the solution of~\eqref{eq:cFL-1}--\eqref{eq:cFL-3} with $\rho = \rho^\star$, with $\lambda$ replaced by the ramp~\eqref{eq:ramp}, and with $q(0) = q^{\mathrm{lo}}(\lambda_0)$.
At $\lambda = \lambda_\flat$, for every $z \in [q_{\min}, q^\ddagger]$ there is a fixed point with $q^*_F = z$; for $z \in [q_{\min}, q^\ddagger)$ its Jacobian has a simple zero eigenvalue, with eigenvector tangent to the curve of these fixed points, and two eigenvalues with negative real part.
In particular, the window is not a saddle-node, and no fixed point is destroyed as $\lambda$ passes $\lambda_\flat$.
Moreover, as $\varepsilon \to 0$:
\begin{enumerate}
    \item[\textup{(a)}] the first time $t_0$ at which $q_F(t_0) = q_{\min}$ satisfies
    \[
        \lambda(t_0) - \lambda_\flat = \Omega_1\,(\beta_\lambda\beta_2)^{-1/3}\,\varepsilon^{2/3}\,\big(1 + o(1)\big),
    \]
    where $\Omega_1$ is the smallest positive zero of $s \mapsto \mathrm{Ai}'(-s)$ and $\beta_\lambda, \beta_2$ are the coefficients~\eqref{eq:sp-coefficients} of the fold of the uncontrolled system at $q_{\min}$;
    equivalently, $q_F$ reaches $q_{\min}$ a time $t_0 - t_\flat = \Omega_1\,(\beta_\lambda\beta_2)^{-1/3}\,\varepsilon^{-1/3}\,(1 + o(1))$ after the launch rate crosses $\lambda_\flat$;
    \item[\textup{(b)}] for every $z_b \in (q_{\min}, q^\ddagger)$, the first time $t_b$ at which $q_F(t_b) = z_b$ satisfies
    \begin{equation}\label{eq:clip-delay}
        \lambda(t_b) - \lambda_\flat = \sqrt{\frac{2\lambda_\flat}{\mu_F}\,\cJ(z_b)}\;\;\varepsilon^{1/2} + O(\varepsilon^{5/6}),
        \qquad
        \cJ(z_b) := \int_{q_{\min}}^{z_b} \frac{c_1^{\rho^\star}(z)}{\sigma_I(z)\sigma_D(z)}\,\frac{dz}{z},
    \end{equation}
    where $c_1^{\rho^\star}(z)$ is the coefficient of Lemma~\ref{lem:clip-linearization} at the fixed point with $q^*_F = z$ and $\lambda = \lambda_\flat$;
    equivalently, $q_F$ reaches $z_b$ a time
    \begin{equation}\label{eq:clip-delay-time}
        t_b - t_\flat = \sqrt{\frac{2\lambda_\flat}{\mu_F}\,\cJ(z_b)}\;\;\varepsilon^{-1/2} + O(\varepsilon^{-1/6})
    \end{equation}
    after the launch rate crosses $\lambda_\flat$, so that the solution transits the window from $q_{\min}$ to $z_b$ in a time of order $\varepsilon^{-1/2}$.
\end{enumerate}
We write $\cJ_{\mathrm{w}} := \cJ(q^\ddagger) = \lim_{z_b \uparrow q^\ddagger}\cJ(z_b)$ for the width of the whole window.
\item[\textup{(ii)}] \emph{Tilted clip.} Let $\eta > 0$, and let $\rho_\eta = (\bar H/H)\,e^{-\eta q_F}$ be the tilted clip of Theorem~\ref{thm:tilted-clip}, whose branch $\lambda \mapsto q^*(\lambda)$ of non-trivial fixed points is continuous and strictly increasing on $\lambda > \lambda_\crit$.
On the window it is
\[
    q_F^*(\lambda) = \eta^{-1}\ln\!\frac{\lambda}{\lambda_\flat}, \qquad \lambda_\flat e^{\eta q_{\min}} \le \lambda \le \lambda_\flat e^{\eta q^\ddagger},
\]
and $c_0^{\rho_\eta} \geq \eta\,\sigma_I\sigma_D\,\mu_F\,q_{\min} > 0$ at each of its fixed points with $\lambda$ in the open interval $\big(\lambda_\flat e^{\eta q_{\min}}, \lambda_\flat e^{\eta q^\ddagger}\big)$; so the window is replaced by a strictly increasing branch, spread over a launch-rate interval of width $\lambda_\flat\big(e^{\eta q^\ddagger} - e^{\eta q_{\min}}\big)$.
Let $[\lambda_a, \lambda_b]$ be a compact subinterval of $\big(\lambda_\flat e^{\eta q_{\min}}, \lambda_\flat e^{\eta q^\ddagger}\big)$ on which these fixed points are locally asymptotically stable (for instance, under the condition of Lemma~\ref{lem:clip-linearization}(iii)), and let $q$ be the solution of~\eqref{eq:cFL-1}--\eqref{eq:cFL-3} with $\rho = \rho_\eta$, with $\lambda$ replaced by $\lambda(t) = \lambda_a + \varepsilon t$, and with $q(0) = q^*(\lambda_a)$, for $t \in [0, (\lambda_b - \lambda_a)/\varepsilon]$.
Then, as $\varepsilon \to 0$,
\[
    \sup\big\{|q(t) - q^*(\lambda(t))| : \lambda(t) \in [\lambda_a, \lambda_b]\big\} = O(\varepsilon);
\]
there is no bifurcation delay, and the fragment count lags the moving fixed point by a time which is bounded uniformly in $\varepsilon$.
\end{enumerate}
The constants depend on $z_b$, $\lambda_0$, $\eta$, $\lambda_a$, $\lambda_b$ and the parameters of the model, but not on $\varepsilon$.
\end{theorem}

The factor $\sigma_I\sigma_D/c_1^{\rho^\star}$ in~(i) accounts for the finite response time of the intact and derelict populations.
It equals one when these populations equilibrate instantaneously, in which case $\cJ_{\mathrm{w}} = \ln(q^\ddagger/q_{\min})$.

\begin{remark}\label{rem:hierarchy}
The transition type, and hence the overshoot, is fixed by the local shape of the controlled fixed-point function at $\lambda_\flat$:
\[
    \underbrace{\;\varepsilon^{2/3}\;}_{\substack{\text{uncontrolled fold}\\[1pt] \text{(quadratic tangency, }H''(q_{\min}))}}
    \;\longrightarrow\;
    \underbrace{\;\varepsilon^{1/2}\;}_{\substack{\text{extremal clip}\\[1pt] \text{(flat window, }\cJ_{\mathrm{w}})}}
    \;\longrightarrow\;
    \underbrace{\;\varepsilon\;}_{\substack{\text{tilted clip}\\[1pt] \text{(hyperbolic, slope }\eta)}} .
\]
Flattening $H$ across the window removes the saddle-node bottleneck: the delay exponent drops from $\tfrac23$ to $\tfrac12$ and its prefactor no longer measures the local curvature but the global width $\cJ_{\mathrm{w}}$ of the window; tilting the clip restores strict hyperbolicity and ordinary adiabatic ($O(\varepsilon)$) tracking.
Measured in time rather than in launch rate, each overshoot is divided by the ramp rate $\varepsilon$, and the hierarchy reads $\varepsilon^{-1/3} \to \varepsilon^{-1/2} \to O(1)$: the uncontrolled system jumps a time of order $\varepsilon^{-1/3}$ after the crossing~\eqref{eq:sp-law-time}, the extremal clip takes a time of order $\varepsilon^{-1/2}$ to cross its window~\eqref{eq:clip-delay-time}, and the tilted clip lags its branch by a time bounded uniformly in $\varepsilon$.
Two consequences are worth noting.
First, because $\varepsilon^{1/2} \gg \varepsilon^{2/3}$ as $\varepsilon\to0$, the fold-free extremal clip transits over a \emph{wider} launch-rate window than the uncontrolled fold delays---the neutral plateau exerts no restoring force, so a larger excursion is needed to drive the crossing---but the crossing is now smooth and reversible: the $O(1)$ static hysteresis $[\lambda_-,\lambda_+]$ is gone, replaced by a \emph{dynamic} loop of width $O(\varepsilon^{1/2})$ about $\lambda_\flat$ that closes as $\varepsilon\to0$.
Second, the two clipped laws meet at the borderline tilt $\eta \sim \varepsilon^{1/2}$: for $\eta \gg \varepsilon^{1/2}$ the branch is effectively hyperbolic and the passage adiabatic, while for $\eta \lesssim \varepsilon^{1/2}$ the tilt is invisible to the ramp and the degenerate law~\eqref{eq:clip-delay} prevails.
\end{remark}

\subsection{Net Useful Throughput}

In this subsection we discuss the \emph{net useful throughput}, the rate of successful mission completions (no catastrophic collisions, and successful post-mission disposal).
Note that this quantity is given by $P\mu_Iq^*_I$.

\begin{proposition}\label{prop:throughput}
For any admissible throttle $\rho$ and launch rate $\lambda \leq \mu_I$, at a non-trivial fixed point $q^*$ of the controlled fluid limit the admitted launch rate and the net useful throughput $T_{\mathrm{net}} := P\mu_I q_I^*$ depend on the policy only through $q_F^*$:
\[
    \lambda\,\rho(q_F^*) = \cM(q^*_F) = \frac{\mu_F}{c\,\kfrags\,H(q_F^*)}, \quad
    T_{\mathrm{net}} = \Phi(q_F^*), \quad
    \Phi(q) := \frac{P\,\mu_I\,\mu_F}{c\,\kfrags\,\big(\alpha(q) + (1-P)\mu_I/(\mu_D + cq)\big)} .
\]
Moreover $\Phi(q) \to L := \dfrac{P\,\mu_I\,\mu_F}{c\,\kfrags}$ as $q \to \infty$, and $q^*_I \to \mu_F/(c\,\kfrags)$ along any family of non-trivial fixed points with $q^*_F \to \infty$.
Suppose that $\mu_F < c\,\kfrags$.
Then along the fixed-point relations $\lambda = \cM_\rho(q_F)$ of the extremal and tilted clips, $q_F^* \sim (\kfrags/\mu_F)\,\lambda$ under the extremal clip and $q_F^* \sim \eta^{-1}\log\lambda$ under the tilted clip as $\lambda \to \infty$; for large $\lambda$ these fixed points have $q^*_I < 1$, so they are fixed points of the reflected controlled fluid limit of Theorem~\ref{thm:clip-fluid}, at which the formulas above continue to hold; and both deliver $T_{\mathrm{net}} \to L$, the \emph{same} limiting net useful throughput, even though their equilibrium debris diverges at very different rates.
\end{proposition}

\begin{figure}[htbp]
\centering
\includegraphics[width=\textwidth]{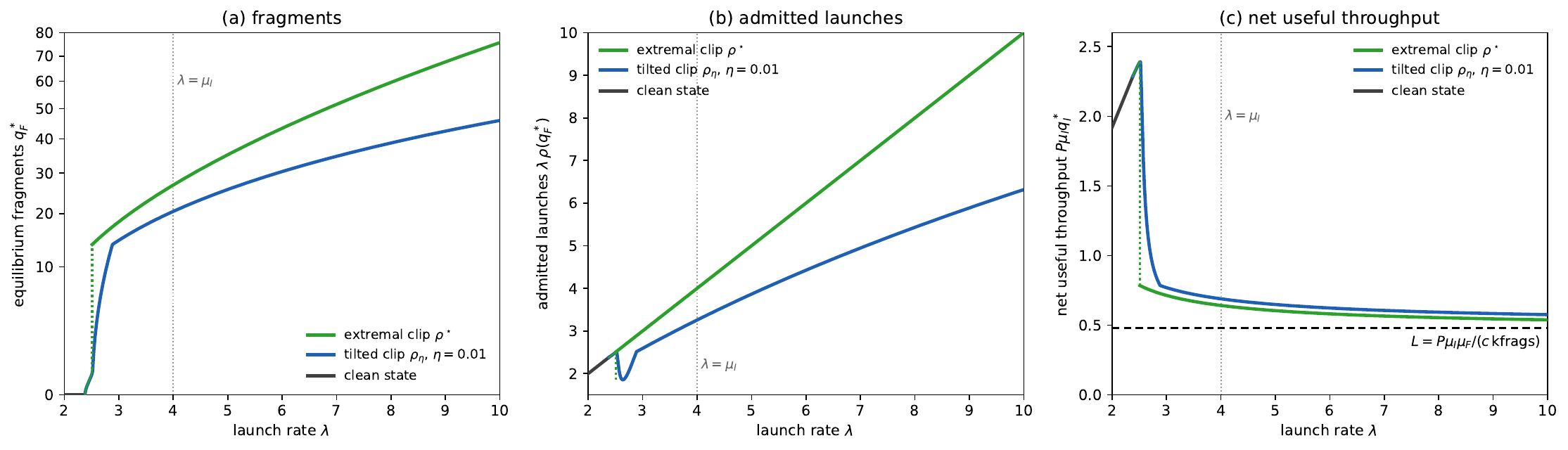}
\caption{Extremal (green) versus tilted (blue, $\eta = 0.01$) clip as functions of the launch rate, with the parameters of Figure~\ref{fig:clip-continuous}. Below $\lambda_\crit \approx 2.38$ both hold the clean state (dark gray). To the right of the dotted line $\lambda = \mu_I$ the non-trivial fixed points have $q^*_I < 1$, so they are fixed points of the reflected controlled fluid limit (Proposition~\ref{prop:throughput}). \textbf{(a)}~Equilibrium fragments (square-root scale): the extremal clip jumps at $\lambda_\flat \approx 2.51$; the tilted clip ramps continuously and holds debris strictly lower, $q_F^* \approx 46$ against $76$ at $\lambda = 10$. \textbf{(b)}~Admitted launches $\lambda\rho(q_F^*) = \cM(q_F^*)$: the extremal clip admits every launch on its branch ($\rho^\star = 1$ off the window), while the tilt blocks a growing fraction, and across the window its admitted rate dips to $\lambda_- = \cM(q_{\max}) \approx 1.85$. \textbf{(c)}~Net useful throughput $P\mu_I q^*_I = \Phi(q_F^*)$: at $\lambda_\flat$ the extremal clip's jump into the collision-dominated regime cuts it from $\Phi(q_{\min}) \approx 2.39$ to $\Phi(q^\ddagger) \approx 0.79$; under the tilt it varies continuously and exceeds the extremal clip's at every $\lambda > \lambda_\flat$; and both tend to the policy-independent limit $L = 0.48$ (dashed) of Proposition~\ref{prop:throughput}, reaching $0.54$ and $0.58$ at $\lambda = 10$.}
\label{fig:clip-throughput}
\end{figure}

\begin{remark}
Proposition~\ref{prop:throughput} sharpens the trade-off.
Raw admitted launches favor the extremal clip, which never throttles on its branch; but those launches are largely wasted, since the extremal clip parks the system deep in the collision-dominated regime ($\alpha c q_F \gg \mu_I$), where most satellites are destroyed before disposal.
Measured in successful disposals the tilt does better: in the example it exceeds the extremal clip at every launch rate above $\lambda_\flat$ (Figure~\ref{fig:clip-throughput}(c)), while holding debris exponentially lower ($\log\lambda$ versus $\lambda$).
The dial $\eta$ therefore trades gross launch volume---not useful output---for a jump-free, lower-debris environment, and at high debris levels every admission policy delivers the same net useful throughput $L$.
\end{remark}

Proposition~\ref{prop:throughput} reduces the comparison of admission policies at equilibrium to a comparison of fragment levels: whatever the policy, an equilibrium with $q^*_F = z > 0$ admits launches at the rate $\cM(z)$ and delivers $T_{\mathrm{net}} = \Phi(z)$.
What remains is to choose the fragment level, and a rule which holds it.
We now show that a \emph{bang-bang} rule, which admits every offered launch while the fragment count is below a threshold and none while it is above, is optimal, and that the optimal threshold is an explicit fragment level $q^\circ$ below the fold window.

Write
\[
    \varphi(q) := \sigma_I(q)H(q) = \alpha(q) + \frac{(1-P)\mu_I}{\sigma_D(q)},
    \qquad\text{so that}\qquad
    \Phi(q) = \frac{P\mu_I\,\mu_F}{c\,\kfrags\,\varphi(q)} = \frac{P\mu_I\,\cM(q)}{\sigma_I(q)}.
\]
At an equilibrium with $q^*_F = q$ we have $q^*_D/q^*_I = (1-P)\mu_I/\sigma_D(q)$, so $\varphi(q) = \alpha(q) + q^*_D/q^*_I$ is the collision load which each intact brings to the shell, through its own collisions and through the derelicts it leaves behind.
The fragment balance fixes the total collision load $\alpha(q^*_F)q^*_I + q^*_D$ at $\mu_F/(c\,\kfrags)$, so the useful throughput is largest where $\varphi$ is smallest.
More fragments make the intacts collide more often, through $\alpha$, but also remove the derelicts faster, through $\sigma_D$; the optimal fragment level balances the two.
Let
\begin{equation}\label{eq:q-circ}
    r := \sqrt{\frac{(1-P)\mu_I}{(1-\alpha_0)\,c\,\Qsat}},
    \qquad
    q^\circ := \frac{r\Qsat - \mu_D/c}{1-r}.
\end{equation}

To compare policies, we describe an equilibrium by its state and its admitted launch rate.
For $J \ge 0$, let $b^J$ be the right-hand side of~\eqref{eq:FL-1}--\eqref{eq:FL-3} with $\lambda$ replaced by $J$.
An \emph{equilibrium under admission at the offered rate $\lambda$} is a pair $(q^*, J)$ with $q^* \in \mathcal S$, $J \in [0,\lambda]$ and $b^J(q^*) = 0$.
Every fixed point of a controlled fluid limit in this section is of this form: under debris feedback with $J = \lambda\rho(q^*_F)$, under occupancy feedback with $\bar\lambda \le \lambda$ and $J = \bar\lambda(q^*_I)$, and under the bang-bang rule below.
At a \emph{clean} equilibrium, $q^*_F = 0$, the fragment row of the Jacobian is $\big(0,\ 0,\ c\,\kfrags\,J H(0) - \mu_F\big)$ under every differentiable admission rule, since the rule enters only the intact equation; so a clean equilibrium with $J > \lambda_\crit$ is unstable, and a single collision there starts a cascade.
We call an equilibrium \emph{debris-stable} if $q^*_F > 0$ or $J \le \lambda_\crit$.

For an offered rate $\lambda$ and a debris cap $\bar q \in (0,\infty]$, we ask for the largest net useful throughput at a debris-stable equilibrium,
\begin{equation}\label{eq:T-opt}
    T^{\mathrm{opt}}(\lambda, \bar q) := \sup\big\{ P\mu_I\,q^*_I : (q^*, J) \text{ is a debris-stable equilibrium under admission at } \lambda,\ q^*_F \le \bar q \big\}.
\end{equation}
This is a steady-state criterion: it ranks admission rules by the useful throughput at their equilibria, and does not value admitted launches, transients, or the time taken to reach the equilibrium.

\begin{theorem}\label{thm:bang-bang}
Assume that
\begin{equation}\label{eq:unimodal}
    \frac{\mu_D}{c\,\Qsat} < r < 1,
\end{equation}
and let $\lambda \le \mu_I$.
\begin{enumerate}
    \item[\textup{(i)}] \emph{(The optimal fragment level.)} $q^\circ > 0$; $\Phi$ is strictly increasing on $[0, q^\circ]$ and strictly decreasing on $[q^\circ, \infty)$; and $\cM$ is strictly increasing on $[0, q^\circ]$, with $\cM(0) = \lambda_\crit$. Moreover $H' < 0$ on $[0, q^\circ]$, so~\eqref{eq:unimodal} rules out the one-fold regime, and $q^\circ < q_{\min}$ in the two-fold regime.
    \item[\textup{(ii)}] \emph{(The bound.)} Let $z_\lambda$ be the solution of $\cM(z) = \lambda$ in $(0, q^\circ)$ if $\lambda_\crit < \lambda < \cM(q^\circ)$, and $z_\lambda := q^\circ$ if $\lambda \ge \cM(q^\circ)$; for $z \in (0, q^\circ]$, let
    \[
        T^\star(\lambda, z) := \begin{cases} P\lambda, & \lambda \le \lambda_\crit, \\ \Phi\big(\min\{z_\lambda, z\}\big), & \lambda > \lambda_\crit. \end{cases}
    \]
    Fix a debris cap $\bar q \in (0, \infty]$.
    Every debris-stable equilibrium under admission at the offered rate $\lambda$ with $q^*_F \le \bar q$ has $T_{\mathrm{net}} = P\mu_I q^*_I \le T^\star\big(\lambda, \min\{\bar q, q^\circ\}\big)$.
    \item[\textup{(iii)}] \emph{(The bang-bang rule attains the bound.)} Let $\hat q := \min\{\bar q, q^\circ\}$, and consider the bang-bang rule with threshold $\hat q$, which admits every offered launch while $q_F < \hat q$ and none while $q_F > \hat q$.
    Its fluid limit, described in Theorem~\ref{thm:bang-bang-fluid} below, has exactly one debris-stable equilibrium:
    the clean state, with $J = \lambda$, if $\lambda \le \lambda_\crit$;
    the low-debris fixed point of the uncontrolled fluid limit, with $q^*_F = z_\lambda < \hat q$ and $J = \lambda$, if $\lambda_\crit < \lambda < \cM(\hat q)$;
    and, if $\lambda \ge \cM(\hat q)$, the point $x^*(\hat q)$ with
    \[
        q_F = \hat q, \qquad q_I = \frac{\cM(\hat q)}{\sigma_I(\hat q)}, \qquad q_D = \frac{(1-P)\mu_I\,q_I}{\sigma_D(\hat q)},
    \]
    and $J = \cM(\hat q)$.
    In each case $T_{\mathrm{net}} = T^\star(\lambda, \hat q)$: the bound of~(ii) is attained, at every offered rate $\lambda \le \mu_I$, by a single rule which does not depend on $\lambda$.
\end{enumerate}
\end{theorem}

The proof is in Section~\ref{sec:admission-proofs}.
Three features of Theorem~\ref{thm:bang-bang} deserve comment.
First, at a fixed equilibrium fragment level $q^*_F = \bar q \le q^\circ$, every policy delivers the same useful throughput $\Phi(\bar q)$ (Proposition~\ref{prop:throughput}), and the bang-bang rule with threshold $\bar q$ is the one which holds the fragments at $\bar q$ at every offered rate $\lambda \ge \cM(\bar q)$.
It uses only the fragment count, not the offered rate or the function $\cM$, and at offered rates below $\cM(\bar q)$ its equilibrium blocks nothing.
Second, the cap $\bar q = \infty$ gives the unconstrained optimum: holding the fragments at $q^\circ$ maximizes the net useful throughput over all debris-stable equilibria of all admission rules.
In the running example, $r \approx 0.130$, $q^\circ \approx 0.342$, $\cM(q^\circ) \approx 2.509$ and $\Phi(q^\circ) \approx 2.39$; at $\lambda = \mu_I = 4$, the extremal clip settles at $q^*_F \approx 26.9$ and delivers $T_{\mathrm{net}} \approx 0.64$, about a quarter of the optimum.
Third, since $q^\circ < q_{\min}$, the optimal rule never lets the fragments enter the fold window.
The bang-bang rule has no equilibrium above its threshold at any launch rate, so it removes the bistability, and above $\cM(q^\circ)$---just below $\lambda_+ \approx 2.511$ in the example---it holds the fragments at $q^\circ$ by blocking the fraction $1 - \cM(q^\circ)/\lambda$ of the offered launches.
When $\lambda_\crit < \lambda < \cM(\hat q)$, the rule admits every launch near its equilibrium, which is then the low-debris fixed point of the uncontrolled fluid limit, with $H'(z_\lambda) < 0$ by~(i), and so locally asymptotically stable under the rate condition~\eqref{eq:rate-condition} by Theorem~\ref{thm:bifurcation-structure}.

Theorem~\ref{thm:bang-bang} also answers a per-launch question: what fraction of the offered launches is \emph{useful}, in the sense that the launch is admitted, the satellite completes its mission without a catastrophic collision, and its post-mission disposal succeeds?
At an equilibrium under admission at the offered rate $\lambda > 0$, this \emph{useful fraction} is the ratio $T_{\mathrm{net}}/\lambda$ of the rate of successful mission completions to the rate of offered launches.

\begin{corollary}\label{cor:useful-fraction}
Let $(q^*, J)$ be an equilibrium under admission at the offered rate $\lambda > 0$.
Then
\begin{equation}\label{eq:useful-fraction}
    \frac{T_{\mathrm{net}}}{\lambda} = \frac{J}{\lambda}\cdot\frac{P\mu_I}{\sigma_I(q^*_F)}:
\end{equation}
the useful fraction is the fraction $J/\lambda$ of the offered launches which are admitted, times the probability $P\mu_I/\sigma_I(q^*_F)$ that a satellite which faces the equilibrium collision rate $\alpha(q^*_F)\,c\,q^*_F$ completes its mission, at rate $\mu_I$, before a catastrophic collision, and is then disposed of successfully.
Assume in addition~\eqref{eq:unimodal} and $\lambda \le \mu_I$, fix a debris cap $\bar q \in (0,\infty]$, and let $\hat q := \min\{\bar q, q^\circ\}$.
Then every debris-stable equilibrium under admission at $\lambda$ with $q^*_F \le \bar q$ has
\[
    \frac{T_{\mathrm{net}}}{\lambda} \le \frac{T^\star(\lambda, \hat q)}{\lambda} =
    \begin{cases}
        P, & \lambda \le \lambda_\crit, \\[2pt]
        \dfrac{P\mu_I}{\sigma_I(z_\lambda)}, & \lambda_\crit < \lambda < \cM(\hat q), \\[8pt]
        \dfrac{\Phi(\hat q)}{\lambda}, & \lambda \ge \cM(\hat q),
    \end{cases}
\]
with $T^\star$ and $z_\lambda$ as in Theorem~\ref{thm:bang-bang}(ii), and the bang-bang rule with threshold $\hat q$ attains this bound at every offered rate $\lambda \in (0, \mu_I]$.
In particular, the bang-bang rule with threshold $q^\circ$ maximizes the useful fraction over all debris-stable equilibria of all admission rules, at every offered rate $\lambda \in (0, \mu_I]$ simultaneously.
\end{corollary}

\begin{proof}
The intact component of $b^J(q^*) = 0$ reads $J = \sigma_I(q^*_F)\,q^*_I$, so $T_{\mathrm{net}}/\lambda = P\mu_I q^*_I/\lambda = (J/\lambda)\,P\mu_I/\sigma_I(q^*_F)$, which is~\eqref{eq:useful-fraction}.
Since $\lambda$ is fixed, dividing the bound of Theorem~\ref{thm:bang-bang}(ii) by $\lambda$ gives $T_{\mathrm{net}}/\lambda \le T^\star(\lambda, \hat q)/\lambda$ at every debris-stable equilibrium with $q^*_F \le \bar q$, and Theorem~\ref{thm:bang-bang}(iii) gives equality for the bang-bang rule with threshold $\hat q$.
It remains to evaluate $T^\star(\lambda, \hat q)/\lambda$.
If $\lambda \le \lambda_\crit$, then $T^\star(\lambda, \hat q) = P\lambda$.
Otherwise $T^\star(\lambda, \hat q) = \Phi(\min\{z_\lambda, \hat q\})$.
By Theorem~\ref{thm:bang-bang}(i), $\cM$ is strictly increasing on $[0, q^\circ]$ with $\cM(0) = \lambda_\crit$, so $\lambda_\crit < \cM(\hat q) \le \cM(q^\circ)$.
If $\lambda_\crit < \lambda < \cM(\hat q)$, then $z_\lambda \in (0, q^\circ)$ solves $\cM(z_\lambda) = \lambda < \cM(\hat q)$, so $z_\lambda < \hat q$, and the identity $\Phi = P\mu_I\cM/\sigma_I$ gives $T^\star(\lambda, \hat q) = \Phi(z_\lambda) = P\mu_I\lambda/\sigma_I(z_\lambda)$.
If $\lambda \ge \cM(\hat q)$, then $z_\lambda \ge \hat q$: either $\lambda \ge \cM(q^\circ)$ and $z_\lambda = q^\circ \ge \hat q$, or $z_\lambda \in (0, q^\circ)$ and $\cM(z_\lambda) = \lambda \ge \cM(\hat q)$.
Hence $T^\star(\lambda, \hat q) = \Phi(\hat q)$.
\end{proof}

The factorization~\eqref{eq:useful-fraction} shows what the bang-bang rule trades.
In the running example at $\lambda = \mu_I = 4$, it admits $\cM(q^\circ)/\lambda \approx 63\%$ of the offered launches, of which $P\mu_I/\sigma_I(q^\circ) \approx 95\%$ are useful, so that about $60\%$ of the offered launches are useful; the extremal clip admits every offered launch, but only about $16\%$ of them are useful.
The fraction of the \emph{admitted} launches which are useful, $P\mu_I/\sigma_I(q^*_F)$, is a different objective, and the bang-bang rule does not maximize it.
It is decreasing in $q^*_F$, because $\alpha(q)\,q$ is increasing, so it is largest, equal to $P$, at a clean equilibrium, and a clean equilibrium with $J \le \lambda_\crit$ is debris-stable at every offered rate.
The bang-bang rule accepts a smaller useful fraction of the admitted launches in exchange for admitting more of them, and the product of the two, the useful fraction of the offered launches, is what it maximizes.

The bang-bang rule is not admissible in the sense of Definition~\ref{def:admissible}, because $\ind_{[0,\bar q)}$ is not Lipschitz, and Theorem~\ref{thm:clip-fluid} does not apply to it.
Fix $\bar q > 0$.
In the $n$-th system of Section~\ref{sec:limiting}, the bang-bang rule with threshold $\bar q$ admits an offered launch if and only if $Q_F < n\bar q$; all other transitions are unchanged.
Denote the resulting process $Q^{n,\bar q}$, and let
\[
    U^n(t) := \int_0^t \ind\{Q^{n,\bar q}_F(s) < n\bar q\}\,ds
\]
be the time up to $t$ during which launches are admitted.
The fluid limit is described by the system
\begin{equation}\label{eq:BB}
\begin{aligned}
    q_I(t) &= q_I(0) + \int_0^t \big(\lambda u(s) - \sigma_I(q_F(s))\,q_I(s)\big)\,ds,
    \\
    q_D(t) &= q_D(0) + \int_0^t b_D(q(s))\,ds,
    \qquad
    q_F(t) = q_F(0) + \int_0^t b_F(q(s))\,ds,
\end{aligned}
\end{equation}
with $b_D, b_F$ as in Section~\ref{sec:limiting}, in which the admitted fraction $u : [0,\infty) \to [0,1]$ is measurable, with $u = 1$ almost everywhere on $\{t : q_F(t) < \bar q\}$ and $u = 0$ almost everywhere on $\{t : q_F(t) > \bar q\}$.
A \emph{solution} of~\eqref{eq:BB} is a continuous $q : [0,\infty) \to \mathcal S$ for which such a control $u$ exists; it is a Filippov solution of~\eqref{eq:cFL-1}--\eqref{eq:cFL-3} with $\rho = \ind_{[0,\bar q)}$~\cite{filippov1988}, in which, on the surface $\Sigma := \{q_F = \bar q\}$ where the launch term jumps, the launch rate may be any convex combination of its two one-sided limits.
An \emph{equilibrium} of~\eqref{eq:BB} is a constant solution.
Let
\[
    \xi(q) := \alpha(q_F)\,q_I + q_D, \qquad \bar\xi := \frac{\mu_F}{c\,\kfrags},
\]
so that $b_F(q) = c\,\kfrags\,q_F\big(\xi(q) - \bar\xi\big)$.
On the curve $\Sigma_0 := \{q \in \Sigma : \xi(q) = \bar\xi\}$, on which $b_F = 0$, let
\[
    u_\Sigma(q) := \frac{1}{\lambda}\Big(\sigma_I(\bar q)\,q_I - \frac{b_D(q)}{\alpha(\bar q)}\Big)
\]
be the admitted fraction which keeps $\xi$ constant, and write $\Sigma_{\mathrm{sl}} := \{q \in \Sigma_0 : 0 < u_\Sigma(q) < 1\}$ for the \emph{sliding set} and $\Sigma_\partial := \Sigma_0 \setminus \Sigma_{\mathrm{sl}}$.

\begin{theorem}[Fluid limit of the bang-bang process]\label{thm:bang-bang-fluid}
    Let $\lambda \le \mu_I$ and $\bar q > 0$, and suppose $Q^{n,\bar q}(0)/n \to q(0) \in \mathcal S$. Write $\bar Q^n := Q^{n,\bar q}/n$.
    \begin{enumerate}
        \item[\textup{(i)}] Almost surely, every subsequence of $(\bar Q^n)$ has a further subsequence which converges, uniformly on compacts, to a solution of~\eqref{eq:BB}. In particular,~\eqref{eq:BB} has a solution.
        \item[\textup{(ii)}] There is $T^\circ \in [0,\infty]$ such that all solutions of~\eqref{eq:BB} coincide on $[0, T^\circ] \cap [0,\infty)$, and $T^\circ$ is the first time at which this common solution $q$ lies in $\Sigma_\partial$ ($T^\circ = \infty$ if it never does).
        On $[0, T^\circ]$, $q$ meets $\Sigma \setminus \Sigma_0$ only at isolated times, at which it crosses $\Sigma$, and the control is determined by the state: almost everywhere, $u = 1$ where $q_F < \bar q$, $u = 0$ where $q_F > \bar q$, and $u = u_\Sigma(q)$ where $q \in \Sigma_0$.
        \item[\textup{(iii)}] Almost surely, $\bar Q^n \to q$ and $U^n \to \int_0^\cdot u(s)\,ds$ uniformly on $[0, T]$, for every finite $T \le T^\circ$.
        \item[\textup{(iv)}] Suppose that $\lambda > \cM(\bar q)$. If $q(t_0) \in \Sigma_{\mathrm{sl}}$ for some $t_0 < T^\circ$, then $T^\circ = \infty$ and $q(t) \in \Sigma_{\mathrm{sl}}$ for all $t \ge t_0$: the fragment count stays at $\bar q$, the admitted fraction is $u_\Sigma(q(t))$, and
        \[
            q_I(t) - q_I^*(\bar q) = \big(q_I(t_0) - q_I^*(\bar q)\big)\,e^{-\omega(t - t_0)}, \qquad q_I^*(\bar q) := \frac{\cM(\bar q)}{\sigma_I(\bar q)}, \qquad \omega := \frac{(1-P)\mu_I + \alpha(\bar q)\sigma_D(\bar q)}{\alpha(\bar q)},
        \]
        so that $q(t) \to x^*(\bar q)$ exponentially fast.
    \end{enumerate}
\end{theorem}

The rule reads the fragment count but acts on the intacts, and the fragments respond to the intacts only through collisions.
The switch is therefore two integrations away from the quantity it controls, and the fluid limit slides, not on the surface $\Sigma$ itself, as it would if the rule acted directly on $q_F$, but on the curve $\Sigma_0 \subset \Sigma$ on which the fragment count is momentarily constant.
At finite $n$ the rule chatters, switching launches on and off as $Q_F$ fluctuates around $n\bar q$, and Statement~(iii) says that the fraction of time during which it admits launches converges to the fluid control, which on $\Sigma_{\mathrm{sl}}$ is the equivalent control $u_\Sigma$.

\section{Conclusions and Future Work}
\label{sec:conclusions}
Low earth orbit (LEO) operations have assumed considerable economic and national security importance. However, sustaining operations of the critical infrastructures present and planned in LEO, including SATCOM constellations and datacenters, is subject to considerable risk from space debris. While this is a problem that has been acknowledged for decades now, models of the debris population have focused largely on the collision physics, with launches, disposal, and collision avoidance taken as inputs, whereas economic models of orbit use endogenize launches but represent the physical environment by stylized collision-risk functions (Section~\ref{sec:related-work}). In particular, collision avoidance enters these models, when it is included at all, as a fixed input, although the capacity to screen conjunctions and plan maneuvers is itself an operational constraint which degrades as the debris population grows. How this operational feedback interacts with the collision physics has received comparatively little attention.

Relative to the existing literature, our contributions are: (i) a layered queueing model
in which collision avoidance is a state-dependent, saturating capability that
couples to the debris dynamics, in contrast to the fixed-parameter treatment in
source--sink and simulation models; (ii) a classification of the resulting
bistability in the fluid limit, with different bifurcation mechanisms from
those in loss networks and in physico-economic models of orbit use; (iii) a
characterization of the slow passage through the fold, which quantifies why
degradation of the debris environment is detected late; and (iv) launch-control
policies which depend on fragment counts. A throttle
indexed by the fragment count can eliminate the bistability, whereas an
occupancy cap, which blocks launches as the shell fills, cannot. Among the
throttles indexed by the fragment count, the extremal clip is the least
restrictive one which eliminates the bistability, and a tilted clip also
removes the delay after the tipping point. Under
condition~\eqref{eq:unimodal}, a bang-bang rule, which admits every launch
while the fragment count is below a threshold and none while it is above,
attains the largest net useful throughput at equilibrium at every offered
launch rate.

There are multiple future directions, that we hope the management science and operations research community can and should pay attention to.
In this paper, we take the launch control as given and do not address its implementation, which
would require coordination among nations and commercial operators.
First, how do we solve the collective action problem of coordinating launches across nations? This may present interesting opportunities for both monetary and non-monetary mechanism design. Second, from a methodological angle, understanding the slow passage through the fold for the pre-limit system will require us to move beyond the fluid analysis of the layered queueing network. In particular, it would be useful to understand what happens as both the constellation size and the launch rate ramping simultaneously vary. This is analogous to the interchangeability analyses for G/G/1 queues, and we conjecture will require some new analytical developments. Third, our current model is limited to a single thin shell. Extending this to multiple shells is a significant challenge and will change the equilibrium analyses. This would also change the operational implications.

\section*{Acknowledgments}
AG and HH acknowledge a generous Abhi Deshmukh IE Frontiers Grant, from the Edwardson School of Industrial Engineering.

\section*{AI Statement}
The authors used generative AI tools to assist with proof reading, code development, and some proof checking. The research questions, theoretical ideas, methodology, analysis, and conclusions are the authors' own. All AI-assisted content was reviewed and verified by the authors. The authors assume all responsibility for the content of this paper.

\bibliographystyle{plainurl}
\bibliography{references,related_work}
\appendix
\section{Proofs of the Limiting Dynamics}\label{sec:limiting}

This section proves the strong-approximation theorem (Theorem~\ref{thm:strong-approx}) and the fluid limit (Theorem~\ref{thm:fluid-limit}).
We prove slightly more than Theorem~\ref{thm:fluid-limit}: for every launch rate $\lambda > 0$, the fluid-scaled process converges to a limit in which the intact equation~\eqref{eq:FL-1} is reflected at the capacity $q_I = 1$ (Proposition~\ref{prop:reflected-fluid}), and the reflection is inactive when $\lambda \leq \mu_I$.

\subsection{The $n$-th system}
In the $n$-th system there are $n$ orbital slots, launches arrive at rate $n\lambda$, and $\Qsat^n = \frac{n\mu_s}{\nu} = n\Qsat$ and $\alpha^n(Q) = \alpha(Q/n)$.
The process $Q = Q^n = (Q_I, Q_D, Q_F)$ is a Markov chain on $\{0,\dots,n\}\times\bZ_+^2$ with the following transitions, where $B \sim \Poi(\kfrags)$ is drawn independently at each collision:
\begin{center}
\begin{tabular}{lll}
    event & transition & rate \\ \hline
    launch & $(Q_I + 1,\ Q_D,\ Q_F)$ & $n\lambda\,\ind_{Q_I < n}$ \\
    collision of an intact & $(Q_I - 1,\ Q_D,\ Q_F + B)$ & $\alpha^n(Q_F)\,\tfrac{c}{n}\,Q_F\,Q_I$ \\
    failed disposal & $(Q_I - 1,\ Q_D + 1,\ Q_F)$ & $(1-P)\mu_I\,Q_I$ \\
    successful disposal & $(Q_I - 1,\ Q_D,\ Q_F)$ & $P\mu_I\,Q_I$ \\
    collision of a derelict & $(Q_I,\ Q_D - 1,\ Q_F + B)$ & $\tfrac{c}{n}\,Q_F\,Q_D$ \\
    decay of a derelict & $(Q_I,\ Q_D - 1,\ Q_F)$ & $\mu_D\,Q_D$ \\
    decay of a fragment & $(Q_I,\ Q_D,\ Q_F - 1)$ & $\mu_F\,Q_F$
\end{tabular}
\end{center}
Define the cumulative intensities
\begin{align*}
    a(t) & := \int_{0}^{t} n\lambda\,\ind_{Q_I(s) < n}\,ds,
    &
    d_{IF}(t) & := \int_{0}^{t}\alpha^n(Q_F(s))\,\tfrac{c}{n}\,Q_F(s)\,Q_I(s)\,ds,
    \\
    d_{ID}(t) & := \int_{0}^{t}(1-P)\mu_I\, Q_I(s)\,ds,
    &
    d_{IR}(t) & := \int_{0}^{t} P\mu_I\, Q_I(s)\,ds,
    \\
    d_{DF}(t) & := \int_{0}^{t}\tfrac{c}{n}\,Q_F(s)\,Q_D(s)\,ds,
    &
    d_{DR}(t) & := \int_{0}^{t}\mu_D\, Q_D(s)\,ds,
    \\
    d_F(t) & := \int_{0}^{t} \mu_F\, Q_F(s)\,ds,
\end{align*}
and $d_I := d_{IF} + d_{ID} + d_{IR}$, $d_D := d_{DF} + d_{DR}$.
Let $A, D_{IF}, D_{ID}, D_{IR}, D_{DF}, D_{DR}, D_F$ be independent unit-rate Poisson processes, and let $(B_i)_{i\in\bN}$ be i.i.d.\ $\Poi(\kfrags)$ random variables which are independent of them.
Each type of transition has its own Poisson process, so the chain has the random time-change representation~\cite[Chapter~6, Section~4]{ethierkurtz1986markov} (see also~\cite{mandelbaum1998b})
\begin{equation}\label{eq:time-change}
\begin{aligned}
    Q_I(t) & = Q_I(0) + A(a(t)) - D_{IF}(d_{IF}(t)) - D_{ID}(d_{ID}(t)) - D_{IR}(d_{IR}(t)),
    \\
    Q_D(t) & = Q_D(0) + D_{ID}(d_{ID}(t)) - D_{DF}(d_{DF}(t)) - D_{DR}(d_{DR}(t)),
    \\
    Q_F(t) & = Q_F(0) + \sum_{i=1}^{N_\mathrm{col}(t)}B_i - D_F(d_F(t)),
    \qquad
    N_\mathrm{col}(t) := D_{IF}(d_{IF}(t)) + D_{DF}(d_{DF}(t)).
\end{aligned}
\end{equation}
We construct the systems for all $n$ on one probability space, by solving~\eqref{eq:time-change} with the same Poisson processes and the same sequence $(B_i)$ for every $n$, and with initial conditions $Q^n(0)$ which are deterministic.
We write $\bar Q^n := Q^n/n$, $\mathcal S := [0,1]\times[0,\infty)^2$, $|\cdot|$ for the maximum norm on $\bR^3$, and $\|x\|_t := \sup_{s\leq t}|x(s)|$.

The transition rates are quadratic in the state, so they are neither bounded nor globally Lipschitz.
The following bounds confine the scaled process to a compact set.

\begin{lemma}\label{lem:apriori}
    Suppose $\bar Q^n(0) \to q(0) \in \mathcal S$, and fix $T > 0$.
    Let
    \[
        R_T := q_F(0) + (1+\kfrags)\big(2q_I(0) + q_D(0) + 2\lambda T\big) + 1,
        \qquad
        C_T := T\big(\lambda + \mu_I + (c + \mu_D + \mu_F)R_T + cR_T^2\big).
    \]
    Then, for every $n$, the chain makes finitely many transitions in $[0,T]$, and almost surely, for all sufficiently large $n$,
    \[
        \sup_{t\leq T}\bar Q^n_D(t) \leq R_T, \qquad \sup_{t\leq T}\bar Q^n_F(t) \leq R_T,
    \]
    and each of $a(T), d_{IF}(T), d_{ID}(T), d_{IR}(T), d_{DF}(T), d_{DR}(T), d_F(T)$ is at most $C_T n$.
\end{lemma}
\begin{proof}
    Since $a(T) \leq n\lambda T$, at most $A(n\lambda T)$ launches are admitted in $[0,T]$, so at most $m^n_I := Q_I(0) + A(n\lambda T)$ distinct intacts are present during $[0,T]$.
    Every derelict was present at time $0$ or was an intact, so at most $m^n_D := Q_D(0) + m^n_I$ distinct derelicts are present during $[0,T]$, and $Q_D(t) \leq m^n_D$.
    Every collision removes an intact or a derelict, and each object is removed once, so $N_\mathrm{col}(T) \leq m^n := m^n_I + m^n_D$, and therefore
    \[
        Q_F(t) \leq Q_F(0) + \sum_{i=1}^{m^n} B_i, \qquad t \leq T.
    \]
    These quantities are finite, and each object decays at most once, which gives the first claim.
    By the strong law of large numbers, almost surely $A(n\lambda T)/n \to \lambda T$, hence $m^n_D/n \to q_I(0) + q_D(0) + \lambda T$ and $m^n/n \to 2q_I(0) + q_D(0) + 2\lambda T$, and $\frac1n\sum_{i\leq m^n}B_i \to \kfrags\,(2q_I(0) + q_D(0) + 2\lambda T)$.
    Since $q_I(0) + q_D(0) + \lambda T < R_T$ and $q_F(0) + \kfrags\,(2q_I(0) + q_D(0) + 2\lambda T) < R_T$, this gives the bounds on $\bar Q^n_D$ and $\bar Q^n_F$.
    Finally, on the event that these bounds hold, $Q_I \leq n$ and $\alpha^n \leq 1$ give
    \[
        a(T) \leq \lambda T n, \quad d_{IF}(T) \leq cR_T T n, \quad d_{ID}(T) + d_{IR}(T) \leq \mu_I T n,
    \]
    \[
        d_{DF}(T) \leq cR_T^2 T n, \quad d_{DR}(T) \leq \mu_D R_T T n, \quad d_F(T) \leq \mu_F R_T T n.
    \]
    \end{proof}

\subsection{Strong approximation}
\begin{theorem}
    \label{thm:strong-approx}
    Suppose $\bar Q^n(0) \to q(0) \in \mathcal S$.
    There exists a probability space which supports the processes $Q^n = (Q_I, Q_D, Q_F)$ for all $n$, seven independent standard Brownian motions \[W_A, W_{IF}, W_{ID}, W_{IR}, W_{DF}, W_{DR}, W_{F},\] and i.i.d.\ random variables $(B_i)_{i \in \bN}$ with distribution $\Poi(\kfrags)$, such that almost surely, for every $T > 0$, the following hold uniformly in $t \in [0, T]$ as $n\to\infty$:
    \begin{align}
        Q_I(t) & = Q_I(0) + a(t) - d_I(t) \nonumber
        \\ & \qquad + W_A(a(t)) - W_{IF}(d_{IF}(t)) - W_{ID}(d_{ID}(t)) - W_{IR}(d_{IR}(t)) + O(\log n),
        \label{eq:SA-I}
        \\
        Q_D(t) & = Q_D(0) + d_{ID}(t) - d_{D}(t) \nonumber
        \\ & \qquad + W_{ID}(d_{ID}(t)) - W_{DF}(d_{DF}(t)) - W_{DR}(d_{DR}(t)) + O(\log n),
        \label{eq:SA-D}
        \\
        Q_F(t) & = Q_F(0) + \kfrags\ d_{IF}(t) + \kfrags\ d_{DF}(t) - d_F(t) \nonumber
        \\ & \qquad + \kfrags\ W_{IF}(d_{IF}(t)) + \kfrags\ W_{DF}(d_{DF}(t)) - W_{F}(d_{F}(t)) + M(t) + O(\log n),
        \label{eq:SA-F}
    \end{align}
    where
    \[
        M(t) := \sum_{i=1}^{N_\mathrm{col}(t)} (B_i - \kfrags),
        \qquad
        N_{\mathrm{col}}(t) = D_{IF}(d_{IF}(t)) + D_{DF}(d_{DF}(t)).
    \]
\end{theorem}

\begin{proof}[Proof of Theorem~\ref{thm:strong-approx}]
    By the Koml\'os--Major--Tusn\'ady coupling~\cite{komlos1975,komlos1976}, a unit-rate Poisson process $N$ can be constructed on a probability space together with a standard Brownian motion $W$ such that, almost surely,
    \begin{equation}\label{eq:KMT}
        \sup_{s\leq u}\big|N(s) - s - W(s)\big| = O(\log u) \qquad \text{as } u\to\infty.
    \end{equation}
    Take seven independent copies of this pair, $(A, W_A), (D_{IF}, W_{IF}), \dots, (D_F, W_F)$, and an independent i.i.d.\ sequence $(B_i)$, and define $Q^n$ for every $n$ by~\eqref{eq:time-change}.
    By Lemma~\ref{lem:apriori}, almost surely, for all sufficiently large $n$ every cumulative intensity is at most $C_T n$ on $[0,T]$, so~\eqref{eq:KMT} gives, uniformly in $t \leq T$,
    \[
        A(a(t)) = a(t) + W_A(a(t)) + O(\log n),
    \]
    and likewise for the other six Poisson processes.
    Substituting in the first two lines of~\eqref{eq:time-change} gives~\eqref{eq:SA-I} and~\eqref{eq:SA-D}.
    For the fragments, $\sum_{i\leq N_\mathrm{col}(t)}B_i = \kfrags N_\mathrm{col}(t) + M(t)$ and
    \[
        N_\mathrm{col}(t) = d_{IF}(t) + d_{DF}(t) + W_{IF}(d_{IF}(t)) + W_{DF}(d_{DF}(t)) + O(\log n),
    \]
    which gives~\eqref{eq:SA-F}.
\end{proof}

\subsection{The fluid limit}
Let $b = (b_I, b_D, b_F)$ be the vector field of~\eqref{eq:FL-1}--\eqref{eq:FL-3},
\begin{align*}
    b_I(q) &:= \lambda - \mu_I q_I - \alpha(q_F)cq_Iq_F,
    \\
    b_D(q) &:= (1-P)\mu_I q_I - \mu_D q_D - cq_Dq_F,
    \\
    b_F(q) &:= \kfrags\,\alpha(q_F) c q_F q_I + \kfrags\, c q_F q_D - \mu_F q_F,
\end{align*}
which is locally Lipschitz on $\mathcal S$.
The capacity constraint is handled by the one-sided Skorokhod map at $q_I = 1$: for a right-continuous function $x$ with left limits and $x(0) \leq 1$, let
\[
    \Gamma(x)(t) := x(t) - \sup_{s\leq t}\big(x(s) - 1\big)^+.
\]
We use two standard facts.
First, $(z, \ell) = (\Gamma(x), x - \Gamma(x))$ is the unique pair such that $z = x - \ell \leq 1$, $\ell$ is non-decreasing with $\ell(0) = 0$, and $\ell$ increases only when $z = 1$, that is, $\int_0^\infty\ind_{z(s) < 1}\,d\ell(s) = 0$.
Second, $\|\Gamma(x) - \Gamma(y)\|_t \leq 2\|x - y\|_t$.

\begin{proposition}[Fluid limit with blocking]\label{prop:reflected-fluid}
    Let $\lambda > 0$, and suppose $\bar Q^n(0) \to q(0) \in \mathcal S$.
    \begin{enumerate}
        \item[\textup{(i)}] There is a unique pair $(q, \ell)$ such that $q : [0,\infty) \to \mathcal S$ is continuous, $\ell$ is continuous and non-decreasing with $\ell(0) = 0$, and
        \begin{equation}\label{eq:RFL}
        \begin{aligned}
            q_I(t) &= q_I(0) + \int_0^t b_I(q(s))\,ds - \ell(t),
            &
            q_D(t) &= q_D(0) + \int_0^t b_D(q(s))\,ds,
            \\
            q_F(t) &= q_F(0) + \int_0^t b_F(q(s))\,ds,
            &
            &\int_0^\infty \ind_{q_I(s) < 1}\,d\ell(s) = 0.
        \end{aligned}
        \end{equation}
        \item[\textup{(ii)}] Almost surely, $\|\bar Q^n - q\|_T \to 0$ for every $T > 0$.
        \item[\textup{(iii)}] If $\lambda \leq \mu_I$, then $\ell \equiv 0$ and $q$ is the solution of~\eqref{eq:FL-1}--\eqref{eq:FL-3}.
        In general,
        \[
            \ell(t) = \int_0^t \big(\lambda - \mu_I - \alpha(q_F(s))cq_F(s)\big)^+\,\ind_{q_I(s) = 1}\,ds,
        \]
        so launches are blocked, at the fluid scale, only while the shell is full and $\lambda > \mu_I + \alpha(q_F)cq_F$.
    \end{enumerate}
\end{proposition}
\begin{proof}
    \emph{Bounds for~\eqref{eq:RFL}.}
    Let $(q,\ell)$ satisfy~\eqref{eq:RFL} on $[0,T]$.
    Since $q_I \leq 1$ and $q_D, q_F \geq 0$, we have $b_D(q) \leq (1-P)\mu_I$, so $q_D(t) \leq q_D(0) + \mu_I T$.
    Moreover,
    \[
        \kfrags\, b_I(q) + \kfrags\, b_D(q) + b_F(q) = \kfrags\ \lambda - \kfrags\  P\mu_I q_I - \kfrags\ \mu_D q_D - \mu_F q_F \leq \kfrags\ \lambda,
    \]
    and $\ell$ is non-decreasing, so $q_F(t) \leq \kfrags\ q_I(t) + \kfrags\ q_D(t) + q_F(t) \leq q_F(0) + \kfrags(q_I(0) + q_D(0) + \lambda T)$.
    Let
    \[
        R := \max\big\{R_T,\ q_D(0) + \mu_I T,\ q_F(0) + \kfrags(q_I(0) + q_D(0) + \lambda T)\big\},
        \qquad
        K_T := [0,1]\times[0,R]^2,
    \]
    and let $L$ be a Lipschitz constant of $b$ on $K_T$ with respect to $|\cdot|$.
    Then $q(t) \in K_T$ for $t \leq T$, and by Lemma~\ref{lem:apriori}, almost surely $\bar Q^n(t) \in K_T$ for all $t \leq T$ and all sufficiently large $n$.

    \emph{Proof of (i).}
    By the first fact about $\Gamma$,~\eqref{eq:RFL} is equivalent to
    \[
        q_I = \Gamma\Big(q_I(0) + \int_0^\cdot b_I(q(s))\,ds\Big), \qquad q_D = q_D(0) + \int_0^\cdot b_D(q(s))\,ds, \qquad q_F = q_F(0) + \int_0^\cdot b_F(q(s))\,ds,
    \]
    with $\ell$ the corresponding regulator.
    Since $\Gamma$ is Lipschitz and $b$ is locally Lipschitz, Picard iteration gives a unique solution on a short time interval, which can be continued for as long as it remains bounded.
    The solution remains in $\mathcal S$: $b_D \geq 0$ when $q_D = 0$, $b_F = 0$ when $q_F = 0$, and $b_I = \lambda > 0$ when $q_I = 0$, while $\ell$ increases only when $q_I = 1$.
    The bounds above then show that the solution exists, and is unique, on $[0,\infty)$.

    We divide~\eqref{eq:SA-I}--\eqref{eq:SA-F} by $n$.
    Since $\alpha^n(Q_F)\tfrac{c}{n}Q_FQ_I = n\,\alpha(\bar Q^n_F)c\,\bar Q^n_F\bar Q^n_I$, and similarly for the other intensities, the intensities other than $a$ reproduce the vector field $b$ exactly, and
    \[
        \frac{a(t)}{n} = \lambda t - \ell^n(t), \qquad \ell^n(t) := \lambda\int_0^t \ind_{\bar Q^n_I(s) = 1}\,ds.
    \]
    Therefore
    \begin{equation}\label{eq:prelimit}
    \begin{aligned}
        \bar Q^n_I(t) &= \bar Q^n_I(0) + \int_0^t b_I(\bar Q^n(s))\,ds - \ell^n(t) + \varepsilon^n_I(t),
        \\
        \bar Q^n_D(t) &= \bar Q^n_D(0) + \int_0^t b_D(\bar Q^n(s))\,ds + \varepsilon^n_D(t),
        \\
        \bar Q^n_F(t) &= \bar Q^n_F(0) + \int_0^t b_F(\bar Q^n(s))\,ds + \varepsilon^n_F(t),
    \end{aligned}
    \end{equation}
    where $\varepsilon^n_I, \varepsilon^n_D, \varepsilon^n_F$ are $1/n$ times the second lines of~\eqref{eq:SA-I},~\eqref{eq:SA-D},~\eqref{eq:SA-F}, respectively.
    The process $\ell^n$ is continuous and non-decreasing, $\ell^n(0) = 0$, it increases only when $\bar Q^n_I = 1$, and $\bar Q^n_I \leq 1$.
    By the first fact about $\Gamma$,
    \[
        \bar Q^n_I = \Gamma(X^n), \qquad X^n(t) := \bar Q^n_I(0) + \int_0^t b_I(\bar Q^n(s))\,ds + \varepsilon^n_I(t).
    \]

    We now show that the fluctuations vanish on the fluid scale.
    Fix $T$.
    Almost surely, $\sup_{s\leq u}|W(s)|/u \to 0$ as $u\to\infty$ for each of the seven Brownian motions, and by Lemma~\ref{lem:apriori} their arguments are at most $C_T n$ on $[0,T]$ for large $n$; hence each term $\tfrac1n W_\bullet(d_\bullet(t))$ tends to $0$ uniformly on $[0,T]$, and so does $O(\log n)/n$.
    For the batch term, $N_\mathrm{col}(T) \leq m^n$ with $m^n/n$ convergent (Proof of Lemma~\ref{lem:apriori}), and $\max_{m\leq m^n}\big|\sum_{i\leq m}(B_i - \kfrags)\big|/n \to 0$ almost surely by the strong law of large numbers, so $\|M\|_T/n \to 0$.
    Thus $\|\varepsilon^n\|_T \to 0$ almost surely, where $\varepsilon^n := (\varepsilon^n_I,\varepsilon^n_D,\varepsilon^n_F)$.

    \emph{Proof of (ii).}
    Let $x(t) := q_I(0) + \int_0^t b_I(q(s))\,ds$, so that $q_I = \Gamma(x)$, and let $\Delta^n(t) := \|\bar Q^n - q\|_t$.
    For $n$ so large that $\bar Q^n$ stays in $K_T$ on $[0,T]$, the Lipschitz property of $\Gamma$ gives, for $t \leq T$,
    \[
        \|\bar Q^n_I - q_I\|_t \leq 2\|X^n - x\|_t \leq 2\Big(|\bar Q^n(0) - q(0)| + \|\varepsilon^n\|_T + L\int_0^t\Delta^n(s)\,ds\Big),
    \]
    and~\eqref{eq:prelimit} and~\eqref{eq:RFL} give the same bound, without the factor $2$, for the other two coordinates.
    Hence $\Delta^n(t) \leq 2\big(|\bar Q^n(0) - q(0)| + \|\varepsilon^n\|_T\big) + 2L\int_0^t\Delta^n(s)\,ds$, and Gr\"onwall's inequality gives
    \[
        \Delta^n(T) \leq 2\big(|\bar Q^n(0) - q(0)| + \|\varepsilon^n\|_T\big)e^{2LT} \to 0.
    \]

    \emph{Proof of (iii).}
    Let $\lambda \leq \mu_I$, and let $\widetilde q$ be the solution of~\eqref{eq:FL-1}--\eqref{eq:FL-3} started at $q(0)$, which remains non-negative.
    On its interval of existence, by~\eqref{eq:FL-1},
    \[
        \frac{d}{dt}(1 - \widetilde q_I) = \mu_I\widetilde q_I - \lambda + \alpha(\widetilde q_F)c\,\widetilde q_I\widetilde q_F \geq -\mu_I(1 - \widetilde q_I),
    \]
    so $1 - \widetilde q_I(t) \geq (1 - q_I(0))e^{-\mu_I t} \geq 0$.
    Thus $\widetilde q$ takes values in $\mathcal S$, the bounds above apply to it with $\ell \equiv 0$, and it exists for all time.
    The pair $(\widetilde q, 0)$ satisfies~\eqref{eq:RFL}; by uniqueness, $q = \widetilde q$ and $\ell \equiv 0$.

    For the formula, note that $x$ is Lipschitz, hence so is $\ell(t) = \sup_{s\leq t}(x(s) - 1)^+$, and $\ell' = 0$ almost everywhere on the open set $\{q_I < 1\}$.
    On the level set $\{q_I = 1\}$ the derivative of the Lipschitz function $q_I$ vanishes almost everywhere, so $\ell' = b_I(q) = \lambda - \mu_I - \alpha(q_F)cq_F$ almost everywhere there; since $\ell' \geq 0$, this equals its positive part.
\end{proof}

In particular, Theorem~\ref{thm:fluid-limit} is the case $\lambda \leq \mu_I$ of Proposition~\ref{prop:reflected-fluid}(ii) and~(iii).

\begin{remark}[Time-varying launch rate]\label{rem:time-varying}
    Lemma~\ref{lem:apriori}, Theorem~\ref{thm:strong-approx}, and Proposition~\ref{prop:reflected-fluid} and their proofs remain valid, with the same arguments, when the launch rate is a continuous function $\lambda(t) > 0$: the $n$-th system has arrival rate $n\lambda(t)$, $a(t) = \int_0^t n\lambda(s)\ind_{Q_I(s) < n}\,ds$, $\lambda T$ is replaced by $\int_0^T\lambda(s)\,ds$ in the bounds, and $b_I$ depends on time through $\lambda(t)$.
    In particular, if $\lambda(t) \leq \mu_I$ for all $t$, then the fluid limit is~\eqref{eq:FL-1}--\eqref{eq:FL-3} with $\lambda$ replaced by $\lambda(t)$.
    This is the fluid limit used in Theorem~\ref{thm:slow-passage}, where the limit $n\to\infty$ is taken first, at a fixed ramp rate $\varepsilon$.
\end{remark}

\section{Proofs: Bifurcation Structure and Slow Passage}

\subsection{Proof of the bifurcation structure}\label{sec:proof-bifurcation}

We prove Theorem~\ref{thm:bifurcation-structure}. Recall $H$ from~\eqref{eq:H} and the fixed-point relation $\lambda = \cM(q_F) = \mu_F/(c\,\kfrags\,H(q_F))$ from~\eqref{eq:cM}.
We isolate the stability analysis in Lemma~\ref{lem:linearization}.

Throughout, for a non-trivial fixed point $q^* = (q^*_I, q^*_D, q^*_F)$ we write $z := q^*_F > 0$ and
\begin{gather*}
    a := \alpha(z) c z, \qquad b := c z, \qquad \sigma_I := \mu_I + a, \qquad \sigma_D := \mu_D + b,
    \\
    e := \frac{z\alpha'(z)}{\alpha(z)}, \qquad \kappa := 1 - \alpha(z) - z\alpha'(z), \qquad \theta := \frac{\alpha(z) q^*_I}{\alpha(z) q^*_I + q^*_D}, \qquad w := e\,\theta\,\mu_F.
\end{gather*}
Here, $\sigma_I$ and $\sigma_D$ are the total rates at which an intact and a derelict leave their stations, and $\theta$ is the fraction of fragment production which is due to intacts.
From $\alpha(q) = (\alpha_0\Qsat + q)/(\Qsat + q)$ we compute
\begin{equation}\label{eq:e-kappa}
    e = \frac{(1-\alpha_0)\Qsat\, z}{(z + \Qsat)(z + \alpha_0\Qsat)},
    \qquad
    \kappa = \frac{(1-\alpha_0)\Qsat^2}{(z + \Qsat)^2},
\end{equation}
so that $e, \kappa \in (0,1)$ and $\theta \in (0,1)$.
Moreover, $z/((z+\Qsat)(z+\alpha_0\Qsat)) = \big(z + (1+\alpha_0)\Qsat + \alpha_0\Qsat^2/z\big)^{-1}$ is maximized at $z = \sqrt{\alpha_0}\,\Qsat$, which gives
\begin{equation}\label{eq:e-bound}
    e \leq \frac{1-\alpha_0}{(1+\sqrt{\alpha_0})^2} = \frac{1-\sqrt{\alpha_0}}{1+\sqrt{\alpha_0}} = \varepsilon_0.
\end{equation}

\begin{lemma}[Linearization at a non-trivial fixed point]\label{lem:linearization}
    Let $q^*$ be a non-trivial fixed point of~\eqref{eq:FL-1}--\eqref{eq:FL-3}, and let $\cL$ be the Jacobian at $q^*$.
    The characteristic polynomial of $\cL$ is $\chi(s) = s^3 + c_2s^2 + c_1s + c_0$, with
    \begin{align*}
        c_2 &= \sigma_I + \sigma_D - w,
        \\
        c_1 &= \sigma_I\sigma_D - w(\sigma_I + \sigma_D) + \mu_F\, b\,(1 - \theta\kappa),
        \\
        c_0 &= -\sigma_I\sigma_D\, w + \mu_F\, b\,\big((1-\theta)\sigma_I + (1-\kappa)\sigma_D\big) = -\sigma_I\sigma_D\,\mu_F\, z\,\frac{H'(z)}{H(z)},
    \end{align*}
    and
    \begin{equation}\label{eq:RH-gap}
        c_1c_2 - c_0 = (\sigma_I + \sigma_D)(\sigma_I - w)(\sigma_D - w) + \mu_F\, b\,\Big(\theta(1-\kappa)(\sigma_I - w) + (1-\theta)(\kappa\sigma_D - w)\Big).
    \end{equation}
    Consequently:
    \begin{enumerate}
        \item[\textup{(i)}] If $H'(z) > 0$, then $\cL$ has a positive real eigenvalue.
        \item[\textup{(ii)}] If $H'(z) < 0$, then all eigenvalues of $\cL$ have negative real part if and only if $c_2 > 0$ and $c_1c_2 - c_0 > 0$.
        \item[\textup{(iii)}] If $e\,\mu_F \leq \min\{\mu_I,\mu_D\}$, then $c_2 > 0$ and $c_1c_2 - c_0 > 0$. In particular, this holds at every non-trivial fixed point under~\eqref{eq:rate-condition}.
    \end{enumerate}
\end{lemma}
\begin{proof}
    \emph{The Jacobian.}
    Differentiating~\eqref{eq:FL-1}--\eqref{eq:FL-3} gives
    \[
        \cL = \begin{pmatrix}
            -\sigma_I & 0 & -c\,(\alpha(z) + z\alpha'(z))\,q^*_I
            \\
            (1-P)\mu_I & -\sigma_D & -c\,q^*_D
            \\
            \kfrags\, a & \kfrags\, b & \kfrags\, c\, q^*_I(\alpha(z) + z\alpha'(z)) + \kfrags\, c\, q^*_D - \mu_F
        \end{pmatrix}.
    \]
    At a non-trivial fixed point,~\eqref{eq:FL-FP-3} gives $\kfrags\, c\,(\alpha(z)q^*_I + q^*_D) = \mu_F$, so that, by the definition of $\theta$,
    \[
        \kfrags\, c\,\alpha(z)\, q^*_I = \theta\mu_F, \qquad \kfrags\, c\, q^*_D = (1-\theta)\mu_F.
    \]
    Hence the matrix entry $\cL_{33} = \kfrags\, c\, q^*_I z\alpha'(z) = e\,\theta\,\mu_F = w$.
    Notice that $\cL_{33} > 0$.
    Similarly, using $\alpha(z) + z\alpha'(z) = 1-\kappa$ and $a = \alpha(z) b$,
    \[
        -\cL_{13}\cL_{31} = \kfrags\, c\, q^*_I(1-\kappa)\,a = \theta\mu_F(1-\kappa)\, b,
        \qquad
        -\cL_{23}\cL_{32} = \kfrags\, c\, q^*_D\, b = (1-\theta)\mu_F\, b.
    \]
    We have $c_2 = -\mathrm{tr}(\cL) = \sigma_I + \sigma_D - w$.
    The coefficient $c_1$ is the sum of the principal $2\times 2$ minors of $\cL$:
    \[
        c_1 = \sigma_I\sigma_D + \big(-\sigma_I w - \cL_{13}\cL_{31}\big) + \big(-\sigma_D w - \cL_{23}\cL_{32}\big)
        = \sigma_I\sigma_D - w(\sigma_I + \sigma_D) + \mu_F\, b\,\big(\theta(1-\kappa) + (1-\theta)\big).
    \]
    Expanding the determinant along the first row and using $\cL_{21}\cL_{32} = (1-P)\mu_I\kfrags\, b$ gives
    \[
        \det \cL = \sigma_I\sigma_D\, w + \sigma_I \cL_{23}\cL_{32} + \cL_{13}\big(\cL_{21}\cL_{32} + \sigma_D \cL_{31}\big).
    \]
    By~\eqref{eq:fp-intact-derelict}, $(1-P)\mu_I q^*_I = \sigma_D q^*_D$, so $\cL_{13}\cL_{21}\cL_{32} = -(1-\kappa)\sigma_D\,\kfrags\, c\, q^*_D\, b = -(1-\kappa)(1-\theta)\sigma_D\mu_F b$, and therefore
    \begin{align*}
        c_0 = -\det \cL &= -\sigma_I\sigma_D\, w + \mu_F\, b\,\big((1-\theta)\sigma_I + (1-\kappa)(1-\theta)\sigma_D + \theta(1-\kappa)\sigma_D\big)
        \\
        &= -\sigma_I\sigma_D\, w + \mu_F\, b\,\big((1-\theta)\sigma_I + (1-\kappa)\sigma_D\big).
    \end{align*}

    \emph{The relation between $c_0$ and $H'$.}
    For $\zeta > 0$ let $(x(\zeta), y(\zeta))$ be the unique solution of~\eqref{eq:FL-FP-1},~\eqref{eq:FL-FP-2} with $q^*_F$ replaced by $\zeta$, which is given by~\eqref{eq:fp-intact-derelict}, and let $g(\zeta)$ be the right-hand side of~\eqref{eq:FL-3} evaluated at $(x(\zeta), y(\zeta), \zeta)$.
    By the definition of $H$,
    \[
        g(\zeta) = \zeta\big(\lambda c\,\kfrags H(\zeta) - \mu_F\big).
    \]
    Write $\cL_{[12]}$ for the upper-left $2\times 2$ block of $\cL$, which is lower triangular with $\det \cL_{[12]} = \sigma_I\sigma_D$.
    Differentiating the identities~\eqref{eq:FL-FP-1},~\eqref{eq:FL-FP-2} in $\zeta$ at $\zeta = z$ gives $(x'(z), y'(z))^\top = -\cL_{[12]}^{-1}(\cL_{13}, \cL_{23})^\top$, so by the chain rule
    \[
        g'(z) = \cL_{33} - (\cL_{31}, \cL_{32})\,\cL_{[12]}^{-1}\,(\cL_{13}, \cL_{23})^\top,
    \]
    which is the Schur complement of $\cL_{[12]}$ in $\cL$. Thus $\det \cL = \sigma_I\sigma_D\, g'(z)$.
    Since $\lambda c\,\kfrags H(z) = \mu_F$ at the fixed point, $g'(z) = z\lambda c\,\kfrags H'(z) = \mu_F z H'(z)/H(z)$, which establishes the second expression for $c_0$.

    \emph{The identity~\eqref{eq:RH-gap}.}
    Write $R := \mu_F b(1-\theta\kappa)$ and $S := \mu_F b((1-\theta)\sigma_I + (1-\kappa)\sigma_D)$, so that $c_1 = \sigma_I\sigma_D - w(\sigma_I+\sigma_D) + R$ and $c_0 = -\sigma_I\sigma_D w + S$.
    Then
    \begin{align*}
        c_1c_2 - c_0 &= (\sigma_I + \sigma_D)\big(\sigma_I\sigma_D - w(\sigma_I + \sigma_D) + w^2\big) + (\sigma_I + \sigma_D - w)R - S
        \\
        &= (\sigma_I + \sigma_D)(\sigma_I - w)(\sigma_D - w) + (\sigma_I + \sigma_D - w)R - S,
    \end{align*}
    and
    \begin{align*}
        \frac{(\sigma_I + \sigma_D - w)R - S}{\mu_F\, b}
        &= (1-\theta\kappa)(\sigma_I + \sigma_D - w) - (1-\theta)\sigma_I - (1-\kappa)\sigma_D
        \\
        &= \theta(1-\kappa)\sigma_I + \kappa(1-\theta)\sigma_D - (1-\theta\kappa)w
        \\
        &= \theta(1-\kappa)(\sigma_I - w) + (1-\theta)(\kappa\sigma_D - w),
    \end{align*}
    where the last step uses $1 - \theta\kappa = \theta(1-\kappa) + (1-\theta)$.

    \emph{Proof of (i).}
    If $H'(z) > 0$ then $\chi(0) = c_0 < 0$, and $\chi(s)\to\infty$ as $s\to\infty$, so $\chi$ has a positive real root.

    \emph{Proof of (ii).}
    If $H'(z) < 0$ then $c_0 > 0$, and the claim is the Routh--Hurwitz criterion for a cubic: all roots of $\chi$ have negative real part if and only if $c_2 > 0$, $c_0 > 0$, and $c_1c_2 > c_0$.

    \emph{Proof of (iii).}
    Let $\mu_{\min} := \min\{\mu_I,\mu_D\}$ and suppose $e\mu_F \leq \mu_{\min}$, so that $w = e\theta\mu_F \leq \theta\mu_{\min}$.
    Then
    \[
        \sigma_I - w \geq a + (1-\theta)\mu_I > 0, \qquad \sigma_D - w \geq b + (1-\theta)\mu_D > 0,
    \]
    so $c_2 > 0$ and the first term of~\eqref{eq:RH-gap} is strictly positive.
    For the second term of~\eqref{eq:RH-gap}, we use $\sigma_I - w \geq (1-\theta)\mu_I$, $\sigma_D \geq \mu_D$, and $w \leq \theta\mu_{\min}$ to obtain
    \begin{align*}
        \theta(1-\kappa)(\sigma_I - w) + (1-\theta)(\kappa\sigma_D - w)
        &\geq (1-\theta)\big(\theta(1-\kappa)\mu_I + \kappa\mu_D - \theta\mu_{\min}\big)
        \\
        &\geq (1-\theta)\,\mu_{\min}\big(\theta(1-\kappa) + \kappa - \theta\big)
        = (1-\theta)^2\kappa\,\mu_{\min} \geq 0.
    \end{align*}
    Hence $c_1c_2 - c_0 > 0$.
    The final assertion follows from~\eqref{eq:e-bound}.
\end{proof}

\begin{proof}[Proof of Theorem~\ref{thm:bifurcation-structure}]
    The number of positive roots of $H'$ determines the fixed-point structure, and Lemma~\ref{lem:linearization} determines the stability.
    By the standing assumption $\lambda \leq \mu_I$ and Proposition~\ref{prop:reflected-fluid}(iii), \eqref{eq:FL-1}--\eqref{eq:FL-3} is the fluid limit and the state space $\mathcal{S}$ is invariant, and by~\eqref{eq:fp-intact-derelict} every non-trivial fixed point has $q^*_I = \lambda/(\mu_I + \alpha(q^*_F)cq^*_F) < 1$, so the local stability of the fixed points is therefore determined by the Jacobian of~\eqref{eq:FL-1}--\eqref{eq:FL-3}.
    Notice that $H > 0$ on $[0,\infty)$ and $H(q)\to 0$ as $q\to\infty$, so $\cM$ is smooth and positive on $[0,\infty)$ with $\cM(q)\to\infty$ as $q\to\infty$, and
    \[
        \cM(0) = \frac{\mu_F}{c\,\kfrags\,H(0)} = \frac{\mu_I\mu_D\mu_F}{\kfrags\, c\,(\alpha_0\mu_D + (1-P)\mu_I)} = \lambda_\crit.
    \]

    \textbf{The trivial fixed point. }
    At the trivial fixed point $(\lambda/\mu_I,\ (1-P)\lambda/\mu_D,\ 0)$ the third row of the Jacobian of~\eqref{eq:FL-1}--\eqref{eq:FL-3} is $(0,\ 0,\ \kfrags\, c(\alpha_0 q^*_I + q^*_D) - \mu_F)$, and the upper-left $2\times 2$ block is lower triangular with diagonal entries $-\mu_I$ and $-\mu_D$.
    The eigenvalues are therefore $-\mu_I$, $-\mu_D$, and
    \[
        \kfrags\, c\,\lambda H(0) - \mu_F = \mu_F\left(\frac{\lambda}{\lambda_\crit} - 1\right),
    \]
    which is negative if $\lambda < \lambda_\crit$ and positive if $\lambda > \lambda_\crit$.

    \textbf{Proof of Statement 1. }
    Notice that the second quantity in the expression~\eqref{eq:H} for $H(q)$,
    \[\frac{(1-P)\mu_I}{(\mu_I + \alpha(q) c q)(\mu_D + cq)},\]
    is strictly decreasing in $q$, since $q\mapsto\alpha(q)q$ is increasing.
    By the quotient rule, the derivative of the first quantity $\frac{\alpha(q)}{\mu_I + \alpha(q) c q}$ has the sign of
    \[\alpha'(q)\mu_I - c\,\alpha(q)^2.\]
    Since $\alpha'(q) = (1-\alpha_0)\Qsat/(q+\Qsat)^2$ is decreasing in $q$ and $\alpha$ is increasing in $q$, this quantity is negative for all $q \geq 0$ if and only if it is negative at $q = 0$, that is, if and only if $(1-\alpha_0)\mu_I/\Qsat < c\,\alpha_0^2$, which is
    \[\Qsat > \overline{Q} = \frac{(1 - \alpha_0)\mu_I}{c\, \alpha_0^2}.\]
    In this case $H' < 0$ on $[0,\infty)$, so $\cM$ is strictly increasing, and it is a bijection from $(0,\infty)$ onto $(\lambda_\crit,\infty)$.
    This establishes the first statement.

    This is only a sufficient condition for the absence of folds.
    Notice that $\overline{Q}$ is only a function of $\mu_I$ and not of $\mu_D$ or $P$, and does not depend on $\kfrags$.
    This is due to the fact that only the first term of $H$ is a non-monotone function of $q$, and the non-monotonicity is the cause of the saddle-node bifurcation.

    \textbf{Proof of Statement 2. }
    For the rest of the proof of the fixed-point structure, we argue via Descartes' rule of signs; the third and fourth statements need further arguments beyond a simple application of the rule.
    Denote \[L(q) := (\mu_I + \alpha(q)cq)(q + \Qsat) = cq^2 + (\alpha_0\Qsat c + \mu_I)q + \mu_I\Qsat.\]
    Then,
    \begin{align*}
        H(q) &= \frac{(\alpha_0\Qsat + q)(cq + \mu_D)}{L(q)(cq + \mu_D)} + \frac{(1-P)\mu_I(q + \Qsat)}{L(q)(cq + \mu_D)}
        =: \frac{N(q)}{D(q)},
    \end{align*}
    and since $D > 0$ on $[0,\infty)$, the roots of $H'$ are the roots of
    \[\pi(q) := N'(q)D(q) - N(q)D'(q) = a_4q^4 + a_3q^3 + a_2q^2 + a_1q + a_0,\]
    with $H'(q) = \pi(q)/D(q)^2$.
    Let
    \[\gamma_N := \mu_D + c\alpha_0\Qsat + (1-P)\mu_I > 0,\]
    \[\gamma_1 := c\left( \alpha_0^2\mu_D^2 + (1-P)\mu_I^2 + (1-P)\alpha_0\mu_I\mu_D \right), \qquad \gamma_2 := \mu_I\mu_D^2(1 - \alpha_0).\]
    Expanding $\pi$ gives
    \[a_4 = -c^3 < 0, \qquad a_3 = -2c^2\gamma_N < 0, \qquad a_0 = -\gamma_1\Qsat^2 + \gamma_2\Qsat,\]
    and since $\Qsat > 0$, the sign of $a_0$ is the sign of $-\gamma_1\Qsat + \gamma_2$.
    Thus $a_0 < 0$ when $\Qsat > \underline{Q} = \gamma_2/\gamma_1$ and $a_0 > 0$ when $\Qsat < \underline{Q}$.
    Also, $\gamma_1 > c\,\alpha_0^2\mu_D^2$ gives $\underline{Q} < \overline{Q}$.

    If $\Qsat > \underline{Q}$, the sign sequence of $(a_4, a_3, a_2, a_1, a_0)$ is $(-,-,\ast,\ast,-)$, which has either zero or two sign changes.
    By Descartes' rule of signs, $\pi$ has either zero or two positive roots, counted with multiplicity; since $\pi(0) < 0$ and $\pi(q)\to-\infty$, a double root is not a sign change of $H'$, so there are either no folds or exactly two folds.

    It remains to show that the transcritical bifurcation is supercritical.
    Since $\operatorname{sign} H'(0) = \operatorname{sign} a_0$, for $\Qsat > \underline{Q}$ we have $H'(0) < 0$ and
    \[
        \cM'(0) = -\frac{\lambda_\crit H'(0)}{H(0)} > 0.
    \]
    The trivial branch $\{q_F = 0\}$ and the non-trivial branch $\{\lambda = \cM(q_F)\}$ therefore cross transversally at $(\lambda, q_F) = (\lambda_\crit, 0)$, and the non-trivial branch emerges into $\lambda > \lambda_\crit$, with
    \[
        q_F^*(\lambda) = \frac{H(0)}{\lambda_\crit\,\lvert H'(0)\rvert}\,(\lambda - \lambda_\crit) + o(\lambda - \lambda_\crit) > 0.
    \]
    The trivial fixed point loses stability at $\lambda_\crit$, as shown above.
    On the emerging branch, as $z = q^*_F \downarrow 0$ we have $w \to 0$ and $b \to 0$, so in the notation of Lemma~\ref{lem:linearization}, $c_2 \to \mu_I + \mu_D$, $c_1 \to \mu_I\mu_D$, and $c_0 \to 0$ with $c_0 > 0$.
    Thus, for $\lambda - \lambda_\crit > 0$ sufficiently small, two eigenvalues are close to $-\mu_I$ and $-\mu_D$ and the third is $-c_0/c_1 + o(c_0) < 0$: the emerging branch is locally asymptotically stable, and the trivial fixed point cedes its stability to it.
    This establishes the second statement.

    \textbf{Proof of Statement 3. }
    A direct computation gives
    \[ a_1 = 2c\,\Qsat\,(\beta_0 - \gamma_0 c\alpha_0\Qsat), \qquad \gamma_0 = \alpha_0 \mu_D + (1-P)\mu_I, \qquad \beta_0 = (\mu_I+\mu_D)(P\mu_I-\alpha_0 \mu_D) - \mu_I^2. \]
    Assume that $a_1 > 0$ at $\Qsat = \underline{Q}$.
    At $\Qsat = \underline{Q}$ we have $a_0 = 0$, so $\pi(q) = a_1q + O(q^2) > 0$ for all small $q>0$, and we fix $q_0>0$ with $\pi(q_0)\big|_{\Qsat=\underline{Q}} > 0$.
    As the coefficients of $\pi$ are continuous in $\Qsat$, there is $\delta>0$ with $\pi(q_0) > 0$ for $\Qsat \in [\underline{Q}, \underline{Q}+\delta)$, and we shrink $\delta$ so that $\underline{Q}+\delta \le \overline{Q}$.
    For any fixed $\Qsat \in (\underline{Q}, \underline{Q}+\delta)$ we have $\pi(0) = a_0 < 0$, $\pi(q_0) > 0$, and $\pi(q)\to-\infty$, hence, by the Proof of Statement 2, $\pi$ has exactly two positive roots $q_{\min} < q_0 < q_{\max}$, both simple.
    Then $H' < 0$ on $(0, q_{\min})$, $H' > 0$ on $(q_{\min}, q_{\max})$, and $H' < 0$ on $(q_{\max},\infty)$, so $\cM$ increases from $\lambda_\crit$ to $\lambda_+ = \cM(q_{\min})$, decreases to $\lambda_- = \cM(q_{\max})$, and increases to $\infty$.
    For $\lambda\in(\max\{\lambda_-,\lambda_\crit\},\lambda_+)$ the equation $\cM(q_F) = \lambda$ thus has exactly three solutions.

    We note that this is only a sufficient condition and that it is possible that other conditions can create two folds.

    \textbf{Proof of Statement 4. }
    Let $0 < \Qsat < \underline{Q}$, so that $a_0 > 0$ and $H'(0) > 0$.
    Then $\cM'(0) < 0$, and the non-trivial branch
    \[ q_F^*(\lambda) = \frac{H(0)}{\lambda_\crit\,H'(0)}\,(\lambda_\crit - \lambda) + o(\lambda_\crit - \lambda) \]
    is positive precisely for $\lambda < \lambda_\crit$: the transcritical bifurcation is subcritical.
    Since $\pi(0) = a_0 > 0$ and $\pi(q)\to-\infty$, $\pi$ has an odd number of positive roots, hence one or three.
    We show that there is exactly one, $q_{\max}$.
    Then $H' > 0$ on $(0, q_{\max})$ and $H' < 0$ on $(q_{\max},\infty)$, so $\cM$ decreases from $\lambda_\crit$ to $\lambda_- = \cM(q_{\max})$ and then increases to $\infty$, which gives the claimed number of fixed points.

    The remainder of the proof is computation-heavy, so we give a brief overview.
    If $a_2 \leq 0$, or if $a_2 > 0$ and $a_1 \geq 0$, then the sign sequence of $(a_4,a_3,a_2,a_1,a_0)$ has exactly one sign change and, by Descartes' rule of signs, $\pi$ has exactly one positive root.
    Thus, to finish the proof, we rule out the situation
    \[a_2>0 \qquad\text{and}\qquad a_1<0.\]
    We pose this as a feasibility problem for a certain linear program and show that it is empty.
    This problem reduces to a finite computation over the rationals which we do via a computer algorithm.

    With $a_1=2c\Qsat(\beta_0-\gamma_0 c\alpha_0\Qsat)$ the set $\{a_1<0\}$ is the half-line $\Qsat>\frac{\beta_0}{\gamma_0 c\alpha_0}$ (all of $\Qsat>0$ when $\beta_0\le0$), and with \[a_2=-c\,\Theta(\Qsat),\]
    \[\Theta(\Qsat)=\alpha_0^2c^2\,\Qsat^2 + c\big(4\alpha_0 \mu_D+2(1+\alpha_0)\mu_I-(3+\alpha_0)P\mu_I\big)\Qsat + (1-P)\mu_I(\mu_I+\mu_D)+\mu_D^2, \]
    the set $\{a_2>0\}=\{\Theta<0\}$. The pattern $a_2 > 0, a_1 < 0$ occurs in $(0,\underline{Q})$ iff $\Theta(\Qsat)<0$ for some $\Qsat\in[L,\underline{Q}]$ with $L=\max\{0,\beta_0/(\gamma_0 c\alpha_0)\}$.
    Obviously, if $L > \underline{Q}$, then we have $a_1 \geq 0$ and $\pi(q)$ only has one positive root.

    Since $\Theta$ is convex ($\partial_{\Qsat}^2\Theta=2\alpha_0^2c^2>0$) and $\Theta(0)=(1-P)\mu_I(\mu_I+\mu_D)+\mu_D^2>0$, its minimum on $[L,\underline{Q}]$ is attained at the vertex clamped to the interval; the pattern is empty iff
    \[ M:=\min_{\Qsat\in[L,\underline{Q}]}\Theta \ge\ 0,\qquad Q_v=\frac{-\big(4\alpha_0 \mu_D+2(1+\alpha_0)\mu_I-(3+\alpha_0)P\mu_I\big)}{2\alpha_0^2 c}, \]
    where $Q_v$ is the vertex of $\Theta$.

    To reduce the polynomials below to the three variables $(\alpha_0, P, r)$, we show that without loss of generality, we can take $c = \mu_I = 1$.

    Put $\widetilde q = cq$, $\widetilde{Q}_\mathrm{sat}=c\Qsat$.
    Then $\alpha(q) =\frac{\widetilde q+\alpha_0\widetilde{Q}_\mathrm{sat}}{\widetilde q+\widetilde{Q}_\mathrm{sat}} =: \widetilde\alpha(\widetilde q)$, $\alpha(q) cq+\mu_I=\widetilde\alpha(\widetilde{q})\widetilde q+\mu_I$, $cq+\mu_D=\widetilde q+\mu_D$, so
    \[H=\frac{\widetilde\alpha(\widetilde q)}{\widetilde\alpha(\widetilde q)\widetilde q+\mu_I}+(1-P)\frac{\mu_I}{(\widetilde q+\mu_D)(\widetilde\alpha(\widetilde q)\widetilde q+\mu_I)},\]
    thus we can take $c = 1$.
    Next, put $q=\mu_I\widehat q$, $\Qsat=\mu_I\widehat Q_\mathrm{sat}$, and let \[r:=\mu_D/\mu_I.\]
    Again, $\alpha(q) =\frac{\widehat q+\alpha_0\widehat Q_\mathrm{sat}}{\widehat q+\widehat Q_\mathrm{sat}} =: \widehat\alpha(\widehat q)$, while $\alpha(q) q+\mu_I=\mu_I(\widehat\alpha(\widehat q)\widehat q+1)$ and $q+\mu_D=\mu_I(\widehat q+r)$.
    Therefore
    \[H=\frac{1}{\mu_I}\Big[\frac{\widehat\alpha(\widehat q)}{\widehat\alpha(\widehat q)\widehat q+1}+(1-P)\frac{1}{(\widehat q+r)(\widehat\alpha(\widehat q)\widehat q+1)}\Big],\]
    and thus we can also take $\mu_I = 1$.

    We use the normalizations $c=\mu_I=1$, $\mu_D=r$, and write $\bar P=1-P$, $\bar\alpha_0=1-\alpha_0$, so that
    \[ \Theta(\Qsat)=\alpha_0^2\Qsat^2+b_1\Qsat+b_0,\quad b_1=4\alpha_0 r+(3+\alpha_0)\bar P-\bar\alpha_0,\quad b_0=\bar P(1+r)+r^2. \]
    Put $g_1=\alpha_0^2 r^2+\bar P(1+\alpha_0 r)$, $g_2=\bar\alpha_0 r^2$ (so $\underline{Q}=g_2/g_1$), $A=\alpha_0 r+\bar P$ (the normalized $\gamma_0$), and $D_{\mathrm{int}}:=\alpha_0 A\,g_2-\beta_0 g_1$.
    Recall that we require $L < \underline{Q}$, and we have $L<\underline{Q} \iff D_{\mathrm{int}}>0$.
    As $\Theta$ is convex with vertex $Q_v=-b_1/2\alpha_0^2$, we separately consider $Q_v \geq \underline{Q}$, $Q_v \leq L$, and $Q_v \in (L, \underline{Q})$.
    (If $b_1\ge0$ then $Q_v\le0$ and $\Theta$ is nondecreasing on $[0,\infty)$ with $\Theta\ge b_0>0$, so assume $b_1<0$.)

    First, assume $Q_v \geq \underline{Q}$.
    We have the polynomial identity
    \begin{equation}\label{eq:certA}
      \Theta(\underline{Q})\,g_1^2 \;=\; R_A \;+\; (\bar P+r\bar\alpha_0)\,D_{\mathrm{int}},
    \end{equation}
    where $R_A=\sum_{k=1}^{6}c_k r^k$ has \emph{every} coefficient nonnegative on $0<\alpha_0,P<1$:
    \[ c_1=\alpha_0\bar P^3,\quad c_2=\bar P^2\big(\alpha_0\bar P(\alpha_0{+}1)+\alpha_0^2+3\bar\alpha_0\big),\quad c_3=\bar P\alpha_0\big(\alpha_0\bar P^2+\bar P(\alpha_0^2{+}3{+}\bar\alpha_0)+2\bar\alpha_0\big), \]
    \[ c_4=\bar P\alpha_0^2\big(\bar P(2\alpha_0{+}1)+2+5\bar\alpha_0\big),\quad c_5=\alpha_0^3\big(\bar P(\alpha_0{+}2)+3\bar\alpha_0\big),\quad c_6=\alpha_0^4 \]
    (this can be checked by re-writing $\alpha_0^2{-}3\alpha_0{+}3=\alpha_0^2{+}3\bar\alpha_0$ and $7{-}5\alpha_0=2{+}5\bar\alpha_0$).
    Hence $D_{\mathrm{int}}>0\Rightarrow\Theta(\underline{Q})\ge0$, and when $Q_v\ge\underline{Q}$ the function $\Theta$ is nonincreasing on $[L,\underline{Q}]$, so $M=\Theta(\underline{Q})\ge0$.

    Next, let $Q_v \leq L$.
    If $\beta_0\le0$ then $L=0$ and $M\ge\Theta(0)=b_0>0$.
    If $\beta_0>0$, the identity
    \begin{equation}\label{eq:certC}
      \alpha_0 A^2\,\Theta(L)\;=\;3A^2\beta_0+D_{\mathrm{int}}
    \end{equation}
    gives $\Theta(L)\ge0$ whenever $\beta_0,D_{\mathrm{int}}\ge0$; when $Q_v\le L$, $\Theta$ is nondecreasing on $[L,\underline{Q}]$, so $M=\Theta(L)\ge0$.

    The identities~\eqref{eq:certA},~\eqref{eq:certC} are checked by expanding both sides of the equalities.

    Finally, if $L<Q_v<\underline{Q}$ then $M=\Theta(Q_v)=b_0-b_1^2/4\alpha_0^2$, and $Q_v<\underline{Q}$ gives $-b_1<2\alpha_0^2\underline{Q}$, hence $b_1^2<4\alpha_0^4\underline{Q}^2$; thus $M\ge0$ follows from $\alpha_0^2\underline{Q}^2\le b_0$, that is, from
    \[ T := b_0 g_1^2-\alpha_0^2 g_2^2 \geq 0. \]

    We now show that $T \geq 0$ on the cell $\{L<Q_v<\underline{Q}\}$.
    Write the cell's defining polynomials as the generators
    \[ h=\big(\alpha_0,\ 1-\alpha_0,\ P,\ 1-P,\ r,\ E,\ G,\ D_{\mathrm{int}}\big), \qquad E=2\alpha_0^2 g_2+b_1 g_1,\quad G=-(b_1 A+2\alpha_0\beta_0), \]
    each nonnegative on the cell ($E\ge0$ is $Q_v\le\underline{Q}$, $G\ge0$ is $Q_v\ge L$). The candidate certificate is
    \begin{equation}\label{eq:cert-b}
      T = b_0 g_1^2-\alpha_0^2 g_2^2 \;=:\; \sigma_0+\sigma_E\,E+\sigma_G\,G+\sigma_D\,D_{\mathrm{int}},
    \end{equation}
    each multiplier $\sigma_\bullet$ being a nonnegative-coefficient combination of monomials in $\alpha_0,1-\alpha_0,P,1-P,r$.
    Here $T$ has degree $10$ in the monomial basis $\alpha_0^i P^j r^k$, and $E$, $G$, $D_{\mathrm{int}}$ have degrees $6$, $4$, $7$.
    We fix degree budgets so that every term of~\eqref{eq:cert-b} has degree at most $\deg T = 10$: the multipliers $\sigma_0$, $\sigma_E$, $\sigma_G$, $\sigma_D$ are combinations of the monomials in $\alpha_0, 1-\alpha_0, P, 1-P, r$ of degree at most $10$, $4$, $6$, and $3$, respectively.
    Writing $M_\beta$ for the resulting products, each a monomial in the five linear generators multiplied by one of $1, E, G, D_{\mathrm{int}}$, the certificate asks for weights $\gamma_\beta\ge0$ with $T=\sum_\beta\gamma_\beta M_\beta$.
    Expanding both sides in the monomial basis $\alpha_0^i P^j r^k$ and matching coefficients turns this into a linear feasibility problem
    \begin{equation}\label{eq:cert-LP}
      \text{find } \gamma\ge0 \ \text{ such that } \ \mathbf{A}\gamma=t,
    \end{equation}
    where the columns of the (rational) matrix $\mathbf{A}$ are the coefficient vectors of the $M_\beta$, and $t$ is the (rational) coefficient vector of $T$.
    That is, one asks whether $t$ lies in the cone generated by the columns of $\mathbf{A}$.
    It has $\binom{15}{5}+\binom{9}{5}+\binom{11}{5}+\binom{8}{5} = 3003+126+462+56 = 3647$ variables and $\binom{13}{3} = 286$ equality constraints, one for each monomial of degree at most $10$.
    These counts come from the fact that there are $\binom{d+m}{m}$ monomials of degree at most $d$ in $m$ variables.
    The variables are the monomials in the $m = 5$ linear generators of degree at most $d = 10, 4, 6, 3$, for $\sigma_0$, $\sigma_E$, $\sigma_G$, $\sigma_D$ respectively, and the constraints are the monomials $\alpha_0^iP^jr^k$ in $m = 3$ variables of degree at most $d = 10$.
    This linear program is small: the dual simplex and the interior-point methods of the solver $\mathsf{HiGHS}$ both find a feasible point in under a second.
    Solving the system restricted to the support of the simplex solution exactly over $\mathbb{Q}$ gives an exact rational certificate with $101$ terms, all with positive weights, for which~\eqref{eq:cert-b} holds as an identity of polynomials; it is recorded in Appendix~\ref{app:cellb}.

    Note that feasibility of~\eqref{eq:cert-LP} certifies the inequality, whereas infeasibility certifies only the absence of a certificate of this form with these degree budgets: this is an easy consequence of Farkas' Lemma.

    This shows that $M \geq 0$ when $0 < \Qsat < \underline{Q}$, and concludes the proof of the fourth statement.

    \textbf{Proof of the stability statements. }
    Statement~(a) is Lemma~\ref{lem:linearization}(i), together with the sign of $H'$ on each branch identified in the Proofs of Statements 3 and 4.
    Statement~(b) is Lemma~\ref{lem:linearization}(ii) and~(iii).

    For Statement~(c), we use Lemma~\ref{lem:linearization}(ii).
    As $z = q^*_F \downarrow 0$, we have $w\to0$ and $b\to0$, so $c_2\to\mu_I+\mu_D > 0$ and $c_1c_2 - c_0 \to (\mu_I+\mu_D)\mu_I\mu_D > 0$.
    As $z\to\infty$, by~\eqref{eq:e-kappa} we have $w \leq e\mu_F \to 0$ and $b\,w \leq c\,z\,e\,\mu_F \leq c\,(1-\alpha_0)\Qsat\,\mu_F$, while $\sigma_I,\sigma_D\to\infty$; dropping the nonnegative terms in~\eqref{eq:RH-gap},
    \[
        c_1c_2 - c_0 \geq (\sigma_I+\sigma_D)(\sigma_I-w)(\sigma_D-w) - \mu_F\, b\, w \to \infty,
    \]
    and $c_2\to\infty$.
    An example in which a fixed point with $H'(q^*_F) < 0$ is unstable is given in Section~\ref{sec:hopf}.

    Finally, under~\eqref{eq:rate-condition}, at a fold $q_\dagger\in\{q_{\min},q_{\max}\}$ we have $c_0 = 0$ and, by Lemma~\ref{lem:linearization}(iii), $c_2 > 0$ and $c_1c_2 > 0$, so $\cL$ has a simple zero eigenvalue and two eigenvalues with negative real part.
    Since the folds are simple roots of $H'$, $\cM''(q_\dagger)\neq0$, and the fixed-point set is locally the graph $\lambda = \cM(q_F)$ with a quadratic turning point.
    The two nondegeneracy conditions of the saddle-node bifurcation~\cite{kuznetsov2023elements} also hold.
    Let $f$ be the right-hand side of~\eqref{eq:FL-1}--\eqref{eq:FL-3}, and let $r$ and $\ell$ be right and left null vectors of $\cL$ with $r_3 = 1$ and $\ell\cdot r = 1$.
    The computation in Step~2 of the proof of Theorem~\ref{thm:slow-passage}, which applies verbatim at either fold, gives the transversality
    \[
        \ell\cdot\partial_\lambda f = \ell_1 = \frac{\mu_F\, q_\dagger\,\sigma_I\sigma_D}{\cM(q_\dagger)\, c_1} > 0,
    \]
    since $c_1 > 0$.
    Differentiating the identity $f\big(x^*(q), \cM(q)\big) = 0$ along the curve of fixed points, where $x^*(q)$ is the fixed point with $q^*_F = q$, gives $x^{*\prime}(q_\dagger) = r$ (as $\cM'(q_\dagger) = 0$ and the third component of $x^{*\prime}$ is $1$), and differentiating once more and applying $\ell$, which annihilates the range of $\cL$, gives
    \[
        \ell\cdot D_x^2 f[r, r] = -\cM''(q_\dagger)\;\ell\cdot\partial_\lambda f \neq 0.
    \]
    Hence this is a saddle-node bifurcation, at which a stable fixed point ($H' < 0$) and an unstable one ($H' > 0$) collide and disappear.
    The bistability and hysteresis assertions follow from Statements~(a) and~(b) and the shape of $\cM$ established in the Proofs of Statements 1, 3, and 4.
\end{proof}

\subsection{Proof of the Hopf threshold}\label{sec:proof-hopf}

We prove Proposition~\ref{prop:hopf}, using the notation of Section~\ref{sec:proof-bifurcation}.
For a non-trivial fixed point with $q^*_F = z$, the relation $(1-P)\mu_I q^*_I = \sigma_D q^*_D$ from~\eqref{eq:fp-intact-derelict} gives
\[
    \theta = \frac{\alpha(z)q^*_I}{\alpha(z)q^*_I + q^*_D} = \frac{\alpha(z)\sigma_D}{\alpha(z)\sigma_D + (1-P)\mu_I},
\]
so the quantities $\sigma_I, \sigma_D, b, e, \kappa, \theta$ of Lemma~\ref{lem:linearization} are the functions of $z$ defined in Section~\ref{sec:hopf}, with $b = cz$, and they do not depend on $\mu_F$.

\begin{proof}[Proof of Proposition~\ref{prop:hopf}]
    Fix $z > 0$ with $H'(z) < 0$, and regard the coefficients $c_0, c_1, c_2$ of Lemma~\ref{lem:linearization} as functions of $m := \mu_F$.
    Since $w = e\theta m$,
    \[
        c_2(m) = \sigma_I + \sigma_D - e\theta m, \qquad c_0(m) = m\,C_0, \qquad C_0 := -\sigma_I\sigma_D\, e\theta + b\big((1-\theta)\sigma_I + (1-\kappa)\sigma_D\big),
    \]
    and $C_0 = -\sigma_I\sigma_D\, z H'(z)/H(z) > 0$ by Lemma~\ref{lem:linearization}.
    Substituting $w = e\theta m$ in~\eqref{eq:RH-gap} and collecting powers of $m$ gives
    \[
        p(m) := c_1(m)c_2(m) - c_0(m) = p_0 + p_1 m + p_2 m^2.
    \]

    \emph{The root.}
    Let $m^* := (\sigma_I + \sigma_D)/(e\theta)$, so that $c_2(m) > 0$ if and only if $m < m^*$.
    We have $p(0) = p_0 > 0$, and since $c_2(m^*) = 0$, $p(m^*) = -c_0(m^*) < 0$.
    A polynomial of degree at most two which is positive at $0$ and negative at $m^*$ has exactly one root in $(0, m^*)$, and this root is simple; we call it $\mu_F^H(z)$.
    Thus $p > 0$ on $[0, \mu_F^H(z))$, $p < 0$ on $(\mu_F^H(z), m^*]$, and $p'(\mu_F^H(z)) < 0$.

    \emph{Proof of Statement 1.}
    If $m < \mu_F^H(z)$, then $c_2 > 0$ and $c_1c_2 - c_0 > 0$, and the claim is Lemma~\ref{lem:linearization}(ii).

    \emph{Proof of Statement 2.}
    If $m = \mu_F^H(z)$, then $c_2 > 0$ and $c_1c_2 = c_0 > 0$, so $c_1 > 0$ and
    \[
        \chi(s) = s^3 + c_2s^2 + c_1s + c_1c_2 = (s + c_2)(s^2 + c_1).
    \]

    \emph{Proof of Statement 3.}
    Let $m > \mu_F^H(z)$.
    If $m \leq m^*$ then $c_1c_2 - c_0 < 0$, and if $m > m^*$ then $c_2 < 0$; in both cases, by Lemma~\ref{lem:linearization}(ii), not all eigenvalues have negative real part.
    No eigenvalue lies on the imaginary axis: zero is not an eigenvalue since $c_0 > 0$, and if $\pm i\omega$ with $\omega > 0$ were eigenvalues, then $\chi(s) = (s + r)(s^2 + \omega^2)$ would give $c_2 = r$, $c_1 = \omega^2 > 0$, and $c_0 = c_1c_2$, hence $c_2 = c_0/c_1 > 0$ and $p(m) = 0$ with $m < m^*$, which forces $m = \mu_F^H(z)$.
    The product of the three eigenvalues is $-c_0 < 0$.
    If there is a pair of complex conjugate eigenvalues, their product is positive, so the real eigenvalue is negative, and the pair must then have positive real part.
    If all the eigenvalues are real, then an odd number of them is negative, and not all three are; hence exactly one is negative and two are positive.

    \emph{Proof of~\eqref{eq:hopf-bounds}.}
    The upper bound is $\mu_F^H(z) < m^*$.
    For the lower bound, if $e\,m \leq \min\{\mu_I,\mu_D\}$ then $c_2 > 0$ and $p(m) > 0$ by Lemma~\ref{lem:linearization}(iii), so $m < \mu_F^H(z)$.

    \emph{Transversality.}
    Now fix $\mu_F$ and $\kfrags$ as well, and write $p(z, m)$ to indicate the dependence on $z$.
    Since $\mu_F^H(z)$ is a simple root of $p(z,\cdot)$, it is a smooth function of $z$ on $\{H' < 0\}$ by the implicit function theorem, and differentiating $p(z, \mu_F^H(z)) = 0$ gives
    \begin{equation}\label{eq:hopf-ift}
        \partial_zp(z, \mu_F^H(z)) = -\partial_mp(z, \mu_F^H(z))\,(\mu_F^H)'(z), \qquad \partial_mp(z, \mu_F^H(z)) < 0.
    \end{equation}
    Let $z_H$ be as in the statement.
    By Statement~2, at $z = z_H$ the Jacobian has the simple eigenvalues $-c_2 < 0$ and $\pm i\omega$ with $\omega = \sqrt{c_1}$, so for $z$ near $z_H$ it has a real eigenvalue $-r(z)$ and a pair $x(z) \pm i\omega(z)$, depending smoothly on $z$, with $x(z_H) = 0$.
    Expanding $\chi(s) = (s + r)\big((s - x)^2 + \omega^2\big)$ gives
    \[
        c_2 = r - 2x, \qquad c_1 = x^2 + \omega^2 - 2rx, \qquad c_0 = r(x^2 + \omega^2),
    \]
    and therefore
    \[
        c_1c_2 - c_0 = -2x\big((r - x)^2 + \omega^2\big).
    \]
    Differentiating in $z$ at $z_H$, where $x = 0$, $r = c_2$, and $\omega^2 = c_1$, and using~\eqref{eq:hopf-ift},
    \[
        x'(z_H) = -\frac{\partial_zp(z_H, \mu_F)}{2(c_2^2 + c_1)} = \frac{\partial_mp(z_H,\mu_F)}{2(c_2^2 + c_1)}\,(\mu_F^H)'(z_H),
    \]
    which is non-zero and has the sign of $-(\mu_F^H)'(z_H)$.
    Finally, since $H'(z_H) < 0$, the map $\cM$ is a smooth increasing change of parameter near $z_H$ with $\cM'(z_H) > 0$, so $\tfrac{dx}{d\lambda}(\lambda_H) = x'(z_H)/\cM'(z_H)$ has the same sign.
\end{proof}

\subsection{Proof of the slow passage}\label{sec:slow-passage}
We prove Theorem~\ref{thm:slow-passage}, using the notation of Section~\ref{sec:slow-passage-result} and of Lemma~\ref{lem:linearization}.
Write $x := (q_I, q_D, q_F)$, $f(x, \lambda)$ for the right-hand side of~\eqref{eq:FL-1}--\eqref{eq:FL-3}, and $\mu := \lambda - \lambda_+$.
Under the ramp~\eqref{eq:ramp}, the fluid limit is the autonomous system
\begin{equation}\label{eq:sp-extended}
    \dot x = f(x, \lambda_+ + \mu), \qquad \dot\mu = \varepsilon, \qquad \dot\varepsilon = 0,
\end{equation}
a fast--slow system in which $x$ is fast and $\mu$ is slow.
Its critical manifold, the set of equilibria at $\varepsilon = 0$, consists of the trivial branch and the curve $\{(x^*(z), \cM(z) - \lambda_+) : z > 0\}$ of non-trivial fixed points, parametrized by $z = q^*_F$.
We use geometric singular perturbation theory~\cite{fenichel1979,krupa2001}.
All launch rates involved are at most $\Lambda \leq \mu_I$, so~\eqref{eq:FL-1}--\eqref{eq:FL-3} with $\lambda = \lambda(t)$ is the fluid limit by Remark~\ref{rem:time-varying}.

\begin{proof}[Proof of Theorem~\ref{thm:slow-passage}]
    \emph{Step 1: away from the fold.}
    For $\lambda \in [\lambda_0, \lambda_+ - \delta]$, the fixed point $q^{\mathrm{lo}}(\lambda)$ has $q^*_F$ in a compact subset of $(0, q_{\min})$, on which $H' < 0$.
    By~(A1) and Lemma~\ref{lem:linearization}(ii), all eigenvalues of the Jacobian there have negative real part, uniformly on this compact piece of the branch, which is therefore normally hyperbolic and attracting.
    By Fenichel's theorem~\cite{fenichel1979}, for small $\varepsilon$ there is a locally invariant attracting slow manifold $S^{\mathrm{lo}}_\varepsilon$ within $O(\varepsilon)$ of this piece, and solutions starting in a fixed neighborhood of it approach it at an exponential rate.
    The initial condition $q^{\mathrm{lo}}(\lambda_0)$ is within $O(\varepsilon)$ of $S^{\mathrm{lo}}_\varepsilon$, which gives~(i).

    \emph{Step 2: the center manifold at the fold.}
    At the fold point, $H'(q_{\min}) = 0$, so $c_0 = 0$ by Lemma~\ref{lem:linearization}, and the characteristic polynomial of the Jacobian $\cL$ is $s(s^2 + c_2s + c_1)$.
    By~(A1), $c_2 > 0$ and $c_1c_2 = c_1c_2 - c_0 > 0$, so $c_1 > 0$: zero is a simple eigenvalue of $\cL$, and the other two eigenvalues have negative real part.
    Hence the linearization of~\eqref{eq:sp-extended} at $(q^\dagger, 0, 0)$ has the eigenvalue $0$ with a three-dimensional generalized eigenspace, spanned by the null vector of $\cL$ and the $\mu$ and $\varepsilon$ directions, and two eigenvalues with negative real part.
    By the center manifold theorem~\cite{carr1981}, there is a three-dimensional, locally invariant, attracting center manifold $W^c$ through $(q^\dagger, 0, 0)$, and nearby solutions approach $W^c$ at an exponential rate.

    Let $r$ and $\ell$ be right and left null vectors of $\cL$.
    The first two rows of $\cL$ give $r = \big(-\cL_{[12]}^{-1}(\cL_{13}, \cL_{23})^\top,\ 1\big)$.
    For the left null vector, write $\tilde\ell = (\tilde\ell_1, \tilde\ell_2, 1)$; the second and first components of $\tilde\ell \cL = 0$ read
    \[
        \tilde\ell_2 = \frac{\kfrags\, c\, q_{\min}}{\sigma_D}, \qquad
        \tilde\ell_1 = \frac{(1-P)\mu_I\tilde\ell_2 + \kfrags\,\alpha(q_{\min})c\, q_{\min}}{\sigma_I} = \kfrags\, c\, q_{\min} H(q_{\min}) = \frac{\mu_F\, q_{\min}}{\lambda_+},
    \]
    where we used $\sigma_I H = \alpha + (1-P)\mu_I/\sigma_D$ and $\lambda_+ c\,\kfrags H(q_{\min}) = \mu_F$.
    The third component of $\tilde\ell \cL$ then vanishes as well, since $0 = \tilde\ell \cL r = (\tilde\ell \cL)_3\, r_3$.
    To normalize, note that since $\cL$ has rank two, its adjugate is $\operatorname{adj}(\cL) = \gamma\, r\tilde\ell^\top$ for a scalar $\gamma$.
    Its $(3,3)$ entry is the principal minor $\det \cL_{[12]} = \sigma_I\sigma_D$, so $\gamma = \sigma_I\sigma_D$ (as $r_3 = \tilde\ell_3 = 1$), and its trace is the sum $c_1$ of the principal $2\times 2$ minors of $\cL$, so $\tilde\ell\cdot r = c_1/(\sigma_I\sigma_D) > 0$.
    We set $\ell := \tilde\ell\,\sigma_I\sigma_D/c_1$, so that $\ell\cdot r = 1$ and $\ell_1 = \beta_\lambda$.

    On $W^c$ write $x = q^\dagger + u\, r + h(u, \mu, \varepsilon)$ with $\ell\cdot h = 0$ and $h = O(u^2 + |\mu| + \varepsilon)$.
    Projecting~\eqref{eq:sp-extended} with $\ell$, and using $\ell \cL = 0$ to eliminate the terms which are linear in $h$, gives the reduced system
    \begin{equation}\label{eq:sp-normalform}
        \dot u = \beta_\lambda\,\mu + b_2\, u^2 + O\big(|u|^3 + |u\mu| + \mu^2 + \varepsilon|u| + \varepsilon|\mu| + \varepsilon^2\big), \qquad \dot\mu = \varepsilon,
    \end{equation}
    where $\ell\cdot\partial_\lambda f = \ell_1 = \beta_\lambda$ and $b_2 := \frac12\,\ell\cdot D^2_x f(q^\dagger,\lambda_+)[r,r]$.
    There is no term of order $\varepsilon$ alone, since $f$ does not depend on $\varepsilon$ and the $\varepsilon$-dependence of $h$ enters only through terms which are at least quadratic.
    To identify $b_2$, compare the equilibria of~\eqref{eq:sp-normalform} at $\varepsilon = 0$, which satisfy $\mu = -(b_2/\beta_\lambda)u^2 + O(|u|^3)$, with the non-trivial fixed points near the fold, which satisfy $\mu = \cM(q_F) - \lambda_+$.
    Since $r_3 = 1$ and $h = O(u^2 + |\mu|)$, we have $q_F = q_{\min} + u + O(u^2)$ along them, and so $\mu = \tfrac12\cM''(q_{\min})u^2 + O(|u|^3)$.
    Hence $b_2 = -\tfrac12\beta_\lambda\cM''(q_{\min})$, and since $H'(q_{\min}) = 0$,
    \[
        \cM''(q_{\min}) = -\frac{\mu_F H''(q_{\min})}{c\,\kfrags H(q_{\min})^2} = -\frac{\lambda_+ H''(q_{\min})}{H(q_{\min})},
    \]
    so $b_2 = \beta_2$ as defined in~\eqref{eq:sp-coefficients}, and $\beta_2 > 0$ because $H''(q_{\min}) > 0$.

    \emph{Step 3: passage through the fold.}
    The substitution $X := \beta_2 u$, $Y := -\beta_\lambda\beta_2\,\mu$, $\tilde\varepsilon := \beta_\lambda\beta_2\,\varepsilon$ transforms~\eqref{eq:sp-normalform} into
    \[
        \dot X = X^2 - Y + O\big(|X|^3 + |XY| + Y^2 + \tilde\varepsilon|X| + \tilde\varepsilon|Y| + \tilde\varepsilon^2\big), \qquad \dot Y = -\tilde\varepsilon,
    \]
    which is the canonical form of a generic fold point of a planar fast--slow system~\cite[Section~2]{krupa2001}.
    Its attracting critical branch $X = -\sqrt{Y}$, $Y > 0$, is the low-debris branch ($u < 0$, $\mu < 0$), and the slow manifold $S^{\mathrm{lo}}_\varepsilon$ of Step~1 continues into $W^c$ as the attracting slow manifold of this planar system.
    By~\cite[Theorem~2.1]{krupa2001} and its proof (see also~\cite{mishchenko1980} for the asymptotic expansion), the continuation of the attracting slow manifold through the fold passes $Y = 0$ at $X = O(\tilde\varepsilon^{1/3})$ and reaches the section $\{X = \tilde\Delta\}$, for fixed small $\tilde\Delta > 0$, at
    \[
        Y = -\Omega_0\,\tilde\varepsilon^{2/3} + O\big(\tilde\varepsilon\log(1/\tilde\varepsilon)\big),
    \]
    and that nearby solutions are attracted to it at an exponential rate.
    The constant $\Omega_0$ comes from the scaling $X = \tilde\varepsilon^{1/3}U$, $Y = -\tilde\varepsilon^{2/3}s$, under which the leading-order equation becomes the Riccati equation $dU/ds = U^2 + s$; its solution which follows the attracting branch $U \sim -\sqrt{-s}$ as $s \to -\infty$ is $U(s) = \mathrm{Ai}'(-s)/\mathrm{Ai}(-s)$, which is finite at $s = 0$ and blows up at $s = \Omega_0$.
  In this scaling, the continuation of the attracting slow manifold follows $U$ up to $o(1)$, uniformly on compact subsets of $(-\infty, \Omega_0)$, by the analysis of the rescaled chart in the proof of~\cite[Theorem~2.1]{krupa2001}.

    By Step~1 and the attractivity of $W^c$, the solution $q$ enters the fold region exponentially close to $S^{\mathrm{lo}}_\varepsilon$, and so exponentially close to $W^c$, and the estimates above apply to it.
    The inner variable $s$ is time on the scale $\varepsilon^{-1/3}$: under the ramp~\eqref{eq:ramp}, $\mu = \lambda(t) - \lambda_+ = \varepsilon(t - t_+)$ with $t_+ := (\lambda_+ - \lambda_0)/\varepsilon$, so that
    \begin{equation}\label{eq:sp-inner-time}
        t - t_+ = \frac{\mu}{\varepsilon} = s\,(\beta_\lambda\beta_2)^{-1/3}\,\varepsilon^{-1/3},
    \end{equation}
    and the blow-up of $U$ at $s = \Omega_0$ is reached a time $\Omega_0(\beta_\lambda\beta_2)^{-1/3}\varepsilon^{-1/3}$ after the launch rate crosses $\lambda_+$.
    Returning to $u = X/\beta_2$ and $\mu = -Y/(\beta_\lambda\beta_2) = s\,(\beta_\lambda\beta_2)^{-1/3}\varepsilon^{2/3}$: the range of $\lambda(t)$ in~(ii) corresponds to $s$ in the compact interval $[-K(\beta_\lambda\beta_2)^{1/3}, \Omega_0 - \delta]$, on which $U$ is bounded, so there $u = O(\varepsilon^{1/3})$ and $\mu = O(\varepsilon^{2/3})$, hence $|q(t) - q^\dagger| \leq |u|\,|r| + |h| = O(\varepsilon^{1/3})$, which is~(ii); and $u$ reaches $\tilde\Delta/\beta_2$ when
    \[
        \mu = \Omega_0\,(\beta_\lambda\beta_2)^{-1/3}\,\varepsilon^{2/3} + O\big(\varepsilon\log(1/\varepsilon)\big).
    \]
    Along $W^c$, $q_F = q_{\min} + u + O(u^2 + |\mu| + \varepsilon)$, and $\dot u \geq \beta_2 u^2/2 > 0$ for $u$ small and positive and $\mu, \varepsilon$ small.
    Choosing $\tilde\Delta$ so that $\tilde\Delta/\beta_2$ is smaller than $\Delta/2$, the solution reaches $q_F = q_{\min} + \Delta$ within a time of order one after $u = \tilde\Delta/\beta_2$, during which $\mu$ changes by $O(\varepsilon)$, and $q_F < q_{\min} + \Delta$ before then.
    This gives~(iii).
    In time: by~\eqref{eq:sp-inner-time}, the range of $s$ used for~(ii) is the range $-K\varepsilon^{-1/3} \leq t - t_+ \leq (\Omega_0 - \delta)(\beta_\lambda\beta_2)^{-1/3}\varepsilon^{-1/3}$ of times, which is the time form of~(ii).
    For~(iii), $u$ reaches $\tilde\Delta/\beta_2$ at the time $t - t_+ = \mu/\varepsilon = \Omega_0(\beta_\lambda\beta_2)^{-1/3}\varepsilon^{-1/3} + O(\log(1/\varepsilon))$, and $q_F$ reaches $q_{\min} + \Delta$ a time of order one later, which is absorbed in the $O(\log(1/\varepsilon))$ term; this is~\eqref{eq:sp-law-time}.

    \emph{Step 4: the jump.}
    At $\mu = \varepsilon = 0$, the set $W^c \cap \{\mu = \varepsilon = 0\}$ is the center manifold of $q^\dagger$ for the system frozen at $\lambda_+$, on which $\dot u = \beta_2 u^2 + O(|u|^3)$; its half $\{u > 0\}$ is the orbit of~(A2).
    At time $t_{\Delta}$ we have $\mu, \varepsilon = O(\varepsilon^{2/3})$, so $q(t_{\Delta})$ is within $O(\varepsilon^{2/3})$ of a point $p$ on this orbit.
    By~(A2), the frozen orbit through $p$ enters, after a time $T_1$ which does not depend on $\varepsilon$, a neighborhood $U$ of $q^{\mathrm{hi}}(\lambda_+)$ which is contained in the basin of attraction of the slow manifold described next.
    For $\lambda \in [\lambda_+ - \delta', \Lambda]$ with $\delta' > 0$ small, the fixed points $q^{\mathrm{hi}}(\lambda)$ have $q^*_F > q_{\max}$, where $H' < 0$, so by~(A1) and Lemma~\ref{lem:linearization}(ii) this piece of the high-debris branch is normally hyperbolic and attracting, and Fenichel's theorem gives an attracting slow manifold $S^{\mathrm{hi}}_\varepsilon$ within $O(\varepsilon)$ of it.
    By continuous dependence on initial conditions and parameters over the time interval $[t_{\Delta}, t_{\Delta} + T_1]$, during which $\lambda$ changes by $\varepsilon T_1$, the solution $q$ enters $U$.
    It then approaches $S^{\mathrm{hi}}_\varepsilon$ at an exponential rate, so it is within $O(\varepsilon)$ of $S^{\mathrm{hi}}_\varepsilon$, and hence of $q^{\mathrm{hi}}(\lambda(t))$, after a further time of order $\log(1/\varepsilon)$, during which $\lambda$ changes by $O(\varepsilon\log(1/\varepsilon))$.
    For $\varepsilon$ small, all of this happens before $\lambda(t)$ reaches $\lambda_+ + \delta$, and the solution then follows $S^{\mathrm{hi}}_\varepsilon$ until $\lambda(t) = \Lambda$.
    This gives~(iv).
\end{proof}

\input{admission-control}

\input{appendix-cellb}

\end{document}

%% file: math.tex
\numberwithin{equation}{section}

\newtheorem{theorem}{Theorem}[section]
\newtheorem{lemma}[theorem]{Lemma}

\newtheorem{corollary}[theorem]{Corollary}
\newtheorem{proposition}[theorem]{Proposition}

\theoremstyle{remark}
\newtheorem{definition}[theorem]{Definition}

\newtheorem{remark}{Remark}

\renewcommand{\le}{\leqslant}
\renewcommand{\ge}{\geqslant}
\renewcommand{\leq}{\leqslant}
\renewcommand{\geq}{\geqslant}

\newcommand{\ind}{\mathbf{1}}

\renewcommand{\iff}{\Leftrightarrow}

\newcommand{\cI}{\mathcal{I}}
\newcommand{\cJ}{\mathcal{J}}\newcommand{\cL}{\mathcal{L}}
\newcommand{\cM}{\mathcal{M}}

\newcommand{\bN}{\mathbb{N}}
\newcommand{\bR}{\mathbb{R}}

\newcommand{\bZ}{\mathbb{Z}}

\newcommand{\Poi}{\mathrm{Poi}}
\newcommand{\Exp}{\mathrm{Exp}}

%% file: intro.tex
\section{Introduction}
Satellite technology has quietly become critical infrastructure for the modern
economy.
The global space economy
reached a record \$686~billion in 2025, up 12\% year over year~\cite{spacefoundation2026},
and is projected to exceed \$1~trillion by the early 2030s; the commercial satellite
sector alone accounts for \$303~billion per year~\cite{sia2025}.
No system embodies this economic impact better than the Global Positioning
System (GPS).
GPS, whose satellites orbit in Medium Earth Orbit (MEO) at an altitude of 20{,}200~km, has
generated an estimated \$1.4~trillion in U.S.\ private-sector economic benefits since
civilian access in the 1980s, roughly \$300~billion per year by 2017.
Indeed, any service
interruption to GPS is estimated to cost roughly \$1~billion per day~\cite{rti2019gps}.
Satellite television and satellite communication (SATCOM) services, delivered from Geostationary Orbit (GEO), 35{,}786~km above the Earth's surface, generated \$108~billion in service revenue in 2024~\cite{sia2024}, with satellite TV alone accounting for \$72~billion.

SATCOM services from GEO or MEO are latency limited, however. For instance, a GEO link's round-trip signal
path is 71{,}600~km, so propagation alone imposes a latency floor on the order of
240~ms, and 500--700~ms in practice. This is high enough to degrade interactive and
latency-sensitive services. SpaceX's deployment of Starlink in \emph{Low Earth Orbit
(LEO)}, the altitude band between 160~km and 2{,}000~km, breaks this physical constraint by cutting the
round-trip latency to 25--60~ms, an order of magnitude below GEO. The high quality of service offered by Starlink is reflected in the firm's rapid growth to
\$11.4~billion in revenue and more than 10~million subscribers in 2025. This is the first
satellite internet service to clear \$10~billion in a single year~\cite{spacex2026s1}.
Starlink has thus demonstrated that high-speed satellite connectivity in LEO is not merely possible but commercially transformative.

In LEO, groups of satellites owned and operated by the same entity, or {\it operator}, are called \emph{constellations}.
These constellations can act as a networked system to provide basic infrastructure, including, for example, SpaceX's Starlink and Amazon's Leo high-speed internet services. As of October 2026, SpaceX's Starlink constellation has more than 11{,}000 satellites, and SpaceX has filed applications for up to 42{,}000 satellites~\cite{orbitalradar2026starlink}.
Amazon's Leo constellation has more than 390 satellites in orbit, and the U.S.\ Federal Communications Commission has authorized Amazon to launch and operate up to 7{,}736 satellites~\cite{orbitalradar2026kuiper}.
The Chinese operators Guowang and Qianfan also have plans for constellations with a total of over 28{,}000 satellites~\cite{orbitalradar2026china}.
Recently, there have also been proposals to launch LEO constellations to serve as AI datacenters~\cite{stanisic2026constella, ohs2026dark, aili2025development}, including Google's experimental Suncatcher project~\cite{ruwitch2026suncatcher}.

LEO constellations are also considered to be {\it resilient} infrastructure, since service is delivered collectively, and denial of service requires
destroying or disabling a large fraction of the satellites, or severing their
links to terrestrial networking infrastructure.
This resilience also carries national security implications, since governments increasingly rely on proliferated LEO systems for communications, missile warning, and
tracking~\cite{gao2026missile}. Since {resiliency}  requires a large number of satellites, launch rates to LEO are also increasing concomitantly. For example, SpaceX now launches a new batch of Starlink satellites every few days~\cite{orbitalradar2026starlink}.

\begin{figure}[ht]
    \centering
    \includegraphics[width=0.85\linewidth]{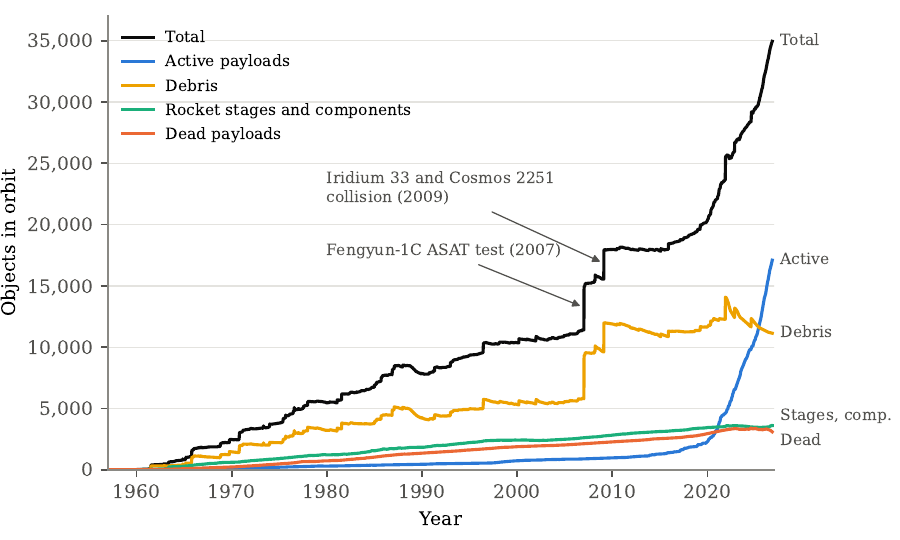}
    \caption{Cataloged objects in Earth orbit, 10~cm or larger in diameter, by type, from 1957 to October 26, 2026, from the satellite and debris population statistics of McDowell~\cite{mcdowell2026population}. The jumps in the debris count in 2007 and 2009 are the Fengyun-1C anti-satellite test and the Iridium~33--Cosmos~2251 collision.}
    \label{fig:tracked-objects}
\end{figure}

These constellations share LEO with a large amount of man-made debris built up over
decades of space operations (Figure~\ref{fig:tracked-objects}).
This includes not only debris left behind from legacy rocket launches, but also
debris from deliberate actions and from accidental collisions.
For example, in 2007, the Chinese armed forces tested an anti-satellite missile on its Fengyun-1C
satellite in the 800--1000~km altitude band, which created more than 3,500 debris
pieces, of which about 2,100 are still in orbit in
2026~\cite{mcdowell2026gcat}.
Two years later, on February 10, 2009, the defunct Russian satellite Cosmos~2251
collided with the operational Iridium~33 at an altitude of about 790~km, the
first accidental catastrophic collision between ``derelict'' and operational
satellites~\cite{nasa-odqn-2009}; the collision produced more than 2,300
cataloged fragments, of which about 700 are still in orbit in
2026~\cite{mcdowell2026gcat}.
These debris pieces have an operational impact.
In November 2024, the International Space Station performed a collision
avoidance maneuver to avoid debris from the Fengyun-1C
satellite~\cite{odqn2025iss}.
Starlink is also lowering many of its satellites with the purpose of ``reducing
the aggregate likelihood of collision'' since there are fewer debris pieces at lower
altitudes~\cite{reuters2026starlink}.
However, this can also lower the operational lifetime of satellites in the constellation, as at lower altitude they are subject to greater atmospheric drag.
With increasing satellite density, the likelihood of a catastrophic collision
with debris increases.

Debris cannot be controlled, and its trajectories are difficult to measure and predict.
Earth-based detection typically only measures debris that are 10~cm or larger in diameter.
At present, there are approximately 10,000 tracked debris pieces, most of them in
LEO~\cite{orbitalradar2026debris}.
However, smaller pieces of debris can still be energetic enough to cause a
catastrophic collision with a satellite.
Depending on the altitude, debris pieces can orbit for tens to hundreds, or even
thousands of years before de-orbiting due to atmospheric drag effects.
A collision between debris and a satellite (es exemplified by the Iridium-Cosmos collision) can be catastrophic and produce
thousands of debris pieces at once, which in turn raise the likelihood of further
collisions; this positive feedback is known as the \emph{Kessler
Syndrome}~\cite{kessler1978}.
The Kessler Syndrome has also entered popular media. For instance, in the Hollywood film \emph{Gravity}, the following is issued as a status update to astronauts shortly before debris destroys the space shuttle they are on: ``Debris chain reaction is out of control and rapidly expanding. Multiple satellites are down and they keep on falling.''
\emph{Debris management} is therefore a fundamental challenge to large-scale operations in LEO. It will be useful to distinguish between two types of debris, viz. {\it fragments} and {\it derelicts}, i.e., a satellite that has completed its mission but has not de-orbited.
Going forward, we will use the words ``intact'' as short form for ``intact satellite'', and ``fragment'' for ``debris fragment.''
We refer to an altitude band within LEO as a {``shell''}.

Most of the solution strategies to handle the challenges posed by large debris counts in LEO fall into three categories: \emph{post-mission disposal (PMD) guidelines}, \emph{active debris removal}, and \emph{operating in the presence of debris}.
PMD guidelines are recommendations that satellites should be removed from orbit within a certain time interval (\textit{e.g.}, 5 years) after their stated mission completion.
These guidelines are designed to prevent future pollution of LEO.
However, PMD only focuses on the prevention of new debris in LEO, instead of on the numerous debris particles which already exist.
Active debris removal is the process of destroying or otherwise removing debris fragments or derelicts from space.
However, active debris removal technologies are largely untested, and NASA
describes them as either ``complex and costly to use,'' or as involving United
Nations treaty violations~\cite{nasa-appel-debris}.
As of 2026, no piece of debris has yet been removed from
orbit~\cite{orbitalradar2026adr}, and the missions closest to doing so each
target one or a few large objects.
Astroscale's ADRAS-J inspected a spent Japanese H-IIA rocket upper stage at close range
in 2024, and its successor ADRAS-J2 is scheduled to capture and deorbit the same
stage, with launch planned for Japan's 2027 fiscal
year~\cite{payload2026adrasj2}.
Astroscale's ELSA-M is designed to capture defunct satellites that were fitted
with magnetic docking plates before launch, and is scheduled to launch in
Astroscale's 2028 fiscal year~\cite{payload2026adrasj2}.
The European Space Agency's ClearSpace-1 will attempt to remove a single
derelict, with launch planned for 2028.
Each of these missions costs on the order of \$100 million, and
none targets fragments, which are far more numerous and much more difficult to capture.

In this paper, we instead focus on the third option of \emph{operating in the presence of debris}, in which intact satellites maneuver to avoid \emph{conjunctions}, or close approaches, with tracked fragments.
Avoidance relies on tracking, and tracking is limited.
For example, SpaceX has developed a satellite-based measurement and debris tracking system called Stargaze~\cite{spacex2026stargaze}, and shares conjunction screening data with other satellite operators free of charge. Starlink uses an automated collision avoidance system that plans risk-mitigating maneuvers without humans in the loop, and satellites assume maneuver responsibility by default~\cite{starlink2026conops}.
However, onboard compute is a limited resource~\cite{onorbitai2026}, and estimating each
object's trajectory, associating observations across multiple observers, and
cataloging the thousands of fragments created by a catastrophic collision are
difficult to do in real time. We therefore model collision avoidance as \emph{saturating}; \textit{i.e.}, as the number of fragments grows, the fraction of conjunctions a satellite can avoid shrinks.

Military and economic uses of LEO depend on large numbers of intacts, and therefore on high launch rates.
However, if the launch rate is too high, the collision rate can increase and the process  becomes self-sustaining, with each catastrophic collision producing
fragments that cause further collisions~\cite{kessler1978}.
Thus, there is a tradeoff between launch rates for new intacts, which we wish to maximize, and collision rates with debris, which we wish to minimize.
In this paper, we assume that intacts' collision avoidance is such that they can successfully avoid other intacts and derelicts, so that all collisions are with fragments.
Thus, the collision rate is driven by the fragment count.
We call the largest launch rate a shell can sustain while keeping its fragment
count steady at a low level (\textit{i.e.}, its equilibrium fragment count) its \emph{carrying capacity}.
The tradeoff is collective among constellation operators, and a potential tragedy of the commons, since all operators in a shell share its debris environment.
Each operator's launches raise the collision risk for everyone.

This paper addresses how the equilibrium fragment count (and hence LEO constellation operations) is shaped by the overall launch rate and the collision avoidance capability of the intacts. We study a single shell, assuming that the flux of intacts,
derelicts, and fragments between shells is negligible.
Within the shell, launches of intacts and de-orbits of all three object
classes determine how the intact,
derelict, and fragment counts vary dynamically over time.
Collision avoidance, which only intacts can perform, adds a second
layer of dynamics that is coupled to the fragment count.
Consequently, we propose a stylized two-layer queueing network model of the shell.
The first layer describes the collision avoidance.
The second layer describes the evolution of the object populations: launches of
intacts, de-orbits of intacts, derelicts, and fragments, and transitions between object
classes (intact to derelict, intact to fragment, and derelict to fragment).

The full model is in Section~\ref{sec:model}, but an overview is as follows, starting with the first layer.
For each intact-fragment pair, a potential collision occurs at rate $c$.
Intacts can perform collision avoidance, but that capability saturates as the number of fragments grows.
We make the simplifying assumption that each intact (or derelict) experiences conjunctions with each fragment at the same rate.
Each intact handles conjunctions as a processor-sharing queue, with its avoidance capacity $\mu_s$ is shared equally among the conjunctions in progress.
Each conjunction has a deadline, namely the time to closest approach, and a conjunction which is not processed before its deadline is an avoidance failure (specifically, the deadlines are distributed as $\Exp(\nu)$).
If a conjunction avoidance failure occurs, the intact suffers a catastrophic collision.
Suppose the instantaneous fragment count is $q$. Then the instantaneous collision avoidance \emph{failure} rate is $\alpha(q)cq$, where $$\alpha(q) = \alpha_0 + (1-\alpha_0)\,q/(q + \Qsat)$$ gives the collision avoidance failure probability, $\alpha_0$ is a baseline failure probability, and $\Qsat = \mu_s/\nu$ is a unitless parameter which measures collision avoidance capability, with larger values indicating greater capability. That is, the instantaneous rate at which an intact becomes a fragment is $\alpha(q)cq$.
In the full model, we assume each catastrophic collision generates a Poisson-distributed batch of new fragments.

For the second layer, we use a queueing network with three queues corresponding to the intacts, derelicts, and fragments.
The intact queue is an $M/\cdot/n/n$ loss system with launch rate $\lambda$, where $n$ is essentially the number of available orbital slots in this shell; the derelict and fragment queues are $\cdot/M/\infty$ infinite server systems.
Denote by $(Q_I, Q_D, Q_F)$ the process consisting of the intact, derelict, and fragment counts, respectively.
Each intact routes to the fragment queue at rate $\alpha(Q_F)cQ_F$.
On the other hand, an intact finishes its mission without catastrophic collision at rate $\mu_I$, and with probability $P$ it departs the system; otherwise it routes to the derelict queue.
Each derelict undergoes catastrophic collision at rate $cQ_F$ and de-orbits at rate $\mu_D$.
Fragments de-orbit at rate $\mu_F$.

As it is the fragment count that couples the layers of the queueing network, we study how the fragment count varies with the launch rate.
Exact constellation-scale analysis is difficult and, therefore, we study a \emph{fluid limit} ordinary differential equation (ODE) approximation of the intact, derelict and fragment counts, in the limit of a large number of orbital slots.
The fluid limit is established as a functional strong law of large numbers.
In the fluid limit, the equilibrium intact and derelict counts are determined as functions of the equilibrium fragment counts. The existing literature on LEO capacity has largely focused on particles-in-a-box (PIB) models which treat satellites, derelicts, and debris like gas particles in a closed box (similar to thermodynamic analyses). Most of these model the counts as solutions of ODEs~\cite{jang2022stability, lavezzi2025stable, lifson2022many, bradley2009space}, though recent work has also developed stochastic differential equation (SDE) models in large-population settings \cite{behera2026stochastic} and in large-volume scaling limits~\cite{choi2026narrow}.
These models typically treat avoidance as a fixed parameter. Consequently, for any sufficiently small launch rate, there is a unique fixed point for the intact, derelict, and fragment counts, i.e., they are {\it monostable}. On the other hand, when launch rates are large the fragment counts are unbounded in these models.

Incorporating collision avoidance fundamentally changes the structure of the fluid limit: it becomes \emph{bistable}, so that there can be two stable fixed points for the intact, derelict, and fragment counts over some ranges of launch rates.
Moreover, the equilibrium fragment count in our model is finite at every launch rate, so the danger is not runaway growth of the fragment count, but a jump from a low-debris to a high-debris equilibrium.

\begin{figure}[ht!]
    \centering
    \begin{subfigure}[b]{0.8\textwidth}
        \centering
        \includegraphics[width=\textwidth]{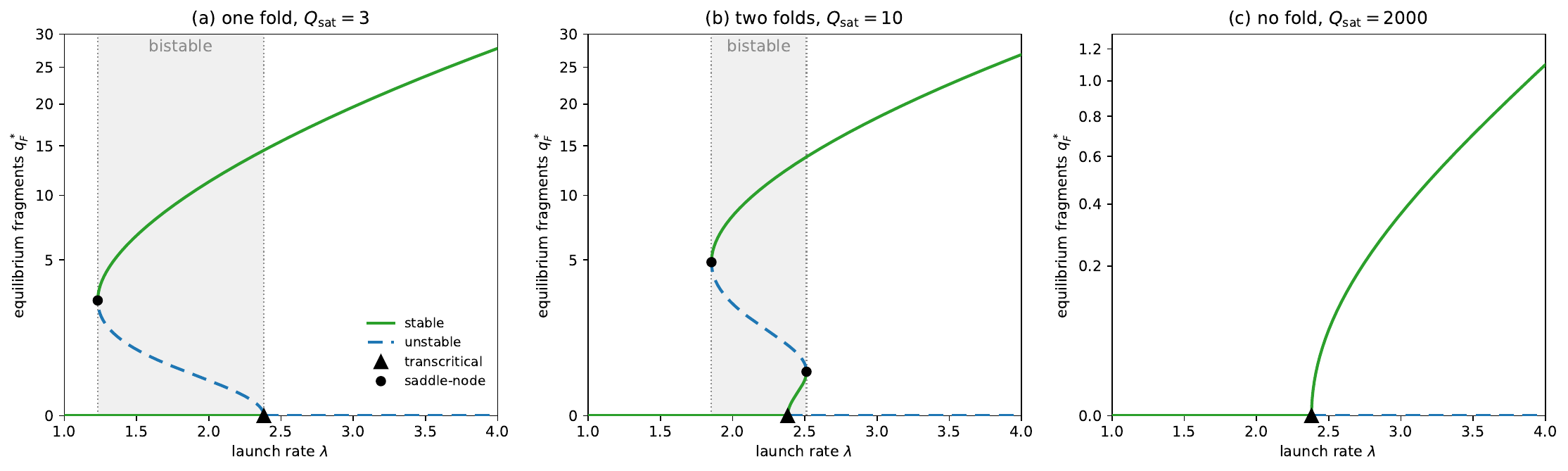}
        \caption{Equilibrium fragment counts as a function of launch rate, with increasing collision avoidance capability $\Qsat$.}
        \label{fig:fragment-bistability}
    \end{subfigure}

    \vspace{0.5cm} %

    \begin{subfigure}[b]{0.8\textwidth}
        \centering
        \includegraphics[width=\textwidth]{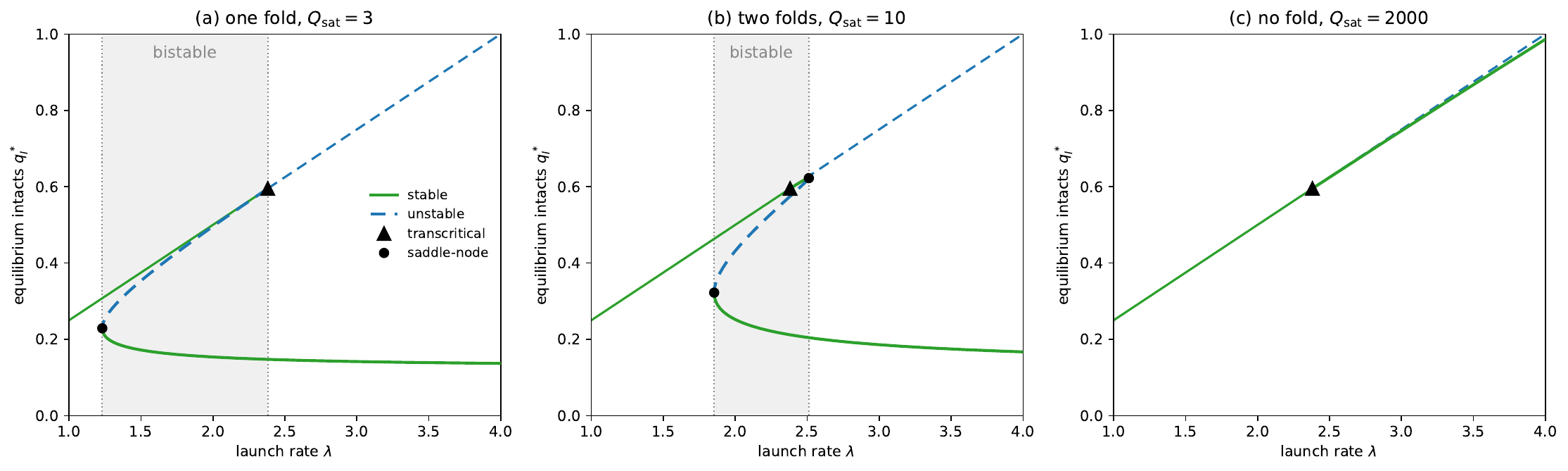}
        \caption{Equilibrium intact counts as a function of launch rate, with increasing collision avoidance capability $\Qsat$.}
        \label{fig:intact-bistability}
    \end{subfigure}

    \caption{Different paths to bistability as a function of $\Qsat$, as shown in Theorem~\ref{thm:bifurcation-structure}.
    Here, we have $\mu_I = 4, \mu_D = \mu_F = c = 1, \alpha_0 = 0.05$, $P = 0.96$ and $\kfrags = 8$.}
    \label{fig:bistability}
\end{figure}

Our first set of results characterizes when the fluid limit admits
multiple stable operating states, as a function of the avoidance
capacity $\Qsat$; Figure~\ref{fig:bistability} illustrates these outcomes.
The bistability in the debris count (Figure~\ref{fig:bistability} (a)) is  reflected in the intact count (Figure~\ref{fig:bistability}(b)), with the high debris equilibrium corresponding to a low intact equilibrium, and vice versa.
Bistability, with one stable state at a low fragment count and another at a
high fragment count, can arise through exactly two bifurcation structures,
involving one or two folds respectively~\cite{kuznetsov2023elements}.
A \emph{fold} (or saddle-node) bifurcation occurs when, as a parameter varies, a
stable and an unstable equilibrium merge and disappear, so that the long-run
behavior of the system changes abruptly.
Figure~\ref{fig:fragment-bistability} illustrates both structures.
With one fold (panel~(a)), the range of launch rates that admits two stable
states is bounded by the fold on one side and by a transcritical bifurcation on
the other; with two folds (panel~(b)), it is bounded by the two folds.
A \emph{branch} is a continuous curve of equilibria as the launch rate is varied.
When $\Qsat$ is large, bistability disappears (panel~(c)).
Instead, the fragment-free state loses stability at a critical launch rate $\lambda_\crit$, and a
branch of equilibria with positive fragment counts emerges (a transcritical
bifurcation).
This branch is stable under an explicit rate condition (Theorem~\ref{thm:bifurcation-structure}, stability condition (b)).
When the condition fails, stability can be lost through a Hopf bifurcation, after
which the fragment count oscillates; see Figure~\ref{fig:hopf} below.
In this regime, a low fragment count may be only one phase of a cycle rather than
a safe operating point.

In this paper, we specifically focus on the bistable regime, about which we now provide more detail.
For small launch rates, there is a stable branch called the \emph{low debris branch} (which could include the state with no fragments), and for larger launch rates, there is a stable branch called the \emph{high debris branch}.
Bistability is when these two branches coexist.
In the bistable regime, the equilibrium behavior of the fluid limit is \emph{hysteretic}.
That is, when the two stable branches coexist for some launch rate, the system has memory of whether it entered the bistable window via the low debris branch or via the high debris branch.
If it entered via the high debris branch, the only way to reach the low debris branch is to decrease the launch rate until the only stable state is the low debris branch (and vice versa).
For example, in Figure~\ref{fig:fragment-bistability} (b), if the launch rate is $2$ and we are on the high debris branch, we need to reduce the launch rate below about $1.85$ to reach the low debris branch.
Because infrastructural resilience for LEO constellations requires higher intact counts, it is important to understand the transition from the low debris branch to the high debris branch, and to develop strategies to remove the bistability.

\begin{figure}[ht]
    \centering
    \includegraphics[width=0.6\linewidth]{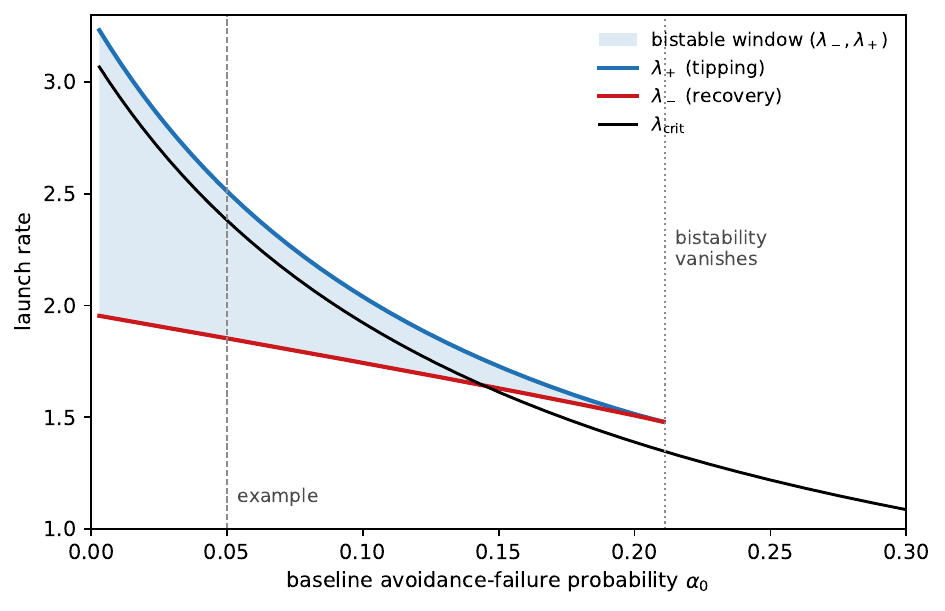}
    \caption{Launch-rate thresholds against the baseline failure probability $\alpha_0$, with $\Qsat = 10$ and the other parameters as in Figure~\ref{fig:bistability}. The shaded band is the bistable window $(\lambda_-, \lambda_+)$, which closes at $\alpha_0 \approx 0.21$.}
    \label{fig:alpha0-thresholds}
\end{figure}

Better maneuver reliability cuts both ways.
Figure~\ref{fig:alpha0-thresholds} shows how the launch-rate thresholds move with the baseline failure probability $\alpha_0$ of collision avoidance, with $\Qsat = 10$ and the other parameters as in Figure~\ref{fig:bistability}.
More reliable maneuvers, with smaller $\alpha_0$, raise the tipping point $\lambda_+$, so the shell can sustain a higher launch rate on its low-debris branch: in the example, $\lambda_+$ rises from $2.04$ at $\alpha_0 = 0.1$ to $2.51$ at $\alpha_0 = 0.05$ and $3.19$ at $\alpha_0 = 0.005$.
But they raise the recovery threshold $\lambda_-$ much less, so the bistable window widens, and the launch-rate cut needed to leave the high-debris branch, $1 - \lambda_-/\lambda_+$, grows from $15\%$ to $26\%$ to $39\%$; the jump at the tipping point is also larger.
Conversely, the bistability disappears altogether if the baseline failure rate is high; for example, above $\alpha_0 \approx 0.21$ in Figure~\ref{fig:alpha0-thresholds}.
The same tension appears in Theorem~\ref{thm:bifurcation-structure}: no fold can occur when the avoidance capability exceeds $\overline{Q} = (1-\alpha_0)\mu_I/(c\,\alpha_0^2)$, which grows like $\alpha_0^{-2}$ as maneuvers become more reliable.

\begin{figure}[ht]
    \centering
    \includegraphics[width=0.8\linewidth]{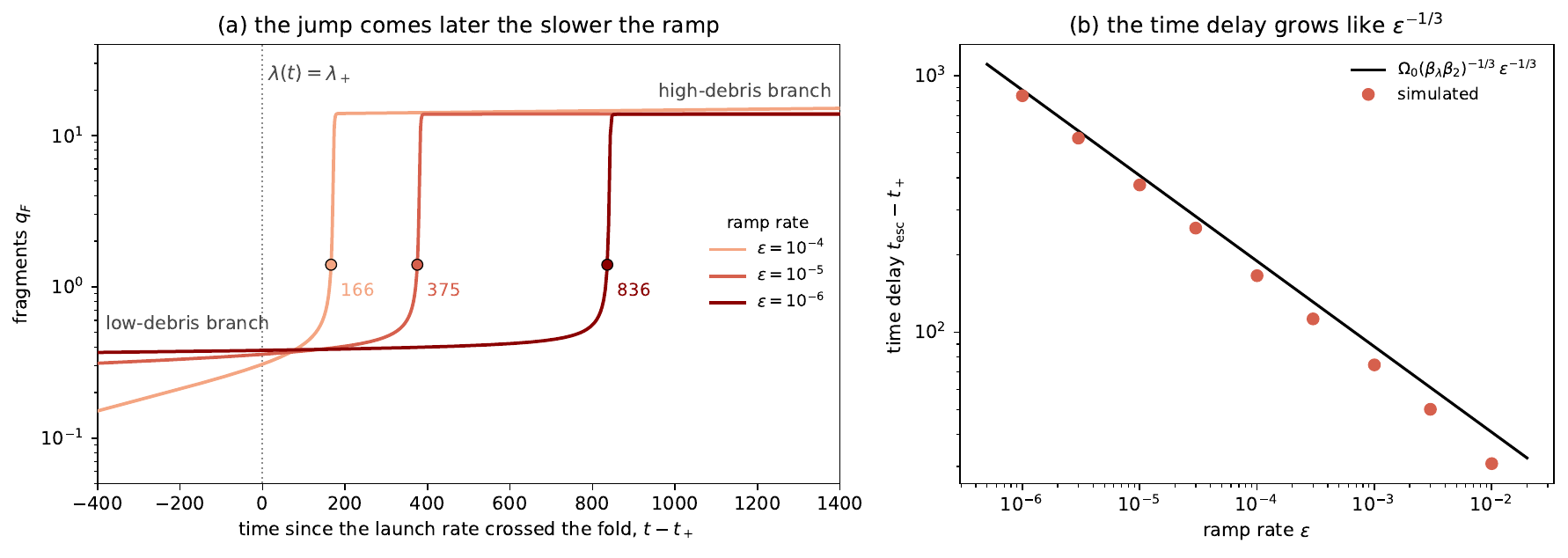}
    \caption{Slow passage through the upper fold in the two-fold regime. Here, the system crosses the fold at time $t_+$, and escapes to the high debris branch at time $t_{\mathrm{esc}}$. (a)~The fragment count against $t - t_+$. (b)~The time delay $t_{\mathrm{esc}} - t_+$ against the ramp rate, compared with the leading-order law $\Omega_0(\beta_\lambda\beta_2)^{-1/3}\varepsilon^{-1/3}$ implied by Theorem~\ref{thm:slow-passage}.}
    \label{fig:slow-passage}
\end{figure}

A distinctive feature of the bistable fluid limit is how the system moves from the
low-debris branch to the high-debris branch as the launch rate increases.
When the launch rate passes the upper end of the \emph{bistable interval}, the range of
launch rates at which both branches are stable, the low-debris state ceases to
exist, but the system does not leave it immediately: it lingers near the vanished
branch before jumping to the high-debris branch.
This is a \emph{slow passage through the fold}~\cite{kuehn2015multiple}.
Specifically, if the launch rate is increased at a small constant ramp rate
$\varepsilon$, starting from the low-debris equilibrium, we show that the delay
between crossing the fold and the departure of the intact, derelict, and fragment
counts to the high-debris branch is
\[
  \Delta t \;\sim\; \Omega_0\,(\beta_\lambda\,\beta_2)^{-1/3}\,\varepsilon^{-1/3}
  \qquad (\varepsilon \to 0),
\]
with $\Omega_0$ the smallest positive zero of $s \mapsto \mathrm{Ai}(-s)$, where $\mathrm{Ai}$ is the Airy function, and $\beta_\lambda$, $\beta_2$ model constants
(Theorem~\ref{thm:slow-passage}; see Figure~\ref{fig:slow-passage}).

This makes the fragment count a \emph{lagging} indicator of an unsustainable
launch rate. A growing shell can be assumed to cross
the upper fold $\lambda_+$ (the launch rate at the upper end of the bistable
window) at a finite ramp rate; for example, annual satellite launches rose
from 2{,}695 in 2024 to 4{,}434 in 2025~\cite{sia2025,sia2026}, driven largely by
LEO broadband megaconstellations~\cite{brycetech2026}. Yet after the launch rate passes $\lambda_+$, the
fragment count remains near its low-branch value for $\Delta t$ time units,
during which the shell \emph{appears} healthy. To leading order in $\varepsilon$,
the number of additional satellites admitted in this interval is
\[
  N_{\mathrm{excess}} \;\approx\; \lambda_+\,\Delta t
  \;=\; \lambda_+\,\Omega_0\,(\beta_\lambda\,\beta_2)^{-1/3}\,\varepsilon^{-1/3},
\]
while the launch rate itself overshoots the threshold by $\Delta\lambda \sim
\varepsilon^{2/3}$ before the fragment count enters the high-debris branch.
The lag is \emph{longer} for \emph{slower} ramps, since
$\Delta t \propto \varepsilon^{-1/3}$, a gradual rise in launch rate produces a
longer window of deceptive health, even though the launch-rate overshoot
$\Delta\lambda \propto \varepsilon^{2/3}$ shrinks, since there are $N_{\mathrm{excess}} \propto
\varepsilon^{-1/3}$ additional launches.
An operator waiting for the fragment count to
rise is reacting only after the sustainable threshold has been crossed; the
launch rate, not the fragment count, therefore, is the leading indicator.

\begin{figure}[ht]
    \centering
    \includegraphics[width=0.8\linewidth]{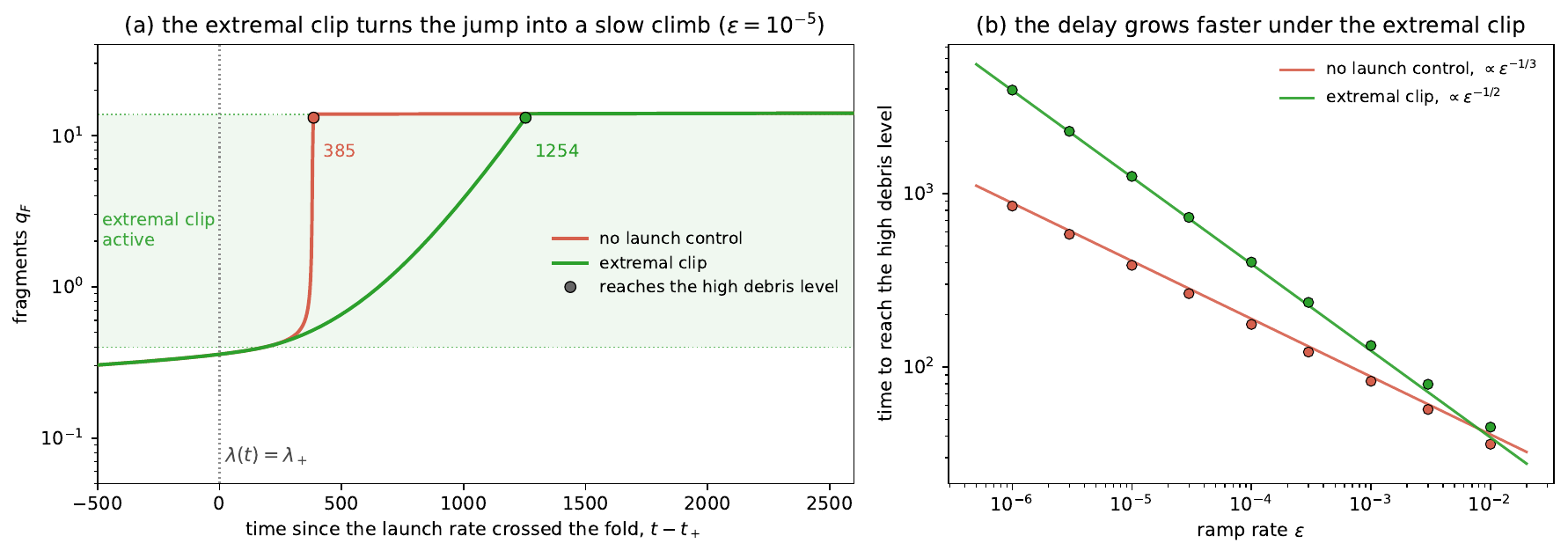}
    \caption{Slow passage under the extremal clip, in the setting of Figure~\ref{fig:slow-passage}. (a)~Without launch control (red), the fragment count jumps to the high debris branch shortly after the launch rate crosses the fold; under the extremal clip (green), it climbs gradually through the range where the clip is active and reaches the high debris level later. (b)~The time to reach the high debris level grows like $\varepsilon^{-1/3}$ without launch control and like $\varepsilon^{-1/2}$ under the extremal clip (Theorems~\ref{thm:slow-passage} and~\ref{prop:clip-slow-passage}).}
    \label{fig:slow-passage-extremal-clip}
\end{figure}

Motivated by this, the second part of our results develops launch control
policies that eliminate the bistability by removing the high-debris branch inside
the bistable window.
We call a policy a \emph{throttle} if it reduces the launch rate as a function of
the state of the shell (specifically, the states of the intact, derelict, and fragment queues).
A natural first idea is to throttle on the number in the intact queue, launching fewer new satellites as the
fleet grows.
We show that this cannot remove the bistability and in fact works against its
purpose. Observe that the intact count falls as debris destroys satellites, and the policy
responds by admitting more launches, thereby feeding the process it is meant to contain.
By contrast, a throttle that depends only on the fragment count can remove the
bistability.
Among such throttles, we characterize the least restrictive one, which we call the \emph{extremal clip}.
The extremal clip only throttles launch rates when fragment counts are above the maximal count on the low-debris branch and below the minimal count on the high-debris branch.
Under this policy, the hysteresis is removed, and the slow passage exhibits a smoother transition from the low-debris branch to the high-debris branch but it still exhibits a slow passage phenomenon, see Figure~\ref{fig:slow-passage-extremal-clip} (a). In fact, under this policy, the slow passage delay grows as $\Delta t \sim \varepsilon^{-1/2}$ (Theorem~\ref{prop:clip-slow-passage}; Figure~\ref{fig:slow-passage-extremal-clip} (b)).

This changes the nature of the lag, but not its existence.
After the launch rate passes $\lambda_+$, the fragment count no longer collapses onto the high-debris branch.
Instead, it begins to rise after a time delay of order $\varepsilon^{-1/3}$ after the crossing, reaching the high-debris level only after a delay of
\[
  \Delta t \;\sim\; \Big(\frac{2\lambda_+\,\mathcal{J}_{\mathrm{w}}}{\mu_F}\Big)^{1/2}\varepsilon^{-1/2}
  \qquad (\text{as}~\varepsilon \to 0),
\]
where $\mathcal{J}_{\mathrm{w}}$ is a model constant which measures the width of the interval of fragment counts where the extremal clip is active.
Correspondingly, $N_{\mathrm{excess}} \approx \lambda_+\,\Delta t \propto \varepsilon^{-1/2}$ additional satellites are admitted before the fragment count reaches the high-debris level, and the launch rate overshoots $\lambda_+$ by $\Delta\lambda \propto \varepsilon^{1/2}$.
Since $\varepsilon^{-1/2} \gg \varepsilon^{-1/3}$ for slow ramps, the fragment count takes even longer to reach the high-debris level than without the throttle: for $\varepsilon = 10^{-5}$ in Figure~\ref{fig:slow-passage-extremal-clip}(a), about $1{,}250$ rather than $385$ time units.
However, two features which made the uncontrolled fragment count dangerous are weakened.
First, the fragment count starts to rise early and increases steadily, so a sustained upward trend is visible long before the high-debris level is reached, instead of a sudden collapse after a period of apparent health.
Second, the transition is reversible, and there is no hysteresis. Therefore, lowering the launch rate back below $\lambda_+$ returns the shell to the low-debris branch, and the satellites admitted during the lag do not trap it on the high-debris branch.

\begin{figure}[ht]
    \centering
    \includegraphics[width=\linewidth]{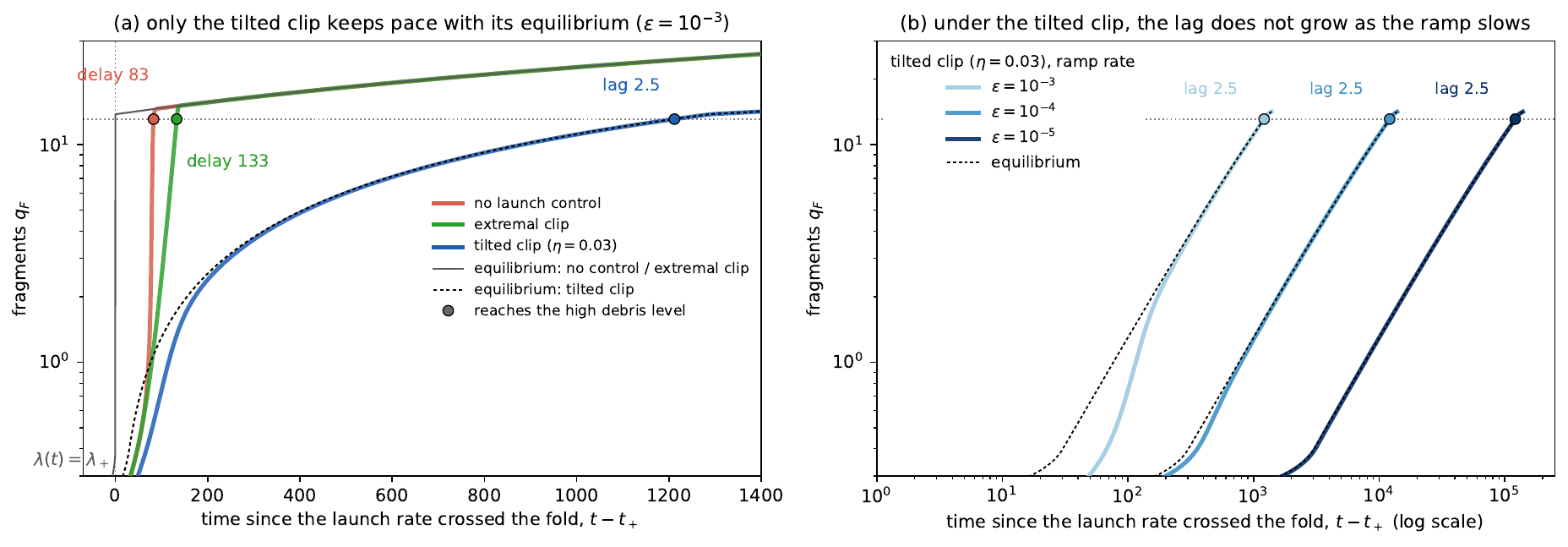}
    \caption{Removing the delay with the tilted clip, in the setting of Figure~\ref{fig:slow-passage}. Thick curves are trajectories, and thin curves are the equilibria they follow as the launch rate increases. (a)~For one ramp rate, the trajectories without launch control and under the extremal clip reach the high debris level well after their equilibria do, whereas the trajectory under the tilted clip stays close to its equilibrium. (b)~Under the tilted clip, the lag behind the equilibrium stays the same as the ramp slows.}
    \label{fig:slow-passage-tilt}
\end{figure}

In order to eliminate the delay $\Delta t$, the throttle must be further changed. We show that any fixed exponential tilt of the extremal clip achieves this, and the fragment count becomes a timely rather than a lagging indicator, giving operators feedback as the threshold is crossed rather than after the fact.
Under a tilted clip, the equilibrium fragment count rises gradually with the launch rate instead of jumping at $\lambda_+$, and the fragment count follows it with a lag which does not grow as the ramp slows.
Specifically, the lag observed is of constant order (see Theorem~\ref{prop:clip-slow-passage}).
At each fragment count $q_F$, the tilted clip admits a fraction $e^{-\eta q_F}$ of the launches the extremal clip admits, so that the logarithm of the ratio of the two admission fractions is proportional to $q_F$, with constant of proportionality $\eta > 0$. Letting $\eta \to 0$ recovers the extremal clip (see Theorem~\ref{prop:heaviest}\,(iv)).
Figure~\ref{fig:slow-passage-tilt} illustrates this via comparisons to the un-controlled launches and the extremal clip (panel (a)) and the removal of the lag as the ramp rate tends to zero (panel (b)).
In the example, the fragment count reaches the high-debris level about $2.5$ time units after its equilibrium does, for every ramp rate from $\varepsilon = 10^{-3}$ to $\varepsilon = 10^{-5}$; at $\varepsilon = 10^{-3}$, the corresponding delays are about $83$ time units without launch control and $133$ under the extremal clip.
The tilted clip also holds the equilibrium fragment count lower at every launch rate above $\lambda_\crit$, at the price of blocking more launches than the extremal clip.

Although we develop these results for orbital debris, the arguments behind the launch-control results rely on a singular structural feature of the fluid limit, namely that every non-trivial equilibrium $(q^*_I, q^*_D, q^*_F)$ is determined by its fragment count $q_F^*$ through a scalar equation of the form
\[
    \lambda\,\rho(q_F^*) = \cM(q_F^*).
\]
Here $\rho(q) \in (0,1]$ is the fraction of offered launches that a throttle admits at fragment count $q$, and $\cM(q)$ is the launch rate at which the uncontrolled system has an equilibrium with fragment count $q$; the equilibrium intact and derelict counts are then explicit functions of $q_F^*$ and $\lambda\rho(q_F^*)$ (Lemma~\ref{lem:G-bifurcation}).
Given this reduction, the controlled system has no fold if and only if $\cM / \rho$ is non-decreasing, and the least restrictive throttle with this property is the extremal clip
\[
    \rho^\star(q) = \frac{\bar H(q)}{H(q)}, \qquad \bar H(q) := \inf_{0 \le s \le q} H(s),
\]
where $H := \mu_F/(c\,\kfrags\,\cM)$ is the function defined in~\eqref{eq:H}, so that $c\,\kfrags\,\lambda H(q)$ is the rate at which a single fragment generates new fragments when no launches are blocked and the intact and derelict counts are at their equilibrium values given fragment count $q$.
Under the extremal clip, $\rho^\star H = \bar H$, and the equilibrium equation becomes $c\,\kfrags\,\lambda\,\bar H(q_F^*) = \mu_F$, which is the uncontrolled equation with $H$ replaced by its running infimum, or equivalently with $\cM$ replaced by its running maximum (Theorem~\ref{prop:heaviest}).
Moreover, the admitted launch rate and the net useful throughput at an equilibrium depend on the policy only through $q_F^*$ (Proposition~\ref{prop:throughput}), so that a comparison of policies at equilibrium reduces to a comparison of the fragment counts which they sustain.
Likewise, the exponent of the delay under a slowly increasing launch rate is determined by the shape of $\cM/\rho$ near the critical fragment count, which has a quadratic turning point without launch control, is flat under the extremal clip, and is strictly increasing under the tilted clip (Remark~\ref{rem:hierarchy}).
Bistability that admission control removes also arises, through different mechanisms, in loss networks with dynamic alternative routing~\cite{gibbens1990bistability,oleskertaylor2022metastability} and in service systems with slowdown~\cite{dong2015service}.
In any model whose equilibria reduce to a scalar equation of the form above (with $\cM \to \infty$ as $q \to \infty)$, with the throttle acting on the arrival rate, the same argument identifies the least restrictive throttle which removes the folds as the one that replaces the equilibrium relation by its running maximum.

Finally, we address the question of what the policies deliver, as measured by the \emph{net useful throughput}, defined as the rate $P\mu_I q_I^*$ at which satellites complete their missions and are disposed of successfully.
At an equilibrium, the net useful throughput depends on the policy only through the equilibrium fragment count (Proposition~\ref{prop:throughput}), because a fraction $\mu_I/(\mu_I + \alpha(q^*_F)cq^*_F)$ of the admitted satellites complete their missions and the remainder are destroyed in collisions.
Admitting more launches can reduce the net useful throughput, as highlighted by the following example.
With the system parameters as in Figure~\ref{fig:bistability} ($\Qsat = 10$), the low-debris equilibrium near the tipping point $\lambda_+$ admits $2.51$ launches per unit time, of which $99\%$ complete their missions, whereas at $\lambda = 4$, only $17\%$ reach mission completion, resulting in about a quarter of the useful throughput.

As the equilibrium debris grows without bound, the net useful throughput of every policy approaches the same limit $P\mu_I\mu_F/(c\,\kfrags)$.
In Theorem~\ref{thm:bang-bang}, we identify an equilibrium fragment level $q^\circ$ which allows for the maximum useful throughput.
We show that a \emph{bang-bang} admission rule, which admits all launches when the fragment count is below $q^\circ$ and none when it is above, attains the maximal useful throughput at every offered launch rate $\lambda$.
The bang-bang rule achieves the maximal ratio of useful throughput to offered launch rate (Corollary~\ref{cor:useful-fraction}).
The bang-bang rule also removes the bistability, with the fraction of time during which launches are admitted in the finite-$n$ system converging to the admitted fraction in the fluid limit (Theorem~\ref{thm:bang-bang-fluid}).

Bistability has precedent in stochastic networks studied in the
mean-field or fluid limit~\cite{antunes2008stochastic,gibbens1990bistability,martirosyan2020equilibrium,oleskertaylor2022metastability, dong2015service}, although it arises there by mechanisms different from
ours; see our literature survey below on bistability, metastability and fluid limits in stochastic networks.
In LEO, the system metrics differ from those more traditionally studied in the stochastic networks literature.
The article~\cite{dong2015service} is perhaps the closest to our work; it studies a tradeoff between quality-driven and efficiency-driven staffing in a queue which has a bistable regime.
However, metrics such as \emph{useful throughput}, the rate of successful mission completions, are more important to LEO operations.
Useful throughput is also studied in congestion control for communication networks.
In that setting, useful throughput is an important metric due to stale and fresh information being transmitted simultaneously on a network link.
In contrast, useful throughput is the rate at which intacts depart the system through PMD.

\paragraph{Operational Implications}
The operational implications of our results, as we see them, are as follows:

\begin{itemize}
    \item Removing bistability should be a paramount concern for a manager. Recovering from a high-debris branch requires a significant decrease in the input rate (i.e., the launch rate) due to the effects of hysteresis. Bistability concerns are further compounded by the slow passage effects. We show that a fragment-indexed input policy (i.e., the extremal clip) can remove the bistability, at the cost of a longer slow passage. The latter can be reduced by an exponential tilt of the input policy.

    \item Recall that $\Qsat$ measures collision avoidance capabilities, and $\alpha_0$ is the baseline collision avoidance failure probability.
    There are two inherent trade-offs that need to be considered for system design with saturating routing probabilities.
    First, monostability can be guaranteed if $\Qsat > \overline{Q} = \Theta(\alpha_0^{-2})$ as $\alpha_0 \to 0$. Second, lower $\alpha_0$ can raise the tipping point, which is good. But it also widens the hysteresis and increases the launch rate cut needed for recovery. In general, when a safeguard of finite capacity protects against a hazard generated by the system's {\it own} failures, reliability gains at a fixed capacity raise the tipping point by a bounded amount and widen the hysteresis, whereas removing the bistability at every load requires a capacity which grows like the inverse square of the baseline failure probability.

    \item A larger launch rate can significantly reduce the useful throughput, since a larger fraction of intacts may undergo catastrophic collisions instead of completing their missions. However, under a unimodality condition, which admits every launch while the fragment count is below a threshold $q^\circ$ and none above the threshold, maximizes the equilibrium rate of successful mission completions at every launch rate. The threshold $q^\circ$ does not depend on the launch rate, so the rule needs the debris level but no forecast of demand. Moreover, $q^\circ$ lies below the fragment count at which the system tips, so maximizing the useful throughput does not conflict with removing the bistability. More generally, if the load admitted to a system generates a hazard which destroys work in progress, and the equilibrium rate of completed work depends on the admission policy only through the level of the hazard, rising and then falling as the hazard grows, then the best policy keeps the hazard at or below the level where that rate peaks, whatever the demand.

\end{itemize}

\input{related_work}

\subsection{Organization}
The rest of this paper is organized as follows.
In Section~\ref{sec:model}, we introduce the model and its fluid limit.
In Section~\ref{sec:bistability}, we characterize the bifurcation structure of the fluid limit, including the Hopf bifurcation (Section~\ref{sec:hopf}) and the slow passage through the fold (Section~\ref{sec:slow-passage-result}).
In Section~\ref{sec:admission-control}, we study launch control.
Concluding remarks are in Section~\ref{sec:conclusions}.
Proofs are in the appendices.

%% file: related_work.tex
\subsection{Related Work}
\label{sec:related-work}

Our literature survey is organized according to several lines of work which surround this paper.

\paragraph{Debris evolution and orbital capacity.}
The idea that collisions between orbiting objects can sustain a growing
population of fragments originates with Kessler and
Cour-Palais~\cite{kessler1978}. Its simplest quantitative descendants are
particle-in-a-box and source--sink models, which track the mean numbers of
intacts, derelicts, and fragments in an altitude shell with ordinary
differential equations~\cite{talent1992analytic,bradley2009space}. Bradley and
Wein~\cite{bradley2009space} calibrate such a model to the 900--1000~km shell and
define a \emph{sustainable risk} as the maximum, over future time, of a
spacecraft's probability of destruction by collision.
They also quantify how the sustainable risk depends on
compliance with post-mission disposal guidelines. Lafleur
\cite{lafleur2011extension} uses a similar ODE model and estimates that the equilibrium
capacity for intacts is only 25\% of the 2011 intact population (present-day for that paper), and that
direct removal of fragments could extend the useful life of LEO by decades. The
MIT Orbital Capacity Assessment Tool (MOCAT) line of work builds multi-shell,
multi-species source-sink ODE models.
For the MOCAT tool, D'Ambrosio~\textit{et al.} define orbital capacity through the stable equilibria of the dynamical system~\cite{dambrosio2022capacity,dambrosio2023novel}.
Further work by Pasiecznik~\textit{et al.} studies how launch rates affect the stability of those
equilibria~\cite{pasiecznik2022dynamical}.
A complementary approach computes the critical number of intact objects at which a shell becomes unstable or enters a runaway regime~\cite{lewis2025critical}; applied to sixteen planned and operational constellations, it finds that six exceed the runaway threshold~\cite{lewis2026critical}. Three-dimensional evolutionary simulators, such as LEGEND~\cite{liou2004legend}, provide higher fidelity at far greater computational cost.

These models share two features. First, they are deterministic, so they cannot represent the randomness of collision times, breakup sizes, and launch schedules.
Second, they treat collision avoidance, when they include it at all, as a fixed input.
In long-term simulations supporting ESA's Zero Debris policy, for example, constellation satellites are assumed to avoid collisions with a success rate of 100\%~\cite{letizia2023zero}. Parametric studies vary avoidance and disposal rates as exogenous parameters~\cite{schaus2021revisiting}.
The critical-size analysis~\cite{lewis2026critical} treats the absence of propulsion as a fixed
design property of a constellation.
Many collision-rate simulations assume perfect avoidance for all maneuverable payloads~\cite{thiele2025orbital}.

Stochastic formulations are more recent. Behera~\textit{et
al.}~\cite{behera2026stochastic} extend the Bradley--Wein model to a density-dependent Markov chain and derive fluid limits, functional central limit theorems, and a Feller diffusion on distinct time-scales, showing that a collisional cascade can be triggered by fluctuations
alone.
Choi~\textit{et al.}~\cite{choi2026narrow} formulate the multi-shell,
multi-species environment as a Markov jump process, recover the source-sink ODE
as a large-volume limit, and show that in narrow shells the mean debris
population can substantially exceed the ODE prediction.

\paragraph{Collision avoidance in large constellations.}
Collision avoidance is increasingly a workload and scalability problem.
Operators receive conjunction data messages (CDMs) from the U.S.\ Space Force.
The CDM volume has reached roughly 600{,}000 per day, about three times the 2020 rate, and Starlink is involved in just over half of the messages~\cite{amos2023lessons,semidecentralized2026}.
Earlier work argued that manual or semi-automated conjunction assessment and avoidance planning is already a sizable human workload and that automation is necessary for space traffic management to scale~\cite{aas2019automating}. Starlink's response has been to fully automate the avoidance system in which satellites assume maneuver responsibility by default~\cite{starlink2026conops}, supported by SpaceX's tracking and conjunction screening service Stargaze~\cite{spacex2026stargaze}.
Empirically, CDM arrivals are consistent with a Poisson process~\cite{cdmpoisson2023}, which motivates the Poisson conjunction model in the first layer of our model.
Avoidance has costs of its own: Chen~\textit{et al.}~\cite{chen2024instability} show, using Starlink's situational awareness data, that pairwise maneuvers can cascade inside a megaconstellation, creating a tension between intact lifetime and system capacity.
Sturza~\textit{et al.}~\cite{sturza2024resilience}
show that resilience of a constellation may not increase monotonically with its size.
The resources that bound avoidance capacity are also studied: satellites operate under size, weight, power, and cost budgets that limit onboard compute~\cite{onorbitai2026}, and operators receive updated tracking information only during intermittent ground
contacts~\cite{semidecentralized2026}.
Finally, the CRASH Clock of Thiele~\textit{et al.}~\cite{thiele2025orbital} measures how quickly a collision would occur if avoidance maneuvers stopped, finding a decline from 164 days in 2018 to 5.5 days in 2025.
This quantifies how much current and future safety depends on collision avoidance technologies.
To our knowledge, none of these works embeds a state-dependent, saturating avoidance capability in a population-level model of debris dynamics, so that the fraction of conjunctions an operator can avoid falls as the fragment count rises and feeds back into the debris dynamics, which is the coupling at the core of our
model.

\paragraph{Economics and governance of orbit use.}
A parallel literature treats orbit as an open-access commons. Adilov~\textit{et
al.}~\cite{adilov2015earth,adilov2018economic} analyze the externality created by
orbital debris and its consequences for the optimal use of orbit. Rouillon
\cite{rouillon2020physico} studies a long-run equilibrium induced by a constant
launch rate and shows that, when collision risk is increasing and convex in the
launch rate, the expected number of functioning satellites is a reversed-U
function of the launch rate, whose maximum defines the carrying capacity of the
orbit.
Guyot and Rouillon~\cite{guyot2023sustainable} extend the analysis to
removal services and show that a combination of taxes and a market for removal
effort can provide the right incentives. Rao~\textit{et al.}~\cite{rao2020orbital} consider debris in an economic model of launches and find that an
orbital-use fee could substantially raise the value of the space industry, while
debris removal can worsen congestion damages by increasing launch incentives.
Rao and Rondina~\cite{raorondina2022economics} show that the Kessler Syndrome can occur both with open-access
and when external collision costs are internalized.
Bernhard~\textit{et
al.}~\cite{bernhard2023large} study large constellations using dynamic games, and Guyot~\textit{et
al.}~\cite{guyot2023oligopoly} study duopoly competition between constellation
operators, where imperfect competition reduces welfare.
Boley and
Byers~\cite{boley2021satellite} warn that rapid constellation deployment risks
tragedies of the commons in orbit, in the atmosphere, and on the ground. These
models endogenize launch behavior but typically represent the physical
environment by stylized collision risk functions.
The multiplicity of
equilibria that appears in~\cite{raorondina2022economics} comes from economic
incentives, whereas in our model it comes from the physical and operational
feedback of saturating avoidance, with launches taken as exogenous. Our launch control results are therefore complementary to this line of work.
We ask how much launch blocking is required to eliminate bistability, and not who should bear its cost. %

\paragraph{Remediation and control of the debris environment.}
Remediation research studies how to remove or prevent debris. For active debris
removal (ADR), the focus is on target selection and mission planning, including
risk-based ranking indices~\cite{medhin2025sustainability} and
reward-maximizing mission planning~\cite{yang2019maximal}. Somma~\textit{et
al.}~\cite{somma2016adaptive,somma2017statistical} embed a feedback controller in a
source--sink model, with the number of removals proportional to the orbital
population.
Rusconi~\textit{et al.}~\cite{rusconi2024feedback} couple a feedback
controller to a density-based model to reach a target scenario. Liou~\textit{et
al.}~\cite{liou2010controlling} study how removal rates affect future LEO
populations. These controls act on debris after it has been
created. Our controls instead act on the \emph{inflow}: we study launch
admission as a function of the state of the shell.
To our knowledge, our finding that
a throttle on the intact count cannot remove bistability, whereas a throttle on
the fragment count can, has no direct counterpart in the ADR or feedback-control
literature.

\paragraph{Resilience of constellations as infrastructure.}
Simulation studies on Starlink-like constellations find
that inter-satellite networks inherit resilience to node failures from their
dynamic topology, and that rerouting is critical for rapid
recovery~\cite{nodefailures2025}.
Analyses of ground hubs find that connectivity and routing
performance can be maintained with as few as twenty functioning ground
stations~\cite{smallworld2025}. This body of work, together with defense-policy
analyses that favor proliferated architectures as a way to deny an adversary the
benefit of an attack~\cite{npec2020war,fiore2023deterrence}, underlies the view
that destroying a constellation's service requires disabling a large fraction of
its satellites; see also the government programs described
in~\cite{gao2026missile}. It typically treats satellite loss as exogenous,
whether random or adversarial; in our model, satellite loss is endogenous,
driven by debris that the constellation's own launches help create.

\paragraph{Bistability, metastability, and fluid limits in stochastic networks.}
The fluid limit for our queueing network arises from a functional law of
large numbers for density-dependent Markov
chains~\cite{kurtz1970solutions,kurtz1981approximation,ethierkurtz1986markov}.
We use ideas similar to those developed by
Hunt and
Kurtz~\cite{huntkurtz1994large} for loss networks (see Kelly~\cite{kelly1991loss} for a review of loss networks). Bistability of this kind has precedent in loss
networks, although it arises there by different mechanisms. In loss networks
with dynamic alternative routing, Gibbens \textit{et al.}~\cite{gibbens1990bistability} show that random alternative routing can
produce instability and hysteresis, and that trunk reservation mitigates it.
Rigorous results on the resulting metastability
include~\cite{tibi2010metastability,martirosyan2020equilibrium, oleskertaylor2022metastability}.
In Antunes \textit{et al.}~\cite{antunes2006analysis,antunes2008stochastic}, coexisting equilibria
arise from blocking under heterogeneous capacity with fixed routing. In our
model a blocked launch is lost rather than rerouted, and the state dependence enters through the collision-avoidance rate; the
fold therefore has a different origin.
Bistability can also arise in systems with mobile servers~\cite{baccelli2018metastability}.

The most closely related work to ours are~\cite{dong2015service, selen2016snowball}.
Dong \textit{et al.}~\cite{dong2015service} study a scalar system where feedback-induced bistability from
state-dependent service rates in a multi-server queue with server slowdown.
They remove the bistability using an admission
control.
Selen \textit{et al.}~\cite{selen2016snowball} identify bistability in a multi-server model where customers who enter service immediately have a faster service rate than those who wait before service, which requires a two-dimensional state space. They note that exactly identifying the multiple stable equilibria in more than one dimension is a difficult problem.

For the stochastic behavior of systems with multiple
stable states, the relevant general theory includes large deviations and
exit-time asymptotics~\cite{freidlinwentzell2012random,shwartzweiss1995large} and
potential-theoretic metastability~\cite{bovier2015metastability}.

\paragraph{Slow passage and delayed bifurcation.}
The slow passage through a bifurcation phenomenon was analyzed for algebraic
bifurcation problems by Haberman~\cite{haberman1979slowly}, for steady
bifurcations by Mandel and Erneux~\cite{mandel1987slow}, and for Hopf
bifurcations by Baer, Erneux, and Rinzel~\cite{baer1989slow}.
We refer the reader to
Neishtadt~\cite{neishtadt1987persistence} for the analytic theory of stability
loss, Kuehn~\cite{kuehn2015multiple} for a modern treatment of slow-fast systems, and
Krupa and Szmolyan~\cite{krupa2001} for the rigorous geometric analysis
of passage through fold points. Berglund and Gentz~\cite{berglund2006noise}
study noisy slow--fast systems and treat the saddle-node
bifurcation and bifurcation delay explicitly. The consequence for
monitoring is that observables lag the transition, which connects to the
literature on early-warning signals for critical transitions
\cite{scheffer2009early,scheffer2001catastrophic}. %
The textbook~\cite{kuznetsov2023elements} is a reference for the fold and transcritical bifurcations used in
our classification results.

%% file: avoidance_layer.tex
\subsection{Per-Satellite Conjunction Avoidance Layer}\label{sec:avoidance-layer}

For each intact, we model collision avoidance as a processor-sharing system in
which the customers are the debris fragments in the shell. The intact has an
avoidance capacity $\mu_s$, the rate at which it can screen fragments and plan
maneuvers. This capacity is shared equally among all fragments, so that when
there are $q$ fragments, each is processed at rate $\mu_s/q$.

Conjunctions between the intact and a given fragment occur at the times of a
Poisson process with rate $c$. Each conjunction carries a \emph{deadline}, the
time remaining until closest approach, which we take to be exponential with
rate $\nu$. The conjunction is avoided if the fragment is processed before
the deadline expires; otherwise the intact suffers a catastrophic collision.
Collision avoidance is thus a race between processing, at rate $\mu_s/q$, and
the deadline, at rate $\nu$. Because both clocks are exponential, a
conjunction that occurs when $q$ fragments are present fails with probability
\[
    \frac{\nu}{\nu + \mu_s/q} = \frac{q}{q + \Qsat},
    \qquad
    \Qsat := \frac{\mu_s}{\nu}.
\]
Deadlines last hours to days, during which the fragment count is essentially
constant, so we evaluate this probability at the current value of $q$.
Avoidance systems are imperfect even when the race is won, so a processed
conjunction still fails with a baseline probability $\alpha_0 \in (0,1)$. The
probability that a conjunction results in a catastrophic collision is therefore
the \emph{collision avoidance failure function}
\[
    \alpha(q) := \alpha_0 + (1 - \alpha_0)\frac{q}{q + \Qsat},
\]
and, since conjunctions with the $q$ fragments occur at total rate $cq$, an
intact suffers catastrophic collisions at rate $\alpha(q)\,c\,q$.

Two features of this construction matter. First, the conjunction rate $c$ and
the deadline rate $\nu$ are distinct: $c$ governs how often an intact meets a
given fragment, and $\nu$ governs how much warning it has when it does.
Second, the capacity is shared among all $q$ fragments, not only among those
currently in a conjunction. Screening must be done in advance, since a
conjunction is detected only by screening, so the capacity is spent maintaining
predictions for every fragment. The parameter $\Qsat$ measures this capacity
relative to the warning time: when $q = \Qsat$, half of all conjunctions are
lost to the deadline, and as $q$ grows beyond $\Qsat$ the capacity is spread
ever more thinly and $\alpha(q) \to 1$.

We assume that collision avoidance dynamics for different intacts are independent.

\paragraph{Custody interpretation.}
The same failure function has an equivalent description in terms of tracking
custody, which clarifies what $\Qsat$ represents. Say that a fragment is
\emph{in custody} if its orbit is known well enough to screen against it.
Processing, at rate $\mu_s/q$, restores custody, and custody is lost at rate
$\nu$ as the orbit prediction degrades. Each fragment then alternates between
being in and out of custody, and in stationarity it is out of custody with
probability
\[
    \frac{\nu}{\nu + \mu_s/q} = \frac{q}{q + \Qsat}.
\]
A conjunction with a fragment out of custody cannot be avoided, so a conjunction
again fails with probability $\alpha(q)$. The two descriptions are time
reversals of each other: in the processor-sharing description the fragment must
be processed within an exponential window \emph{after} the conjunction begins,
and in the custody description it must have been processed within an
exponential window \emph{before}. Under this interpretation, $\Qsat$ is the
\emph{custody capacity}, the number of fragments that the avoidance system can
hold in custody at once relative to how quickly custody decays. This interpretation applies
equally whether screening is done on board each satellite or by a shared
catalog such as a ground- or space-based surveillance network; in the latter
case, $\mu_s$ is the throughput of the catalog rather than of a single
satellite. In either reading, collision avoidance saturates because a fixed
capacity is spread over a growing fragment population, not because any single
conjunction is hard to avoid.

%% file: admission-control.tex
\section{Admission Control: Proofs}\label{sec:admission-proofs}

This section collects the proofs of the admission-control results of Section~\ref{sec:admission-control}.
The occupancy/debris feedback dichotomy is Propositions~\ref{prop:occupancy}--\ref{prop:debris}, Corollary~\ref{cor:exp}, and Lemma~\ref{lem:clip-linearization}; the largest fold-free throttle is Theorem~\ref{prop:heaviest}; the throttled process and its fluid limit are Lemma~\ref{lem:rho-lip} and Theorems~\ref{thm:clip-fluid} and~\ref{thm:tilted-clip}; and the slow passage under the extremal and tilted clips is Theorem~\ref{prop:clip-slow-passage}.
Throughout, $\sigma_I(z) = \mu_I + \alpha(z)cz$, $\sigma_D(z) = \mu_D + cz$, and we use the identity $\sigma_I H = \alpha + (1-P)\mu_I/\sigma_D$.

The throttled process uses the following Lipschitz property of $\rho^\star = \bar H/H$.
\begin{lemma}[The clip is Lipschitz]\label{lem:rho-lip}
    $\rho^\star$ is Lipschitz on $[0,\infty)$, takes values in $(0,1]$, and equals $1$ outside the compact window ($[q_{\min}, q^\ddagger]$ in the two-fold regime, $[0, q^\ddagger]$ in the one-fold regime), with Lipschitz constant $L_\rho = \sup_{[q_{\min}, q^\ddagger]} |H'|\,/\,H(q_{\min})$.
\end{lemma}
\begin{proof}
    $H$ is $C^1$ on $[0,\infty)$ with bounded derivative (it is a rational function of $q$ whose denominator does not vanish, and $H' = O(q^{-2})$), so $H$ is globally Lipschitz, positive and bounded.
    On the window $\bar H \equiv H(q_{\min})$, so $\rho^\star = H(q_{\min})/H$ with $H \ge H(q_{\min}) > 0$, giving $|\rho^{\star\prime}| = H(q_{\min})|H'|/H^2 \le |H'|/H(q_{\min})$; off the window $\bar H = H$, so $\rho^\star \equiv 1$.
    As $\bar H \le H$ we have $\rho^\star \le 1$, and $\bar H > 0$ on every $[0,x]$ gives $\rho^\star > 0$.
\end{proof}

\begin{proof}[Proof of Theorem~\ref{thm:clip-fluid}]
    The throttle changes only the launch intensity, which becomes $$a(t) = \int_0^t n\lambda\,\rho(\bar Q^{n,\rho}_F(s))\,\ind_{Q_I(s) < n}\,ds \leq n\lambda t.$$
    Lemma~\ref{lem:apriori} and Theorem~\ref{thm:strong-approx} therefore hold for $Q^{n,\rho}$ with the same proofs, since they use only the bound $a(t) \leq n\lambda t$.
    Dividing by $n$ as in the proof of Proposition~\ref{prop:reflected-fluid}, $a(t)/n = \int_0^t \lambda\rho(\bar Q^{n,\rho}_F)\,ds - \ell^n(t)$ with $\ell^n(t) := \int_0^t \lambda\rho(\bar Q^{n,\rho}_F)\ind_{\bar Q^{n,\rho}_I = 1}\,ds$, which is continuous, non-decreasing, and increases only when $\bar Q^{n,\rho}_I = 1$.
    Hence~\eqref{eq:prelimit} holds with $b_I$ replaced by $b^\rho_I(q) := \lambda\rho(q_F) - \mu_I q_I - \alpha(q_F)cq_Iq_F$.
    Since $\rho$ is bounded and Lipschitz, $b^\rho$ is locally Lipschitz, the bounds of the proof of Proposition~\ref{prop:reflected-fluid} hold for the reflected controlled system (they use only $b^\rho_I \leq \lambda - \mu_I q_I$), and the remainder of that proof applies here.
    If $\lambda \leq \mu_I$, then $\widetilde q$, the solution of~\eqref{eq:cFL-1}--\eqref{eq:cFL-3}, satisfies $\frac{d}{dt}(1 - \widetilde q_I) \geq \mu_I\widetilde q_I - \lambda \geq -\mu_I(1 - \widetilde q_I)$, so it stays in $\mathcal S$, and by uniqueness it is the limit.

    For $\lambda \leq \mu_I$, \eqref{eq:cFL-1} and~\eqref{eq:cFL-2} at a fixed point give $q_I^* = \lambda\rho(z)/\sigma_I(z)$ and $q_D^* = (1-P)\mu_I q_I^*/\sigma_D(z)$, with $z = q^*_F$, and for $z > 0$, \eqref{eq:cFL-3} divided by $z$ gives $\kfrags\,c\,q^*_I\sigma_I(z)H(z) = \kfrags\,c\,\lambda\rho(z)H(z) = \mu_F$.
\end{proof}

\begin{proof}[Proof of Proposition~\ref{prop:occupancy}]
    \emph{(i).} Fix $z \geq 0$. The map $q \mapsto \bar\lambda(q) - \sigma_I(z)q$ is continuous and strictly decreasing on $[0,1]$, positive at $q = 0$, and non-positive at $q = 1$ because $\bar\lambda \leq \mu_I \leq \sigma_I(z)$; so it has a unique root $q_I^\circ(z) \in (0,1]$.
    Since $\bar\lambda \leq \mu_I$, the fluid limit is~\eqref{eq:FL-1}--\eqref{eq:FL-3} with $\lambda$ replaced by $\bar\lambda(q_I)$.
    At a fixed point with $q^*_F = z$, the intact equation gives $q^*_I = q_I^\circ(z)$, the derelict equation gives $q^*_D = (1-P)\mu_I q^*_I/\sigma_D(z)$, and for $z > 0$ the fragment equation divided by $z$ gives
    \[
        \kfrags\, c\, q^*_I\Big(\alpha(z) + \frac{(1-P)\mu_I}{\sigma_D(z)}\Big) = \kfrags\, c\, q^*_I\,\sigma_I(z)H(z) = \kfrags\, c\, J(z)H(z) = \mu_F.
    \]
    If $z_1 < z_2$, then $\sigma_I(z_1) < \sigma_I(z_2)$, so at $q = q_I^\circ(z_1)$ we have $\bar\lambda(q) - \sigma_I(z_2)q < 0$, and therefore $q_I^\circ(z_2) < q_I^\circ(z_1)$; as $\bar\lambda$ is non-increasing, $J(z_2) = \bar\lambda(q_I^\circ(z_2)) \geq \bar\lambda(q_I^\circ(z_1)) = J(z_1)$.

    \emph{(ii).} $J$ is positive and non-decreasing, so $JH$ is increasing wherever $H$ is.
    For the instability, let $q^*$ be a fixed point with $q^*_F = z$ in such an interval, and $\bar\lambda$ continuously differentiable near $q^*_I$.
    The Jacobian has the same form as in Lemma~\ref{lem:linearization}, except that its $(1,1)$ entry is $\bar\lambda'(q^*_I) - \sigma_I$, so its upper-left $2\times2$ block has determinant $(\sigma_I - \bar\lambda'(q^*_I))\sigma_D > 0$.
    As in the proof of Lemma~\ref{lem:linearization}, the Schur complement of this block is the derivative of the fragment equation with $q_I, q_D$ determined by $q_F$, namely of $\zeta \mapsto \zeta(c\,\kfrags\,J(\zeta)H(\zeta) - \mu_F)$, which at the fixed point equals $z\,c\,\kfrags\,(JH)'(z) \geq z\,c\,\kfrags\,J(z)H'(z) > 0$ (by the implicit function theorem $J$ is differentiable at $z$, with $J' = \sigma_I'\,q_I^\circ\,\bar\lambda'/(\bar\lambda' - \sigma_I) \geq 0$).
    Hence $\det \cL > 0$, the constant coefficient of the characteristic polynomial is negative, and there is a positive real eigenvalue.
    For the existence statement, $\zeta \mapsto c\,\kfrags\,J(\zeta)H(\zeta) - \mu_F$ is continuous, negative at $q_{\min}$ and positive at $q_{\max}$ by hypothesis, and negative for large $\zeta$ because $J \leq \mu_I$ and $H(\zeta)\to0$; so it has a root in $(q_{\min}, q_{\max})$, which is unstable by the above, and a root in $(q_{\max}, \infty)$.

    \emph{(iii).} For fixed $z$, the root $q_I^\circ$, and with it $J_s(z) = \sigma_I(z)q_I^\circ$, is continuous and increasing in $s$, and $J_s(z) \to 0$ as $s \to 0$.
    Let $s_* := \inf\{s \in (0,\bar s] : c\,\kfrags\,J_s(q_{\max})H(q_{\max}) \geq \mu_F\}$; by hypothesis and continuity, $s_* \in (0, \bar s)$ and $c\,\kfrags\,J_{s_*}(q_{\max})H(q_{\max}) = \mu_F$.
    Since $J_{s_*}$ is non-decreasing and $H(q_{\min}) < H(q_{\max})$, we have $c\,\kfrags\,J_{s_*}(q_{\min})H(q_{\min}) < \mu_F$.
    By continuity in $s$, both inequalities of~(ii) hold for all $s$ in an interval $(s_*, s_* + \delta)$.
\end{proof}

\begin{proof}[Proof of Proposition~\ref{prop:debris}]
    The fixed-point equation is $c\,\kfrags\,\lambda\,G(q_F) = \mu_F$ with $G = \rho H$.
    Differentiate: $G' = \rho'H + \rho H'$. Off $\cI$ we have $H' \le 0$ and, since $\rho'<0$ and $H>0$, $\rho'H < 0$, so $G' < 0$. On $\cI$ the displayed inequality is exactly $\rho'H + \rho H' < 0$. Hence $G' < 0$ on $[0,\infty)$, so $G$ is strictly decreasing and the level set $\{G = \mu_F/(c\,\kfrags\,\lambda)\}$ has at most one point.
    Finally, $\{H' > 0\}$ is bounded because $H'$ has the sign of the quartic $\pi$ of the Proof of Theorem~\ref{thm:bifurcation-structure}, whose leading coefficient is negative.
\end{proof}

\begin{proof}[Proof of Lemma~\ref{lem:clip-linearization}]
    At a fixed point of~\eqref{eq:cFL-1}--\eqref{eq:cFL-3}, the  identities $\kfrags\,c\,\alpha(z)q^*_I = \theta\mu_F$, $\kfrags\,c\,q^*_D = (1-\theta)\mu_F$ and $(1-P)\mu_I q^*_I = \sigma_D q^*_D$ of the proof of Lemma~\ref{lem:linearization} hold, with the admitted rate $\lambda\rho(z)$ in place of $\lambda$.
    The Jacobian differs from the one computed there only in its $(1,3)$ entry, which acquires the extra term $\lambda\rho'(z) = -\psi_\rho\,\lambda\rho(z)$.
    In particular the trace, and so $c_2$, is unchanged, and the upper-left $2\times2$ block is unchanged.

    For $c_0^\rho$: as in the proof of Lemma~\ref{lem:linearization}, $\det \cL = \sigma_I\sigma_D\,g'(z)$, where now $(x(\zeta), y(\zeta)) := \big(\lambda\rho(\zeta)/\sigma_I(\zeta),\ (1-P)\mu_I x(\zeta)/\sigma_D(\zeta)\big)$ is the unique point at which the right-hand sides of~\eqref{eq:cFL-1} and~\eqref{eq:cFL-2} vanish when $q_F = \zeta$, and $g(\zeta) = \zeta\big(c\,\kfrags\,\lambda\,G(\zeta) - \mu_F\big)$ is the right-hand side of~\eqref{eq:cFL-3} evaluated at $(x(\zeta), y(\zeta), \zeta)$; at the fixed point $c\,\kfrags\,\lambda\,G(z) = \mu_F$, so $g'(z) = \mu_F z G'(z)/G(z)$.

    For the Routh--Hurwitz gap: write $\delta := \psi_\rho\,\lambda\rho(z)$, $a := \alpha(z)cz$, $b := cz$, $p := (1-P)\mu_I$.
    The extra term $-\delta$ in $\cL_{13}$ changes $c_1$ by $-(-\delta)\cL_{31} = \delta\,\kfrags\, a$ and, expanding $\det \cL$ along the first row, changes $c_0 = -\det \cL$ by $\delta\,(\cL_{21}\cL_{32} + \sigma_D \cL_{31}) = \delta\,\kfrags\,(pb + \sigma_D a)$.
    Hence $c_1^\rho c_2 - c_0^\rho$ exceeds the corresponding uncontrolled expression $c_1c_2 - c_0$ of Lemma~\ref{lem:linearization} by
    \[
        \delta\,\kfrags\big(a(\sigma_I + \sigma_D - w) - pb - \sigma_D a\big) = \delta\,\kfrags\,c\,z\big(\alpha(\sigma_I - w) - p\big).
    \]
    At the fixed point $\kfrags\,c\,\lambda\rho(z) = \mu_F/H(z)$, so $\delta\,\kfrags\,c = \psi_\rho\,\mu_F/H$; and from $H = (\alpha\sigma_D + p)/(\sigma_I\sigma_D)$ and $\theta = \alpha\sigma_D/(\alpha\sigma_D + p)$ we get $\alpha/H = \theta\sigma_I$ and $p/H = (1-\theta)\sigma_I\sigma_D$.

    \emph{(i) and (ii)} follow as in Lemma~\ref{lem:linearization}(i) and~(ii).

    \emph{(iii).} Under the rate condition, the proof of Lemma~\ref{lem:linearization}(iii) gives $c_2 > 0$, $w \leq \theta\min\{\mu_I,\mu_D\}$, and the lower bound $c_1c_2 - c_0 \geq (\sigma_I + \sigma_D)(\sigma_I - w)(\sigma_D - w)$ (the second term of~\eqref{eq:RH-gap} is non-negative).
    Moreover $\sigma_I - w \geq \sigma_I - \mu_I = \alpha c z$, $\sigma_D - w \geq \sigma_D - \mu_D = cz$, and $\sigma_D - w \geq b + (1-\theta)\mu_D \geq (1-\theta)\sigma_D$.
    Write $\psi := \psi_\rho(z)$ and suppose $|\psi|\mu_F \leq c\alpha$.
    If $\psi \geq 0$, drop the non-negative term $\psi\mu_F z\sigma_I\theta(\sigma_I - w)$:
    \[
        c_1^\rho c_2 - c_0^\rho \geq (1-\theta)\sigma_D\big((\sigma_I + \sigma_D)(\sigma_I - w) - \psi\mu_F z\sigma_I\big)
        \geq (1-\theta)\sigma_D\big(\sigma_D(\sigma_I - w) + \sigma_I z(c\alpha - \psi\mu_F)\big) > 0,
    \]
    using $\theta < 1$ (as $P < 1$).
    If $\psi < 0$, drop the non-negative term $|\psi|\mu_F z\sigma_I(1-\theta)\sigma_D$:
    \[
        c_1^\rho c_2 - c_0^\rho \geq (\sigma_I - w)\big((\sigma_I + \sigma_D)(\sigma_D - w) - |\psi|\mu_F z\sigma_I\theta\big)
        \geq (\sigma_I - w)\big(\sigma_D(\sigma_D - w) + \sigma_I z(c - |\psi|\mu_F)\big) > 0,
    \]
    using $\theta \leq 1$ and $|\psi|\mu_F \leq c\alpha \leq c$.
\end{proof}

\begin{proof}[Proof of Theorem~\ref{prop:heaviest}]
    \emph{(i).} The function $G = \rho H$ is continuous and positive, and $G(q) \to 0$ as $q \to \infty$, because $\rho \leq 1$ and $H(q) \to 0$.
    If $G$ is non-increasing, every level set $\{G = L\}$ is an interval, possibly a single point or empty, so $\rho$ is fold-free; if $G$ is strictly decreasing, every level set has at most one point.
    Conversely, if $G$ is not non-increasing, there are $q_1 < q_2$ with $G(q_1) < G(q_2)$; for $L \in (G(q_1), G(q_2))$, $G - L$ is negative at $q_1$, positive at $q_2$, and negative for large $q$, so $\{G = L\}$ meets both $(q_1,q_2)$ and $(q_2, \infty)$ but not $q_2$, and it is not connected; with $\lambda := \mu_F/(c\,\kfrags\,L)$, $\rho$ is not fold-free.
    Finally, a monostable throttle is fold-free, so its $G$ is non-increasing, and as no level set contains an interval, $G$ is strictly decreasing.

    \emph{(ii).} $\rho^\star$ is admissible by Lemma~\ref{lem:rho-lip}, and $\rho^\star H = \bar H$ is non-increasing, so $\rho^\star$ is fold-free.
    If $\rho$ is fold-free, then $G = \rho H$ is non-increasing and $G \le H$, so for every $s \le q$, $G(q) \le G(s) \le H(s)$, and therefore $G(q) \le \bar H(q)$; dividing by $H(q)$ gives $\rho \le \rho^\star$.

    \emph{(iii).} In the two-fold regime, $G = \bar H$ is continuous and non-increasing with $\bar H(0) = H(0)$; it is strictly decreasing on $[0, q_{\min}]$ and on $[q^\ddagger, \infty)$, constant equal to $H(q_{\min})$ on the window, and tends to $0$.
    The equation $\bar H(z) = \mu_F/(c\,\kfrags\,\lambda)$ therefore has no solution $z > 0$ if $\lambda \leq \lambda_\crit$, the unique solution $z \in (0, q_{\min})$ if $\lambda_\crit < \lambda < \lambda_\flat$, the window if $\lambda = \lambda_\flat$, and the unique solution $z \in (q^\ddagger,\infty)$ if $\lambda > \lambda_\flat$.
    In the second and fourth of these cases, $z < q_{\min}$ or $z > q^\ddagger$, so $\rho^\star \equiv 1$ near $z$, so $\psi_{\rho^\star}(z) = 0$ and $G'(z) = H'(z) < 0$; under the rate condition~\eqref{eq:rate-condition}, Lemma~\ref{lem:clip-linearization}(iii) gives local asymptotic stability.

    \emph{(iv).} If $H$ is non-increasing then $\bar H = H$ and $\rho^\star \equiv 1$.
    Otherwise there are $s < q$ with $H(s) < H(q)$, so $\bar H(q) < H(q)$; the set $\{\bar H < H\}$ is open and non-empty, and $\bar H$ is constant on each of its components (near a point where $H > \bar H$, the infimum defining $\bar H$ is attained to the left), so $\rho^\star H = \bar H$ is not strictly decreasing and $\rho^\star$ is not monostable.
    Every monostable throttle is fold-free, so $\rho \le \rho^\star$ by~(ii).
    The tilted clip $\rho^\star e^{-\eta q}$ is admissible and monostable, since $G = \bar H e^{-\eta q}$ is strictly decreasing (proof of Theorem~\ref{thm:tilted-clip}), and it increases to $\rho^\star$ as $\eta \downarrow 0$.
    If there were a largest monostable throttle $\hat\rho$, then $\hat\rho \geq \rho^\star e^{-\eta q}$ for every $\eta > 0$, so $\hat\rho \geq \rho^\star$, hence $\hat\rho = \rho^\star$, which is not monostable: a contradiction.

    For the closed form, $\inf_{s\le q_F}H$ tracks $H$ down to its minimizer (at $q_{\min}$ in the two-fold regime, at the origin in the one-fold regime), stays flat at that value while $H$ exceeds it---i.e.\ across $[q_{\min}, q^\ddagger]$---and rejoins $H$ beyond $q^\ddagger$; dividing by $H(q_F)$ gives $\rho^\star$.
    The deepest cut is at the maximizer $q_{\max}$, where $\rho^\star = H(q_{\min})/H(q_{\max})$.
\end{proof}

\begin{proof}[Proof of Theorem~\ref{thm:tilted-clip}]
    Throughout, $H : [0,\infty) \to (0,\infty)$ is continuous with $H(q_F) \to 0$ as $q_F \to \infty$, and $\bar H(q_F) = \inf_{s \le q_F} H(s)$; no other property of $H$ is used for the fold and plateau statements.

    \emph{Admissibility.} $\rho_\eta = (\bar H/H)\,e^{-\eta q_F}$ is a product of factors in $(0,1]$ (as $\bar H \le H$ and $\eta > 0$), so $0 < \rho_\eta \le 1$; and it is Lipschitz, since $\bar H/H = \rho^\star$ is Lipschitz and bounded (Lemma~\ref{lem:rho-lip}), and $e^{-\eta q_F}$ is Lipschitz and bounded. Thus Theorem~\ref{thm:clip-fluid} applies, and the non-trivial fixed points solve $c\,\kfrags\,\lambda\,G_\eta(q_F) = \mu_F$ with $G_\eta = \rho_\eta H = \bar H\,e^{-\eta q_F}$.

    The key claim is that $G_\eta$ is a continuous, strictly decreasing bijection $[0,\infty) \to (0, H(0)]$.
    As the running infimum of a continuous function, $\bar H$ is continuous and non-increasing, and $\bar H(x) = \inf_{s \le x} H(s) > 0$ because $H$ attains a positive minimum on the compact $[0,x]$; hence $G_\eta$ is continuous and positive, with $G_\eta(0) = \bar H(0) = H(0)$ and $G_\eta(q_F) \to 0$. For $x_1 < x_2$,
    \[
        G_\eta(x_2) = \bar H(x_2)\,e^{-\eta x_2} \;\le\; \bar H(x_1)\,e^{-\eta x_2}
        \;<\; \bar H(x_1)\,e^{-\eta x_1} = G_\eta(x_1),
    \]
    the first inequality by $\bar H(x_2) \le \bar H(x_1)$ with $e^{-\eta x_2} > 0$, the \emph{strict} one by $\bar H(x_1) > 0$ with $e^{-\eta x_2} < e^{-\eta x_1}$. So $G_\eta$ is strictly decreasing; being continuous with $G_\eta(0) = H(0)$ and limit $0$, it is a bijection onto $(0, H(0)]$. The exponential is what supplies strictness where $\bar H$ is flat (the fold window $[q_{\min}, q^\ddagger]$), which the extremal clip $G = \bar H$ lacks.

    \emph{(i).} The non-trivial fixed-point equation reads $G_\eta(q_F) = \mu_F/(c\,\kfrags\,\lambda)$; the right side lies in $(0, H(0))$ if and only if $\lambda > \lambda_\crit = \mu_F/(c\,\kfrags\,H(0))$, and then the bijection gives exactly one solution $q_F^* > 0$; for $\lambda \leq \lambda_\crit$ there is none.
    At the clean state $\rho_\eta(0) = 1$, so the controlled Jacobian there is the uncontrolled one, whose eigenvalues are $-\mu_I$, $-\mu_D$ and $\mu_F(\lambda/\lambda_\crit - 1)$ (Proof of Theorem~\ref{thm:bifurcation-structure}); the clean state is therefore locally asymptotically stable for $\lambda < \lambda_\crit$ and unstable for $\lambda > \lambda_\crit$.

    \emph{(ii).} The branch is the graph of $\cM_{\rho_\eta}(q_F) = \mu_F/(c\,\kfrags\,G_\eta(q_F))$. Since $G_\eta$ is strictly decreasing and positive, $\cM_{\rho_\eta}$ is continuous and strictly increasing, mapping $(0,\infty)$ onto $(\lambda_\crit, \infty)$; its inverse $q_F^*(\lambda)$ is therefore continuous and strictly increasing. A fold would require $\cM_{\rho_\eta}$ to be non-monotone, which it is not, so there are no saddle-nodes; the branch rises continuously from $0$ at $\lambda_\crit$, eliminating the plateau and jump of the extremal clip.

    \emph{(iii).} Where $\rho_\eta$ is differentiable, $G_\eta'/G_\eta = \bar H'/\bar H - \eta \le -\eta$ (as $\bar H' \le 0$), so Lemma~\ref{lem:clip-linearization} gives
    \[
        c_0^{\rho_\eta} = -\sigma_I\sigma_D\,\mu_F\,q_F^*\,\frac{G_\eta'(q^*_F)}{G_\eta(q^*_F)} \;\ge\; \eta\,\sigma_I\sigma_D\,\mu_F\,q_F^* \;>\; 0.
    \]
    The stability statement is Lemma~\ref{lem:clip-linearization}(iii), since $G_\eta' < 0$.
\end{proof}

\begin{proof}[Proof of Theorem~\ref{prop:clip-slow-passage}]
    Write $\mu := \lambda - \lambda_\flat$ and $\mu(t) = \lambda(t) - \lambda_\flat$.
    The right-hand side of~\eqref{eq:cFL-1}--\eqref{eq:cFL-3} is affine in $\lambda$: it equals $g(x) + \mu F_\lambda(x)$, where $g$ is the right-hand side at $\lambda = \lambda_\flat$ and $F_\lambda(x) := (\rho(q_F), 0, 0)$.

    \emph{(i) Extremal clip. Step 1: the fixed points of the window.}
    For $z \in [q_{\min}, q^\ddagger]$, $G(z) = \bar H(z) = H(q_{\min})$ and $c\,\kfrags\,\lambda_\flat H(q_{\min}) = \mu_F$, so the point $x^*(z)$ with $q_F = z$, $q_I = \lambda_\flat\rho^\star(z)/\sigma_I(z)$, $q_D = (1-P)\mu_I q_I/\sigma_D(z)$ is a fixed point of $g$ (Theorem~\ref{thm:clip-fluid}).
    On $(q_{\min}, q^\ddagger)$, $\rho^\star = H(q_{\min})/H$ is smooth, $G' = 0$, and $\psi_{\rho^\star} = H'/H$; so Lemma~\ref{lem:clip-linearization} and the hypothesis give $c_0^{\rho^\star} = 0$, $c_2 > 0$, and $c_1^{\rho^\star}c_2 > 0$.
    At $z = q_{\min}$, $\rho^{\star\prime}(q_{\min}) = -H(q_{\min})H'(q_{\min})/H(q_{\min})^2 = 0$, so $\rho^\star$ is $C^1$ there (with a Lipschitz derivative) and the Jacobian of $g$ at $x^*(q_{\min}) = q^\dagger$ is the uncontrolled Jacobian at the fold, for which $c_0 = 0$, $c_2 > 0$, $c_1 > 0$ by Lemma~\ref{lem:linearization}(iii).
    In both cases the characteristic polynomial is $s(s^2 + c_2s + c_1)$ with $c_1, c_2 > 0$, which gives the eigenvalue statement.
    For $\mu \neq 0$, Theorem~\ref{thm:clip-fluid} shows that the fixed point lies off the window, so no fixed point is destroyed.
    This proves the first part of~(i).

    \emph{Step 2: arrival at the window.}
    On $\{q_F \leq q_{\min}\}$ we have $\rho^\star = 1$, so the clipped and uncontrolled vector fields coincide there.
    Hence, up to the first time $t_0$ at which $q_F = q_{\min}$, the clipped solution coincides with the uncontrolled solution of Theorem~\ref{thm:slow-passage}, with the same initial condition and ramp; the hypotheses of Statements~(i)--(iii) of that theorem, whose proof does not use~(A2), hold because the rate condition implies~(A1) by Lemma~\ref{lem:linearization}(iii).
    In the scaling of Step~3 of its proof, the uncontrolled solution follows the special solution $U(s) = \mathrm{Ai}'(-s)/\mathrm{Ai}(-s)$ of $dU/ds = U^2 + s$ up to $o(1)$, uniformly on compact subsets of $(-\infty, \Omega_0)$.
    This solution is negative for $s < \Omega_1$ and vanishes first at $s = \Omega_1$, where $dU/ds = \Omega_1 > 0$, so it crosses zero transversally.
    Since $q_F - q_{\min} = u + O(u^2 + |\mu| + \varepsilon)$ on the center manifold $W^c$, and the solution is exponentially close to $W^c$, it follows that $t_0$ exists, that $\mu(t_0) = \Omega_1(\beta_\lambda\beta_2)^{-1/3}\varepsilon^{2/3}(1 + o(1))$, and that $|q(t_0) - q^\dagger| = O(\varepsilon^{2/3})$.
    In time, since $\mu(t) = \lambda(t) - \lambda_\flat = \varepsilon(t - t_\flat)$, this is $t_0 - t_\flat = \mu(t_0)/\varepsilon = \Omega_1(\beta_\lambda\beta_2)^{-1/3}\varepsilon^{-1/3}(1 + o(1))$; equivalently, by~\eqref{eq:sp-inner-time}, the inner time $s$ reaches $\Omega_1$.
    This proves~(i)(a).

    \emph{Step 3: the window as a normally attracting manifold.}
    Fix $z_b \in (q_{\min}, q^\ddagger)$ and $\bar m > 0$, and let $\delta := \varepsilon^{1/2}$ and $m := \mu/\delta$, so that the ramped system reads
    \[
        \dot x = g(x) + \delta\, m\, F_\lambda(x), \qquad \dot m = \delta.
    \]
    To accommodate base points slightly to the left of $q_{\min}$, and to have a smooth vector field, we work with the auxiliary throttle $\hat\rho := H(q_{\min})/H$ on $[0, q^\ddagger]$: it is smooth, it equals $\rho^\star$ on the window, so the auxiliary field $\hat g$ equals $g$ on $\{q_{\min} \leq q_F \leq q^\ddagger\}$, and on $\{q_F < q_{\min}\}$ it differs from $\rho^\star = 1$ by $1 - H(q_{\min})/H(q_F) = O((q_F - q_{\min})^2)$, because $H'(q_{\min}) = 0$.
    For $\hat g$ the curve of fixed points extends to $z \in [q_{\min} - \delta_1, z_b]$, for a small fixed $\delta_1 > 0$.
    At $\delta = 0$ the set
    \[
        M_0 := \big\{(x^*(z), m) : z \in [q_{\min} - \delta_1, z_b],\ |m| \leq \bar m\big\}
    \]
    is a compact two-dimensional manifold of fixed points of $\hat g$, and by Step~1 (and continuity of the eigenvalues in $z$, for $\delta_1$ small) the linearization at each of its points has a two-dimensional kernel, tangent to $M_0$, and two eigenvalues with negative real part bounded away from zero: $M_0$ is normally attracting.
    By Fenichel's theorem~\cite{fenichel1979}, applied after the usual modification of the vector field near $\partial M_0$, for small $\delta$ there is a locally invariant manifold $M_\delta$ within $O(\delta)$ of $M_0$, which attracts nearby solutions at an exponential rate with asymptotic phase.
    In the coordinates $(z, m)$, with $z = q_F$, the flow on $M_\delta$ is the projection of the perturbation onto the tangent space of $M_0$ along the stable directions, up to higher order:
    \[
        \frac{dz}{dt} = \delta\, m\, R(z) + O(\delta^2), \qquad \frac{dm}{dt} = \delta, \qquad R(z) := \ell(z)\cdot F_\lambda(x^*(z)),
    \]
    where $\ell(z)$ is the left null vector of the Jacobian at $x^*(z)$, normalized by $\ell\cdot r = 1$ against the tangent $r$ of the curve with $r_3 = 1$; the remainder is $O(\delta^2)$ because $\hat g$ is smooth.
    We take $\bar m$ larger than $2\Psi(z_b)^{1/2} + 1$, with $\Psi$ as in Step~4, so that the transit below takes place in $\{|m| \leq \bar m\}$.
    The left null vector involves only the first two columns of the Jacobian, which the throttle does not change, so the computation in Step~2 of the proof of Theorem~\ref{thm:slow-passage} applies and gives $\ell_1 = \kfrags\,c\,z\,H(z)\,\sigma_I\sigma_D/c_1^{\rho^\star}$, whence
    \[
        R(z) = \ell_1\,\rho^\star(z) = \kfrags\,c\,z\,H(q_{\min})\,\frac{\sigma_I\sigma_D}{c_1^{\rho^\star}} = \frac{\mu_F}{\lambda_\flat}\,\frac{\sigma_I(z)\sigma_D(z)}{c_1^{\rho^\star}(z)}\,z > 0.
    \]

    \emph{Step 4: the transit.}
    Let $\Psi(z) := \int_{q_{\min}}^{z} d\zeta/R(\zeta) = (\lambda_\flat/\mu_F)\,\cJ(z)$.
    In the time variable $t$, the flow on $M_\delta$ reads $dz/dt = \mu R(z) + O(\varepsilon)$, so along it
    \[
        \frac{d}{dt}\Psi(z(t)) = \mu(t) + O(\varepsilon),
    \]
    uniformly, since $R$ is bounded below by a positive constant.
    By Step~2 and the exponential attraction of $M_\delta$, after a time of order $\log(1/\varepsilon)$ following $t_0$, at a time $t_1$ with $\mu(t_1) = O(\varepsilon^{2/3})$, the solution is within $O(\varepsilon)$ of $M_\delta$, and its asymptotic phase is a solution on $M_\delta$ with $|z(t_1) - q_{\min}| = O(\varepsilon^{1/2})$ (as $M_\delta$ lies within $O(\delta)$ of $M_0$).
    The actual vector field differs from $\hat g$ only where $q_F < q_{\min}$, and on the part of the solution with $q_F \geq q_{\min} - O(\varepsilon^{1/2})$ by $O(\varepsilon)$.
    Since $\mu(t_1)$ is of order $\varepsilon^{2/3}$, $dz/dt \geq \mu R/2 \geq c\,\varepsilon^{2/3}$ from $t_1$ on, so this part of the solution lasts a time $O(\varepsilon^{-1/6})$, and the difference of the vector fields shifts $\Psi$ by $O(\varepsilon^{5/6})$.
    Since $d\Psi/dt \geq \mu/2 > 0$ after $t_1$, $z$ reaches $z_b$, at a time $\tau_b$ at which $\mu(\tau_b) = O(\varepsilon^{1/2})$, so that $\tau_b - t_1 = O(\varepsilon^{-1/2})$.
    Integrating from $t_1$ to $\tau_b$, with $\dot\mu = \varepsilon$,
    \[
        \Psi(z_b) - \Psi(z(t_1)) = \frac{\mu(\tau_b)^2 - \mu(t_1)^2}{2\varepsilon} + O\big(\varepsilon(\tau_b - t_1)\big),
    \]
    and $\Psi(z(t_1)) = O(\varepsilon^{1/2})$, $\mu(t_1)^2/(2\varepsilon) = O(\varepsilon^{1/3})$, and $\varepsilon(\tau_b - t_1) = O(\varepsilon^{1/2})$.
    Hence $\mu(\tau_b)^2 = 2\varepsilon\Psi(z_b) + O(\varepsilon^{4/3})$, that is, $\mu(\tau_b) = \sqrt{2\varepsilon\Psi(z_b)} + O(\varepsilon^{5/6})$.
    Finally, $q_F$ along the actual solution differs from $z$ along its asymptotic phase by a quantity which decays exponentially in $t - t_1$, while $dz/dt \geq \mu R/2$ is of order $\varepsilon^{1/2}$ near $\tau_b$; so the first time $t_b$ at which $q_F(t_b) = z_b$ differs from $\tau_b$ by an exponentially small amount.
    Therefore $\mu(t_b) = \mu(\tau_b) + \varepsilon(t_b - \tau_b) = \sqrt{2\varepsilon\Psi(z_b)} + O(\varepsilon^{5/6})$, which is~\eqref{eq:clip-delay}, since $\Psi(z_b) = (\lambda_\flat/\mu_F)\,\cJ(z_b)$.
    Dividing by $\varepsilon$, with $\mu(t) = \varepsilon(t - t_\flat)$, gives~\eqref{eq:clip-delay-time}: $t_b = t_\flat + \sqrt{2\Psi(z_b)}\,\varepsilon^{-1/2} + O(\varepsilon^{-1/6})$.
    By~(i)(a), $t_0 - t_\flat = O(\varepsilon^{-1/3})$, so
    \[
        t_b - t_0 = \sqrt{2\Psi(z_b)}\;\varepsilon^{-1/2}\,\big(1 + o(1)\big),
    \]
    and $\Psi(z_b) > 0$ because $R > 0$; the transit from $q_{\min}$ to $z_b$ therefore takes a time of order $\varepsilon^{-1/2}$.
    This proves~(i)(b).

    \emph{(ii) Tilted clip.}
    By Theorem~\ref{thm:tilted-clip}, $\rho_\eta$ is admissible and its branch of non-trivial fixed points is continuous and strictly increasing.
    On the window $G_\eta = \bar H e^{-\eta q_F} = H(q_{\min})e^{-\eta q_F}$, so the fixed-point relation $c\,\kfrags\,\lambda\,G_\eta(q^*_F) = \mu_F$ reads $\lambda = \lambda_\flat e^{\eta q^*_F}$, i.e.\ $q^*_F(\lambda) = \eta^{-1}\ln(\lambda/\lambda_\flat)$, which sweeps $[q_{\min}, q^\ddagger]$ as $\lambda$ runs over $[\lambda_\flat e^{\eta q_{\min}}, \lambda_\flat e^{\eta q^\ddagger}]$.
    For $\lambda$ in the open interval $\big(\lambda_\flat e^{\eta q_{\min}}, \lambda_\flat e^{\eta q^\ddagger}\big)$, $q^*_F(\lambda)$ lies in the open window, on which $\rho_\eta = H(q_{\min})e^{-\eta q_F}/H$ is smooth; at $q^*_F = q^\ddagger$ it has a corner, inherited from $\rho^\star$, and Lemma~\ref{lem:clip-linearization} does not apply there.
    At the interior points $G_\eta'/G_\eta = -\eta$, so Lemma~\ref{lem:clip-linearization} gives $c_0^{\rho_\eta} = \eta\,\sigma_I\sigma_D\,\mu_F\,q^*_F \geq \eta\,\sigma_I\sigma_D\,\mu_F\,q_{\min} > 0$.
    For $\lambda \in [\lambda_a, \lambda_b]$ the fixed points have $q^*_F$ in a compact subset of the open window, on which $\rho_\eta$ is smooth; as they are locally asymptotically stable, this piece of the branch is a normally hyperbolic, attracting critical manifold of the ramped system, and Fenichel's theorem~\cite{fenichel1979}, as in Step~1 of the proof of Theorem~\ref{thm:slow-passage}, gives tracking with an $O(\varepsilon)$ lag and no delay.
    In time: $q^*_F(\lambda)$ is strictly increasing with $dq^*_F/d\lambda = 1/(\eta\lambda) \geq 1/(\eta\lambda_b)$ on $[\lambda_a, \lambda_b]$, so if $|q_F(t) - q^*_F(\lambda(t))| \leq C\varepsilon$, then $q_F(t) = q^*_F(\lambda(t - \tau))$ for some $\tau$ with $|\tau| \leq C\eta\lambda_b$ (away from the ends of the interval); the fragment count lags the moving fixed point by a time which is bounded uniformly in $\varepsilon$.
    This proves~(ii).
    Finally, when $\eta$ is small the slowest eigenvalue on the branch is of order $c_0^{\rho_\eta}/c_1^{\rho_\eta} = O(\eta)$, so the lag is of order $\varepsilon/\eta$; equating it with the overshoot $O(\varepsilon^{1/2})$ of~(i) gives the crossover $\eta \sim \varepsilon^{1/2}$ of Remark~\ref{rem:hierarchy}.
\end{proof}

\begin{proof}[Proof of Proposition~\ref{prop:throughput}]
    By Theorem~\ref{thm:clip-fluid}, a non-trivial fixed point satisfies $c\,\kfrags\,\lambda\rho(q^*_F)H(q^*_F) = \mu_F$, so the admitted launch rate is $\lambda\rho(q^*_F) = \cM(q^*_F)$, a function of $q^*_F$ alone, and
    \[
        T_{\mathrm{net}} = P\mu_I q^*_I = P\mu_I\,\frac{\lambda\rho(q^*_F)}{\sigma_I(q^*_F)} = \frac{P\mu_I\,\mu_F}{c\,\kfrags\,\sigma_I(q^*_F)H(q^*_F)} = \Phi(q^*_F),
    \]
    by the identity $\sigma_I H = \alpha + (1-P)\mu_I/\sigma_D$. Neither expression refers to $\rho$ except through $q_F^*$; this is the claimed policy-independence.

    \emph{Limits.} As $q \to \infty$, $\alpha(q) \to 1$ and $(1-P)\mu_I/(\mu_D + cq) \to 0$, so $\Phi(q) \to P\mu_I\mu_F/(c\,\kfrags) = L$, and $q^*_I = \Phi(q^*_F)/(P\mu_I) \to \mu_F/(c\,\kfrags)$.

    \emph{The two clips.} For both clips the fixed-point function $G$ is continuous, positive and tends to $0$, so the solution of $c\,\kfrags\,\lambda\,G(q^*_F) = \mu_F$ on the branch satisfies $q_F^*(\lambda) \to \infty$ as $\lambda \to \infty$. Beyond the window $G = H$ for the extremal clip and $G = He^{-\eta q_F}$ for the tilted clip, and $H(q) = (1 + o(1))/(cq)$ as $q \to \infty$, the derelict term of $H$ being $O(q^{-2})$. Hence the extremal relation $c\,\kfrags\,\lambda H(q^*_F) = \mu_F$ gives $\kfrags\lambda/q^*_F \sim \mu_F$, i.e.\ $q_F^* \sim (\kfrags/\mu_F)\,\lambda$, and the tilted relation gives $q^*_F e^{\eta q^*_F} \sim (\kfrags/\mu_F)\,\lambda$, i.e.\ $q_F^* \sim \eta^{-1}\log\lambda$.
    If $\mu_F < c\,\kfrags$, then $q^*_I \to \mu_F/(c\,\kfrags) < 1$, so for large $\lambda$ these fixed points lie strictly below the wall $q_I = 1$; there the reflection is inactive and they are fixed points of the reflected controlled fluid limit of Theorem~\ref{thm:clip-fluid}, at which the computation above applies. Therefore $T_{\mathrm{net}}(\lambda) = \Phi(q_F^*(\lambda)) \to L$ for both clips.
\end{proof}

\begin{proof}[Proof of Theorem~\ref{thm:bang-bang}]
    \emph{(i).} Since $\alpha'(q) = (1-\alpha_0)\Qsat/(q + \Qsat)^2$,
    \[
        \varphi'(q) = \frac{(1-\alpha_0)\Qsat}{(q + \Qsat)^2} - \frac{(1-P)\mu_I\,c}{\sigma_D(q)^2}.
    \]
    Both terms are positive, so $\varphi'(q)$ has the sign of $\sqrt{(1-\alpha_0)\Qsat}\,\sigma_D(q) - \sqrt{(1-P)\mu_I c}\,(q + \Qsat)$, that is, of
    \[
        \sigma_D(q) - c\,r\,(q + \Qsat) = c\,(1-r)\,(q - q^\circ).
    \]
    Under~\eqref{eq:unimodal}, $1 - r > 0$ and $q^\circ > 0$, so $\varphi$ is strictly decreasing on $[0, q^\circ]$ and strictly increasing on $[q^\circ, \infty)$, which is the claim about $\Phi = P\mu_I\mu_F/(c\,\kfrags\,\varphi)$.
    Since $\sigma_I' = c(\alpha + q\alpha') > 0$, the function $\cM = \Phi\sigma_I/(P\mu_I)$ is a product of two positive, strictly increasing functions on $[0, q^\circ]$, and $\cM(0) = \mu_F/(c\,\kfrags\,H(0)) = \lambda_\crit$.
    Finally $H = \varphi/\sigma_I$, so $H' = (\varphi'\sigma_I - \varphi\sigma_I')/\sigma_I^2 < 0$ on $[0, q^\circ]$.
    In particular $H'(0) < 0$, which excludes the one-fold regime, and in the two-fold regime the first positive zero of $H'$, which is $q_{\min}$, exceeds $q^\circ$.

    \emph{(ii).} Let $(q^*, J)$ be an equilibrium under admission. The equations $b^J_I(q^*) = 0$ and $b^J_D(q^*) = 0$ give $q^*_I = J/\sigma_I(q^*_F)$ and $q^*_D = (1-P)\mu_I q^*_I/\sigma_D(q^*_F)$, and if $q^*_F > 0$, then $b_F(q^*) = 0$ divided by $q^*_F$ gives $c\,\kfrags\,J H(q^*_F) = \mu_F$, so that, as in the proof of Proposition~\ref{prop:throughput},
    \[
        J = \cM(q^*_F), \qquad T_{\mathrm{net}} = \frac{P\mu_I\,J}{\sigma_I(q^*_F)} = \Phi(q^*_F).
    \]
    Now suppose that $(q^*, J)$ is debris-stable and $q^*_F \le \bar q$, and let $\hat q := \min\{\bar q, q^\circ\}$ and, when $\lambda > \lambda_\crit$, $z := \min\{z_\lambda, \hat q\} \in (0, q^\circ]$.

    If $q^*_F = 0$, then $q^*_I = J/\mu_I$ and $T_{\mathrm{net}} = PJ$, with $J \le \min\{\lambda, \lambda_\crit\}$. If $\lambda \le \lambda_\crit$, this is at most $P\lambda = T^\star(\lambda, \hat q)$; if $\lambda > \lambda_\crit$, it is at most $P\lambda_\crit = \Phi(0) < \Phi(z)$, by~(i).

    If $y := q^*_F > 0$, then $\cM(y) = J \le \lambda$, $T_{\mathrm{net}} = \Phi(y)$, and, since $\Phi(y) = P\mu_I\cM(y)/\sigma_I(y)$,
    \begin{equation}\label{eq:Tnet-bound}
        T_{\mathrm{net}} \le \frac{P\mu_I\,\lambda}{\sigma_I(y)}.
    \end{equation}
    Suppose first that $y \le q^\circ$. As $\cM$ is strictly increasing on $[0, q^\circ]$ with $\cM(0) = \lambda_\crit$, the inequality $\cM(y) \le \lambda$ forces $\lambda > \lambda_\crit$ and $y \le z_\lambda$ (if $\lambda < \cM(q^\circ)$ because $\cM(z_\lambda) = \lambda$, and otherwise because $z_\lambda = q^\circ$). Since also $y \le \bar q$, we get $y \le z$, and so $T_{\mathrm{net}} = \Phi(y) \le \Phi(z)$ by~(i).
    Suppose next that $y > q^\circ$, which requires $\bar q > q^\circ$, so that $\hat q = q^\circ$ and $z = z_\lambda$.
    If $\lambda \le \lambda_\crit$, then~\eqref{eq:Tnet-bound} and $\sigma_I(y) > \mu_I$ give $T_{\mathrm{net}} < P\lambda$.
    If $\lambda_\crit < \lambda < \cM(q^\circ)$, then $z_\lambda < q^\circ < y$, and~\eqref{eq:Tnet-bound} gives $T_{\mathrm{net}} < P\mu_I\lambda/\sigma_I(z_\lambda) = P\mu_I\cM(z_\lambda)/\sigma_I(z_\lambda) = \Phi(z_\lambda)$.
    If $\lambda \ge \cM(q^\circ)$, then $T_{\mathrm{net}} = \Phi(y) < \Phi(q^\circ) = \Phi(z)$ by~(i).
    In every case $T_{\mathrm{net}} \le T^\star(\lambda, \hat q)$.

    \emph{(iii).} Let $q \equiv q^*$ be an equilibrium of~\eqref{eq:BB}, with threshold $\hat q$ and control $u$.
    If $q^*_F > \hat q$, then $u = 0$ almost everywhere, so the intact equation gives $q^*_I = 0$, the derelict equation then gives $q^*_D = 0$, and $b_F(q^*) = -\mu_F q^*_F \neq 0$, which is impossible.
    If $q^*_F < \hat q$, then $u = 1$ almost everywhere, so $(q^*, \lambda)$ is an equilibrium under admission: either the clean state, which is debris-stable exactly when $\lambda \le \lambda_\crit$, or a point with $\cM(q^*_F) = \lambda$ and $0 < q^*_F < \hat q \le q^\circ$, which by~(i) exists exactly when $\lambda_\crit < \lambda < \cM(\hat q)$, and then $q^*_F = z_\lambda$.
    If $q^*_F = \hat q$, the intact equation shows that $\lambda u(t) = \sigma_I(\hat q)\,q^*_I =: J$ for almost every $t$, so $(q^*, J)$ is an equilibrium under admission, and by the first paragraph of the proof of~(ii), $J = \cM(\hat q) \le \lambda$ and $q^* = x^*(\hat q)$; conversely, if $\cM(\hat q) \le \lambda$, then $x^*(\hat q)$ is a constant solution of~\eqref{eq:BB}, with $u \equiv \cM(\hat q)/\lambda \in (0,1]$.
    As $\cM(\hat q) > \lambda_\crit$, these are the three cases of the statement, and $T_{\mathrm{net}}$ equals $P\lambda$, $\Phi(z_\lambda)$ and $\Phi(\hat q)$ in them.
    This is $T^\star(\lambda, \hat q)$, because $\min\{z_\lambda, \hat q\}$ equals $z_\lambda$ when $\lambda < \cM(\hat q)$ and equals $\hat q$ when $\lambda \ge \cM(\hat q)$.
\end{proof}

\begin{proof}[Proof of Theorem~\ref{thm:bang-bang-fluid}]
    \emph{(i).} The rule changes only the launch intensity, which becomes $a(t) = \int_0^t n\lambda\,\ind\{Q_F(s) < n\bar q\}\ind\{Q_I(s) < n\}\,ds \le n\lambda t$.
    As in the proof of Theorem~\ref{thm:clip-fluid}, Lemma~\ref{lem:apriori} and Theorem~\ref{thm:strong-approx} therefore hold for $Q^{n,\bar q}$, and dividing by $n$ gives~\eqref{eq:prelimit} with its first line replaced by
    \[
        \bar Q^n_I(t) = \bar Q^n_I(0) + \lambda U^n(t) - \int_0^t \sigma_I(\bar Q^n_F)\,\bar Q^n_I\,ds - \ell^n(t) + \varepsilon^n_I(t),
        \qquad
        \ell^n(t) := \lambda\int_0^t \ind\{\bar Q^n_F < \bar q,\ \bar Q^n_I = 1\}\,ds.
    \]
    Almost surely, for every $T \in \bN$, $\|\varepsilon^n\|_T \to 0$ and, for all large $n$, $\bar Q^n$ takes values on $[0,T]$ in the compact set $K_T$ of the proof of Proposition~\ref{prop:reflected-fluid}.
    Fix an outcome in this event, and let $C_T$ bound $|b|$ on $K_T$, where $b$ is the uncontrolled vector field of Section~\ref{sec:limiting}.
    For all large $n$, the continuous processes $Y^n := \bar Q^n - \varepsilon^n$, $U^n$ and $\ell^n$ are Lipschitz on $[0,T]$, with a constant which does not depend on $n$, since $|\sigma_I(q_F)q_I| = |\lambda - b_I(q)| \le \lambda + C_T$ on $K_T$.
    By the Arzel\`a--Ascoli theorem and a diagonal argument over $T \in \bN$, every subsequence has a further subsequence along which $(Y^n, U^n, \ell^n) \to (q, U, \ell)$ uniformly on compacts, and then $\bar Q^n \to q$ as well.
    Since $b$ is Lipschitz on $K_T$, we may pass to the limit in the integrals: $q$ satisfies the second line of~\eqref{eq:BB}, and
    \[
        q_I(t) = q_I(0) + \lambda U(t) - \int_0^t \sigma_I(q_F)\,q_I\,ds - \ell(t).
    \]
    The limit $U$ is non-decreasing and $1$-Lipschitz, so $U = \int_0^\cdot u\,ds$ for a measurable $u$ with values in $[0,1]$.
    If $q_F \le \bar q - \delta$ on an interval $[t_1, t_2]$ for some $\delta > 0$, then $\bar Q^n_F < \bar q$ on $[t_1, t_2]$ for all large $n$ in the subsequence, so $U^n(t_2) - U^n(t_1) = t_2 - t_1$, and in the limit $U(t_2) - U(t_1) = t_2 - t_1$.
    The open set $\{t : q_F(t) < \bar q\}$ is a countable union of such intervals, so $u = 1$ almost everywhere on it, and in the same way $u = 0$ almost everywhere on $\{t : q_F(t) > \bar q\}$.
    In the same way, $\ell$ is non-decreasing and constant on every interval on which $q_I < 1$, and $q_I \le 1$; on the level set $\{q_I = 1\}$, where $\dot q_I = 0$ almost everywhere, $\dot\ell = \lambda u - \sigma_I(q_F) \le \lambda - \mu_I \le 0$ almost everywhere.
    Hence $\dot\ell = 0$ almost everywhere, $\ell \equiv 0$, and $q$ solves~\eqref{eq:BB}.

    \emph{Two identities.} Let $q$ be a solution of~\eqref{eq:BB} with control $u$, and let $y := q_F - \bar q$ and $w := \xi(q) - \bar\xi$.
    These are Lipschitz functions of time, and almost everywhere
    \begin{equation}\label{eq:BB-yw}
        \dot y = c\,\kfrags\,q_F\,w,
        \qquad
        \dot w = (1-u)\,g_0(q) + u\,g_1(q),
    \end{equation}
    where, for $v \in \{0, 1\}$,
    \[
        g_v(q) := \alpha'(q_F)\,b_F(q)\,q_I + \alpha(q_F)\big(\lambda v - \sigma_I(q_F)\,q_I\big) + b_D(q)
    \]
    is the rate of change of $\xi$ along the vector field with $u \equiv v$.
    On $\Sigma_0$, where $b_F = 0$, the definition of $u_\Sigma$ gives
    \begin{equation}\label{eq:g-on-Sigma0}
        g_v(q) = \lambda\,\alpha(\bar q)\big(v - u_\Sigma(q)\big), \qquad q \in \Sigma_0.
    \end{equation}

    \emph{(a) Motion on $\Sigma_0$.}
    For almost every $t$ with $q(t) \in \Sigma_0$, we have $y(t) = w(t) = 0$, hence $\dot y(t) = \dot w(t) = 0$, since the derivative of a Lipschitz function vanishes almost everywhere on each of its level sets; then~\eqref{eq:BB-yw} and~\eqref{eq:g-on-Sigma0} give $u(t) = u_\Sigma(q(t))$.
    If $q(t) \in \Sigma_0$ for all $t$ in an interval, then on that interval $q_F = \bar q$, $q_D = \bar\xi - \alpha(\bar q)\,q_I$, and
    \[
        \dot q_I = \lambda u_\Sigma(q) - \sigma_I(\bar q)\,q_I = -\frac{b_D(q)}{\alpha(\bar q)} = -\omega\,\big(q_I - q_I^*(\bar q)\big),
    \]
    with $\omega$ and $q_I^*(\bar q)$ as in~(iv): substituting $q_D = \bar\xi - \alpha(\bar q)q_I$ gives $b_D(q) = \alpha(\bar q)\,\omega\,q_I - \sigma_D(\bar q)\,\bar\xi$, and $\sigma_D(\bar q)\,\bar\xi/\big(\alpha(\bar q)\,\omega\big) = \bar\xi/\varphi(\bar q) = \cM(\bar q)/\sigma_I(\bar q)$.
    Along $\Sigma_0$, $u_\Sigma$ is an affine function of $q_I$, which takes the value $\cM(\bar q)/\lambda$ at $q_I = q_I^*(\bar q)$, where $b_D = 0$.

    \emph{(b) Capture by $\Sigma_0$.}
    Suppose that $p := q(t_0) \in \Sigma_{\mathrm{sl}}$. We show that $q(t) \in \Sigma_0$ for all $t$ in some interval $[t_0, t_0 + \delta]$.
    By~\eqref{eq:g-on-Sigma0}, $g_0(p) = -\lambda\alpha(\bar q)u_\Sigma(p) < 0 < \lambda\alpha(\bar q)\big(1 - u_\Sigma(p)\big) = g_1(p)$.
    As $g_0, g_1$ are continuously differentiable, there is a neighborhood $N$ of $p$ on which $q_F > 0$, and on which $\beta := -g_0/(c\,\kfrags\,q_F)$ and $\gamma := g_1/(c\,\kfrags\,q_F)$ are bounded below by some $m > 0$ and have bounded gradients.
    On $N$, let
    \[
        V(q) := \tfrac12\big(\xi(q) - \bar\xi\big)^2 + \beta(q)\,(q_F - \bar q)^+ + \gamma(q)\,(q_F - \bar q)^-,
    \]
    which is Lipschitz and non-negative, and vanishes exactly on $\Sigma_0 \cap N$.
    Choose $\delta > 0$ such that $q(t) \in N$ for $t \in [t_0, t_0 + \delta]$; then $V(q(t))$ is absolutely continuous on this interval.
    For almost every $t$ with $y(t) > 0$, we have $u(t) = 0$, and by~\eqref{eq:BB-yw}
    \[
        \frac{d}{dt}V(q(t)) = w\,g_0(q) + \beta(q)\,c\,\kfrags\,q_F\,w + y\,\nabla\beta(q)\cdot\dot q = y\,\nabla\beta(q)\cdot\dot q,
    \]
    the first two terms canceling by the choice of $\beta$; in the same way, for almost every $t$ with $y(t) < 0$, $u(t) = 1$ and $\frac{d}{dt}V(q(t)) = -y\,\nabla\gamma(q)\cdot\dot q$.
    For almost every $t$ with $y(t) = 0$, $\dot y(t) = 0$, so $w(t) = 0$ and $V(q(t)) = 0$, and the derivative of the non-negative function $V(q(\cdot))$ vanishes at $t$.
    Since $\dot q$ is bounded and $V \ge m|y|$ on $N$, there is $C$ such that $\frac{d}{dt}V(q(t)) \le C|y(t)| \le (C/m)\,V(q(t))$ for almost every $t \in [t_0, t_0+\delta]$, and Gr\"onwall's inequality and $V(p) = 0$ give $V(q(t)) = 0$, that is, $q(t) \in \Sigma_0$, on $[t_0, t_0 + \delta]$.

    \emph{(ii).} Let $q$ and $\widetilde q$ be solutions of~\eqref{eq:BB} with $q(0) = \widetilde q(0)$ which do not coincide, let $\tau := \sup\{t : q = \widetilde q \text{ on } [0,t]\}$, and let $p := q(\tau) = \widetilde q(\tau)$. We show that $p \in \Sigma_\partial$.
    If $p_F \neq \bar q$, then both solutions stay on the same side of $\Sigma$ on some interval $[\tau, \tau + \delta]$, where their controls are both equal to $1$, or both equal to $0$, almost everywhere; so both solve the same differential equation with a locally Lipschitz right-hand side, and they coincide on $[\tau, \tau + \delta]$, contradicting the definition of $\tau$.
    If $p \in \Sigma \setminus \Sigma_0$, then $b_F(p) \neq 0$, say $b_F(p) > 0$; by continuity, $b_F(q(t))$ and $b_F(\widetilde q(t))$ are positive on some interval $[\tau, \tau + \delta]$, so $q_F$ and $\widetilde q_F$ exceed $\bar q$ on $(\tau, \tau + \delta]$, and the same argument applies; the case $b_F(p) < 0$ is the same.
    If $p \in \Sigma_{\mathrm{sl}}$, then by~(b) both solutions stay in $\Sigma_0$ on some interval $[\tau, \tau + \delta]$, and by~(a) both solve the same linear equation there, again a contradiction.
    Hence $p \in \Sigma_\partial$, so any two solutions coincide up to the first time at which one of them, and hence each of them, lies in $\Sigma_\partial$; by~(i) there is at least one solution.
    At a time $t \le T^\circ$ with $q(t) \in \Sigma \setminus \Sigma_0$, $q_F - \bar q$ is strictly monotone near $t$, because $b_F(q(t)) \neq 0$; so such times are isolated, and $q$ crosses $\Sigma$ at each of them.
    Together with~(a), this gives the description of the control.

    \emph{(iii).} Fix an outcome as in the proof of~(i), and a finite $T \le T^\circ$.
    By~(i), every subsequence of $(\bar Q^n, U^n)$ has a further subsequence which converges uniformly on $[0, T]$ to $\big(\widetilde q, \int_0^\cdot \widetilde u\,ds\big)$, where $\widetilde q$ is a solution of~\eqref{eq:BB} with control $\widetilde u$.
    By~(ii), $\widetilde q = q$ on $[0, T]$, and since the control is determined by the state almost everywhere, $\widetilde u = u$ almost everywhere on $[0, T]$.
    Hence the whole sequence converges.

    \emph{(iv).} Suppose that $\lambda > \cM(\bar q)$, $t_0 < T^\circ$ and $q(t_0) \in \Sigma_{\mathrm{sl}}$, and let $t_1 := \sup\{t \ge t_0 : q(s) \in \Sigma_0 \text{ for all } s \in [t_0, t]\}$, which exceeds $t_0$ by~(b).
    On $[t_0, t_1)$, by~(a), $q_I(t)$ lies between $q_I(t_0)$ and $q_I^*(\bar q)$, so $u_\Sigma(q(t))$, which is affine in $q_I$, lies between $u_\Sigma(q(t_0))$ and $\cM(\bar q)/\lambda$, both of which are in $(0,1)$; hence $q(t) \in \Sigma_{\mathrm{sl}}$.
    If $t_1 < \infty$, then $q(t_1) \in \Sigma_0$, as $\Sigma_0$ is closed, and the same bound gives $q(t_1) \in \Sigma_{\mathrm{sl}}$, so that~(b) extends the interval beyond $t_1$, a contradiction.
    Thus $q(t) \in \Sigma_{\mathrm{sl}}$ for all $t \ge t_0$; as $q$ does not lie in $\Sigma_\partial$ before $t_0 < T^\circ$ either, $T^\circ = \infty$, and~(a) gives the formula for $q_I$, with $q_D = \bar\xi - \alpha(\bar q)\,q_I$ and $q_F = \bar q$.
\end{proof}

%% file: appendix-cellb.tex
\section{Linear Program Certificate for Equation~\eqref{eq:cert-b}}\label{app:cellb}

This appendix records the nonnegative-combination certificate for equation~\eqref{eq:cert-b} in the proof of Theorem~\ref{thm:bifurcation-structure}. We use the normalized variables of the Proof of Statement 4 there: $c=\mu_I=1$, $\mu_D=r$, $\bar P=1-P$, $\bar\alpha_0=1-\alpha_0$, with
\[ g_1=\alpha_0^2 r^2+\bar P(1+\alpha_0 r),\quad g_2=\bar\alpha_0 r^2,\quad A=\alpha_0 r+\bar P, \]
and $E=2\alpha_0^2 g_2+b_1 g_1$, $G=-(b_1A+2\alpha_0\beta_0)$, $D_{\mathrm{int}}=\alpha_0 A g_2-\beta_0 g_1$, and $R_A$ as in~\eqref{eq:certA}. The inequality to be certified is
\[ T:=b_0 g_1^2-\alpha_0^2 g_2^2\ \ge\ 0 \qquad\text{on the cell}\quad \{E\ge0,\ G\ge0\}\subseteq(0,1)^2\times(0,\infty), \]
where $E\ge0$, $G\ge0$ encode $Q_v\le\underline{Q}$, $Q_v\ge L$ respectively.

\paragraph{A structural reduction.}
Expanding $EA+G g_1$ gives the identity $EA+Gg_1=2\alpha_0 D_{\mathrm{int}}$; combined with~\eqref{eq:certA} it yields the exact identity
\begin{equation}\label{eq:cellb-clean}
  A\,T \;=\; A\,R_A \;+\; \big(\bar P^2+\bar P r-\alpha_0\bar\alpha_0 r^2\big)\,D_{\mathrm{int}} \;+\; g_1 g_2\,G .
\end{equation}
On the cell $A>0$, $R_A\ge0$ (the coefficients displayed after~\eqref{eq:certA}), $g_1,g_2\ge0$, and $G,D_{\mathrm{int}}\ge0$, so two of the three terms on the right are manifestly nonnegative and the sign of $T$ is reduced to the single bracket $\beta:=\bar P^2+\bar P r-\alpha_0\bar\alpha_0 r^2$. As $\alpha_0\bar\alpha_0\le\tfrac14$, we have $\beta\ge \bar P^2+\bar P r-\tfrac14 r^2\ge0$ whenever $\bar P\ge\tfrac{\sqrt2-1}{2}\,r$, which is all of the cell save a thin sliver on which $\bar P$ is small relative to $r$.

\paragraph{The certificate.}
For a proof valid on the entire cell---in particular on the sliver where~\eqref{eq:cellb-clean} is not manifestly signed---we exhibit the explicit nonnegative-combination certificate of~\eqref{eq:cert-b}. Solving the linear program~\eqref{eq:cert-LP}, with the degree budgets of the proof, over the generators $\{\alpha_0,\,1-\alpha_0,\,P,\,1-P,\,r,\,E,\,G,\,D_{\mathrm{int}}\}$ (the sign-indefinite $\beta_0$ is \emph{not} admitted as a generator) yields positive rational weights $\gamma_j$ and generator products $M_j$, each a monomial in $\alpha_0, 1-\alpha_0, P, 1-P, r$ times one of $1, E, G, D_{\mathrm{int}}$, with
\[ T=\sum_j \gamma_j\,M_j,\qquad \gamma_j>0, \]
an identity in $\mathbb{Q}[\alpha_0,P,r]$. Feasibility was confirmed by both the dual simplex and the interior-point methods of \textsc{HiGHS}; taking the simplex vertex and solving the resulting system exactly over $\mathbb{Q}$ gives a particular solution with $101$ nonzero terms ($66$ in $\sigma_0$, $2$ in $\sigma_E$, $22$ in $\sigma_G$, and $11$ in $\sigma_D$), for which the polynomial identity holds with zero residual under exact expansion. Clearing denominators, $D\,T=\sum_j N_j M_j$ with $D = 82484487281218$ and $N_j$ positive integers; the explicit list of $(N_j,M_j)$ is provided in the ancillary file \texttt{cell\_b\_certificate.txt}. Since each $M_j\ge0$ on the cell and each $N_j>0$, this establishes $T\ge0$ throughout $\{E\ge0,\ G\ge0\}$, closing case~(b).